%% file: 0paper.tex
\documentclass[12pt]{article}
\usepackage{epsfig}
\usepackage{amsfonts}
\title{Divide and Conquer: A Distributed Approach to Five Point Energy
  Minimization}
\author{Richard Evan Schwartz}

\newtheorem{theorem}{Theorem}[section]

\newtheorem{lemma}[theorem]{Lemma}

\newtheorem{corollary}[theorem]{Corollary}

\def\startproof{{\bf {\medskip}{\noindent}Proof: }}

\def\endproof{$\spadesuit$  \newline}

\def\K{\mbox{\boldmath{$K$}}}%
\def\N{\mbox{\boldmath{$N$}}}%
\def\Q{\mbox{\boldmath{$Q$}}}%
\def\R{\mbox{\boldmath{$R$}}}%
\def\Z{\mbox{\boldmath{$Z$}}}%

\DeclareFontFamily{U}{rcjhbltx}{}
\DeclareFontShape{U}{rcjhbltx}{m}{n}{<->rcjhbltx}{}
\DeclareSymbolFont{hebrewletters}{U}{rcjhbltx}{m}{n}
\let\aleph\relax\let\beth\relax
\let\gimel\relax\let\daleth\relax
\DeclareMathSymbol{\aleph}{\mathord}{hebrewletters}{39}
\DeclareMathSymbol{\beth}{\mathord}{hebrewletters}{98}
\DeclareMathSymbol{\gimel}{\mathord}{hebrewletters}{103}
\DeclareMathSymbol{\daleth}{\mathord}{hebrewletters}{100}
\DeclareMathSymbol{\lamed}{\mathord}{hebrewletters}{108}
\DeclareMathSymbol{\mem}{\mathord}{hebrewletters}{109}
\DeclareMathSymbol{\ayin}{\mathord}{hebrewletters}{96}
\DeclareMathSymbol{\tsadi}{\mathord}{hebrewletters}{118}
\DeclareMathSymbol{\qof}{\mathord}{hebrewletters}{114}
\DeclareMathSymbol{\shin}{\mathord}{hebrewletters}{152}

\begin{document}
\maketitle

\begin{abstract}
  We prove the 1977
  Melnyk-Knop-Smith phase transition
  conjecture for $5$-point energy minimization.
  This result contains, as a special case, the solution
  of Thomson's $5$-electron problem from 1904.
  \end{abstract}

  \input{1intro.tex}

  \input{2back.tex}
  \input{3proof.tex}
  \input{4hermite.tex}
  \input{5local.tex}
  \input{6symm.tex}

\input{7endgame.tex}
  \input{8energy.tex}

\input{9energy_proof.tex}
  \input{10calc.tex}
  \input{refs}

\end{document}

\input{refs}
\end{document}

\end{document}

%% file: 1intro.tex
\section{Introduction}

\subsection{History and Context}

Let $S^2$ be the unit sphere in $\R^3$.
Given a configuration $\{\widehat p_0,...,\widehat p_{N-1}\} \subset  S^2$ of $N$ distinct points
and a function $F: (0,2] \to \R$, define
\begin{equation}
  F(P)=\sum_{0\leq i<j <N} F(\|\widehat p_i-\widehat p_j\|).
\end{equation}
This quantity is commonly called the
$F$-{\it potential\/} or the $F$-{\it energy\/} of $P$.
A configuration $P$ is a {\it minimizer\/} for $F$ if
$F(P) \leq F(P')$ for all other $N$-point
configurations $P'$.
The question of finding energy minimizers has a long literature;
the classic case goes back to Thomson [{\bf Th\/}] in 1904.

The classic choice for this question is $F=R_s$, where
  \begin{equation}
    R_s(d)=\frac{{\rm sign\/}(s)}{d^s}.
  \end{equation}
  Here ${\rm sign\/}(s)=1$ if $s>0$ and ${\rm sign\/}(s)=-1$ if $s<0$.
  When $s>0$ this potential is called the {\it Riesz potential\/}. When
  $s<0$ it is called  the {\it Fejes T\'oth potential\/}.
The case $s=1$ is specially called
the {\it Coulomb potential\/} or the
{\it electrostatic potential\/}.
This case of the energy minimization
problem is known as {\it Thomson's problem\/}.
See [{\bf Th\/}].

There is a large literature on the energy 
minimization problem. See [{\bf F\"o\/}] and
[{\bf C\/}] for some early local results.
See [{\bf MKS\/}] for a definitive numerical study
on the minimizers of the Riesz potential for $N$ 
relatively small. The website
[{\bf CCD\/}] has a compilation of experimental
results which stretches all the way up to 
about $N=1000$. 
The paper [{\bf SK\/}] gives a nice survey of results,
with an emphasis on the case when $N$ is large. 
See also [{\bf RSZ\/}].
The paper [{\bf BBCGKS\/}] gives a survey
of results, both theoretical and experimental, about
highly symmetric configurations in higher dimensions.

When $N=2,3$ the problem is fairly trivial.
The cases $N=4,6,12$ are well understood.
The regular tetrahedron, octahedron, and icosahedron are
the unique minimizers for all power law potentials.
See [{\bf KY\/}], [{\bf A\/}], [{\bf Y\/}].
More generally,
[{\bf CK\/}, Theorem 1.2], proves that any so-called
sharp configuration minimizes every Riesz potential
(and a more general class of potentials as well,
called {\it monotone potentials\/}.)

The case $N=5$ has been notoriously intractable.
First let me introduce the two main players.
    The {\it Triangular Bi-Pyramid\/} (TBP) is the $5$-point configuration having
    one point at the north pole, one point at the south pole, and $3$ points
    arranged in an equilateral triangle on the equator.  A {\it Four Pyramid\/}
    (FP) is a $5$-point configuration having one point at the north pole and
    $4$ points arranged in a square equidistant from the north pole.
Here is a run-down on what is known so far:
\begin{itemize}
  \item The paper [{\bf HZ\/}] has a rigorous computer-assisted
proof that the TBP is the unique minimizer
for the potential $R_{-1}$.  That is, the TBP maximizes
the sum of all pairwise distances among all $5$-point configurations.
\item My paper [{\bf S1\/}] has a
rigorous computer-assisted proof that the TBP
is the unique minimizer for the Riesz potentials $R_1$ (Thomson's problem)
and $R_2$.
\item 
The paper [{\bf DLT\/}] gives a traditional
proof that the TBP is the unique minimizer
for the logarithmic potential $F(r)=-\log(r)$.

\item In [{\bf BHS\/}, Theorem 7] it is shown that,
as $s \to \infty$, any sequence of
$5$-point minimizers w.r.t. $R_s$ must converge (up to
rotations) to the FP having one point at the north pole
and the other $4$ points on the equator.
In particular, the TBP is not a minimizer
w.r.t. $R_s$ when $s$ is sufficiently large.

\item Define
  $G_k(r)=(4-r^2)^k$.
  In [{\bf T\/}], A. Tumanov proves that
the TBP is the unique minimizer for $G_2$.
The minimizers for $G_1$ are those configurations
(including the TBP)
whose center of mass is the origin.
This paper was a big inspiration for me.
\end{itemize}

\subsection{The Main Results}

Our main result
verifies the phase-transition for $5$-point
energy minimization first observed in
[{\bf MKS\/}], in $1977$, by
T. W. Melnyk, O. Knop, and W. R. Smith. 
        Define
        \begin{equation}
          15_{+}  = 15 + \frac{25}{512}.
        \end{equation}

        \begin{theorem}[Phase Transition]
          There exists $\shin \in (15,15_{+})$ such that:
  \begin{enumerate}
  \item For $s \in (0,\shin)$ the TBP
    is the unique minimizer for $R_s$.
  \item For $s=\shin$ the TBP and a unique FP are the two minimizers
    for $R_s$.
  \item For each $s \in (\shin,15_{+})$ some FP is the unique minimizer for
    $R_s$.
  \end{enumerate}
\end{theorem}

We also prove the following result.
\begin{theorem}
  \label{EXP}
  For $s \in (-2,0)$ the TBP is the unique
  minimizer for $R_s$.
\end{theorem}

We will concentrate on the Phase Transition Theorem, but in
\S \ref{POLYA} we will explain how to prove Theorem \ref{EXP}
from the same ingredients.

It is also worth mentioning that along the way to proving
the Phase Transition Theorem we will also prove the following
result.

\begin{theorem}
  \label{GG}
  The TBP is the unique minimizer for $G_3,G_4,G_5,G_6$, and
  $$G_5^{\flat}=G_5-25G_1, \hskip 30 pt
  G_{10}^{\sharp\sharp} = G_{10}+28 G_5 + 102 G_2.$$
\end{theorem}
The TBP is not a minimizer for $G_k$ for $k \geq 7$.
The last two potentials are carefully contrived, and look
rather artificial, but it seems worth recording their
properties here.  Their role in the proof of the
Phase Transition Theorem will make more sense later on.
Another potential that is crucial to us is
$G_{10}^{\sharp}=G_{10}+13 G_5 + 68 G_2$.  The TBP is
almost but not quite the
minimizer for $G_{10}^{\sharp}$; see Part 6 of
the next section for a precise statement.

\subsection{Verification}
\label{outline}

To make the proof
easier to verify, I have divided it
up into $7$ self-contained units.
These units could potentially be verified by
$7$ independent readers.
\newline
\newline
{\bf Part 0,  Assembly:\/}
This part of the proof deduces the
main results from results
proved in Parts 1--6.
\newline
\newline
{\bf Part 1,  Interpolation:\/}
We introduce potentials which we
call {\it hybrid triples\/}:
\begin{equation}
  a_0 G_{b_0}(r) + a_1 G_{b_1}(r)+ a_2 G_{b_2}(r), \hskip 20 pt
  a_k \in \Z, \hskip 20 pt G_b(r)=(4-r^2)^b.
\end{equation}
The ones listed in Theorem \ref{GG}, and
also $G_{10}^{\sharp}$, are the ones we use.
This part of the proof
establishes Lemma \ref{LEMMA} in
\S \ref{HermiteMain}, which says that
if the TBP minimizes
various collections of hybrid triples,
then it also minimizes the
Riesz potentials within certain ranges.
For example, Tumanov remarks in
[{\bf T\/}] that if the TBP is a minimizer for
$G_3$ and $G_5$ it is also the unique
minimizer for $R_s$ with $s \in (-2,0) \cup (0,2]$.
The crazy potentials in Theorem \ref{GG} are 
designed to handle larger parameter ranges.
One can view our method here as an {\it ad hoc\/}
variant of Hermite Interpolation, as discussed
e.g. in [{\bf CK\/}].
\newline
\newline
{\bf Part 2, Local Analysis:\/}
In this part of the proof, we
show that there
is an explicitly defined neighborhood $\Omega_0$ of the
TBP in which the TBP minimizes certain
hybrid triples.   See \S \ref{definite}.
Parts 1 and 2 combine
to prove the Phase Transition Theorem
for all (configuration, exponent) pairs
in $\Omega_0 \times (0,15_+]$.
\newline
\newline
{\bf Part 3, Symmetrization:\/}
Let $\K_4$ denote the set of
$5$-point configurations which have
$4$-fold dihedral symmetry.   The
dihedral symmetry group fixes one
of the points.
This part of the proof deals with
a small open subset $\Upsilon$ of configurations
near $\K_4$, and
power law exponents $s \in [12,16]$.
See \S \ref{special}.
Here we produce a retraction
$\Upsilon \to \K_4$
and show that it is energy-decreasing
on $\Upsilon-\K_4$.  This is the deepest and
luckiest part of the paper.
\newline
\newline
{\bf Part 4,  Symmetric Configurations:\/}
This part of the proof treats configurations
in $\Upsilon \cap \K_4$ for $s \in [12,16]$. This is
a $2$-dimensional problem for each power law
parameter $s$.
Our work here combines with Part 3 to prove the
Phase Transition Theorem for all
(configuration, exponent) pairs in
$\Upsilon \times [13,15_+]$.
This is the region where the phase transition
actually occurs.
\newline
\newline
{\bf Part 5, Energy Estimate:\/}
This part of the proof establishes
an estimate which allows us
to prove, just using finitely
many calculations, that an entire
open subset of the configuration space
consists of configurations having
larger $F$-energy than the TBP.
Here $F$ is one of the hybrid triple
potentials of interest to us.
This part of the proof is
completely theoretical.  There
are no computer calculations involved.
The estimate is geometrically motivated
and reasonably sharp.
\newline
\newline
{\bf Part 6, The Big Calculation:\/}
Parts 1,2,3,4 of the proof wipe out
all the pairs (configuration, exponent) in the sets
$\Omega_0 \times (0,15_+]$ and
$\Upsilon \times [13,15_+]$.   The remaining
pairs do not pose a serious threat to the TBP.
In this part of the proof we establish
Theorem \ref{GG} and also the statement
that the TBP is the unique minimizer for
$G_{10}^{\sharp}$ on $\Omega-\Upsilon$.
This calculational result combines with
Part 1 to finish off the proof of
the Phase Transition Theorem and
Theorem \ref{EXP}.   To prove this,
we use the energy
estimate from Part 5 to deal with all
the remaining configurations.  
This part of the proof is the hardest
to verify because it relies on a
massive computer calculation.
On the other hand, the computer calculation
is just doing the same small
thing over and over again.
\newline
\newline
\noindent
{\bf Verification Summary:\/}
As discussed above, each of the hypothetical
$7$ readers only needs to read a relatively
small part of the paper.  Here is a precise
breakdown of what is needed for each part.
\begin{itemize}
\setlength{\itemsep}{1pt}
\setlength{\parskip}{0pt}
\item Reader 0 (assembly): \S 2, \S 3.
\item Reader 1 (interpolation): \S 2.2, \S 3.1, \S 4.
\item Reader 2 (local analysis): \S 2.1,\S 2.2,\S 2.4, \S 3.4, \S 5.
\item Reader 3 (symmetrization): \S 2.1, \S 2.5, \S 2.6, \S 3.5, \S 6.
\item Reader 4 (symmetric configs.): \S 2.1, \S 2.5, \S 2.6, \S 3.6, \S 7.
\item Reader 5 (energy estimate): \S 2.1, \S 8, \S 9.
\item Reader 6 (big calculation): \S 2.1--2.5, \S 3.2, \S 3.6, \S 3.7, \S 8, \S 10.
\end{itemize}

\noindent
{\bf Computer Specifications:\/}
Many of the parts above call auxiliary
computer programs I wrote.  
The computer code is all written in Java and
Mathematica.  The Java code runs on
Temurin Java 21.0.2+13 (LTS).
I ran everything on a 2025 MacBook Pro
running macOS Tahoe 26.5.
The Mathematica code seems to run on all
modern versions of Mathematica.
One can download the computer code from
\newline
\newline
{\bf http://www.math.brown.edu/$\sim$res/Papers/TBP.tar\/}
\newline
\newline
The download has a
main directory with a README
file that explains how the code is
organized. Also, this should be
obvious from the names.
There are $7$ subdirectories, containing the code
relevant for the corresponding part of the paper.
A few of these subdirectories are empty, and
I have included them just to preserve the
structure of the paper.
I use 
 the convention that all computer
files and programs mentioned in
Part $N$ of the paper are contained
in the directory starting with $N$.  Thus,
for instance, the subdirectory
{\bf 6BigCalculation\/} has all the code
for Part 6.
Each subdirectory has further
subdirectories and a README file
that gives additional guidance on
how to access the relevant programs.

In 2026, the reader can feed these programs
into an AI and have the machine explain how
they work and what they do. I tested this claim
extensively.
\newline
\newline
{\bf Deep Versus Shallow Calculations:\/}
Part 6 of our proof relies on a big calculation, made rigorous
by interval arithmetic.
      One important point about the interval arithmetic calculations
      done
      in this paper is that they are {\it shallow\/}.
      To put a finer point on this observation, let me contrast
      it with what might happen with iterative floating point calculations.
      In a setting like dynamical systems, interval arithmetic might
      be needed to control the floating point errors that could
      accumulate
      during a long iteration.  There is some error in the output of
      the
      calculation, and then this feeds back into the calculation, and
      so on.  In such a setting, involving deep and recursive
      calculations,
      interval arithmetic is a very important mechanism for controlling
      round-off error and thereby deriving
      rigorous conclusions.

      In contrast, the calculations here are
      just doing the same small thing over and over again, with
      each calculation independent from the all the others.
      While our calculation does involve a recursive
      subdivision tree, the subdivision part is done with
      exact integers.  In the end, our
      calculation is a million miles wide and one millimeter deep.
      There
      is no accumulation of error.  The interval arithmetic here
      mostly serves as a kind of slowly incanted blessing over
      the floating point calculation.

\subsection{AI Usage}

I proved all the results in this paper and
wrote all the accompanying software without
using AI tools.  The original version of this
paper goes all the way back to 2016, though
since then I have made a number of
changes and improvements. (The first version
was about $200$ pages long!)  I used
AI for the first time
in the final
pass through the paper in Sept, 2026.
I used  ChatGPT Sol 5.6,
on ``high intelligence'', for a number of
mostly routine (but quite extensive) tasks. 

In Sept, 2026, ChatGPT and I went through my
paper, section by section, to verify that it was
correct.  ChatGPT helped implement the suggestions made by the referee,
found additional typos and notational inconsistencies,
checked the logic of the proofs, and
went through the computer code and checked that it
really functions as it is supposed to.   Getting
ChatGPT to say that the paper is mathematically
clean and consistent is a pretty exacting task.

There were a few  important things that ChatGPT
helped me with:
\begin{enumerate}
\item In \S 5, ChatGPT helped me optimize the constants.
  This was just a lot of tedious arithmetic.
\item  In \S 7, ChatGPT added some nitty gritty details
    in Lemma C4.  Previously I
    had just said you could eyeball the result.
  \item In \S 9, ChatGPT found 
many small notational problems and a few glitches in the
formulas.  Getting all the notation here really took a
collaboration between myself and ChatGPT.
  \item For the code related to \S 10,
    ChatGPT fixed an inefficiency in my
    interval arithmetic implementation, where I had
    been repeatedly computing some bounds that
    could be hardcoded.  This small change sped up
    the calculations from \S 10 by a factor of about $6$.  This was
    great.
\end{enumerate}

\subsection{Acknowledgements}

I would like to thank Alejandro Adem, Peter Doyle,
Doug Hardin, Ed Saff,
Javier G\'omez-Serrano, and Stephen D. Miller
for their helpful comments and encouragement.
I would also like to thank the
Institut des Hautes \'{E}tudes Scientifiques
for their support, in the form of the Gretchen and
Barry Mazur Chair.  It was a pleasure to finish
up this long project at I.H.E.S. in fall 2026.

\newpage

%% file: 2back.tex
\section{Preliminaries}
\label{prelim}

   \subsection{Avatars}
   \label{avatar}

    Let $S^2 \subset \R^3$ be the unit $2$-sphere.
    {\it Stereographic projection\/} is the map $\Sigma: S^2 \to \R^2 \cup \infty$ given
    by the following formula.
\begin{equation}
\label{stereoBasic}
\Sigma(x,y,z)=\bigg(\frac{x}{1-z}, \frac{y}{1-z}\bigg).
\end{equation}
Here is the inverse map:
\begin{equation}
\label{inversestereo}
\Sigma^{-1}(x,y)=\bigg(\frac{2x}{1+x^2+y^2},\frac{2y}{1+x^2+y^2},
1-\frac{2}{1+x^2+y^2}\bigg).
\end{equation}
$\Sigma^{-1}$ maps circles in $\R^2$ to circles
in $S^2$ and $\Sigma^{-1}(\infty)=(0,0,1)$.
\newline

Stereographic projection gives us a correspondence
between $5$-point configurations on $S^2$ having
$(0,0,1)$ as the last point and planar configurations:
\begin{equation}
\widehat p_0,\widehat p_1,\widehat p_2,\widehat p_3,(0,0,1) \in S^2 \hskip 6 pt
\Longleftrightarrow \hskip 6 pt
p_0,p_1,p_2,p_3 \in \R^2, \hskip 20 pt
\widehat p_k=\Sigma^{-1}(p_k).
\end{equation}
We call the planar configuration the
{\it avatar\/} of the corresponding
configuration in $S^2$.
We call $2$ avatars {\it isomorphic\/} if
the corresponding $5$-point configurations on $S^2$ are isometric.
We mean this in the widest possible sense: The isometry taking
one configuration to the other need not preserve the labelings
of the points.

We write $F(p_0,p_1,p_2,p_3)$
when we mean the $F$-potential of the corresponding
$5$-point configuration.  If
$\xi=(p_0,p_1,p_2,p_3)$ then we will write
\begin{equation}
  F(\xi)=F(p_0,p_1,p_2,p_3).
\end{equation}
Isomorphic
avatars have the same $F$-potential.

A calculation shows that
\begin{equation}
  \|p_k\| \leq 1/2 \hskip 30 pt
  \Longleftrightarrow \hskip 30 pt \|\widehat p_k - (0,0,1)\| \geq
  \frac{4}{\sqrt 5}.
\end{equation}
We call an avatar {\it odd\/} or {\it even\/} according as
it has an odd or even number of points $p_k$ with
$\|p_k\| \leq 1/2$.  The cutoff of $1/2$ is not really
canonical; it is just a convenient choice.

 Figure 2.1 shows two avatars of the triangular bi-pyramid,
 first separately and then
superimposed.  We call the one on the left the {\it even avatar\/}, and
the one in the middle the {\it odd avatar\/}. 
Let $\xi_0$ denote the even avatar.
The points of $\xi_0$ are $(\pm 1,0)$ and $(0,\pm \sqrt
3/3)$.

\begin{center}
\resizebox{!}{1.7in}{\includegraphics{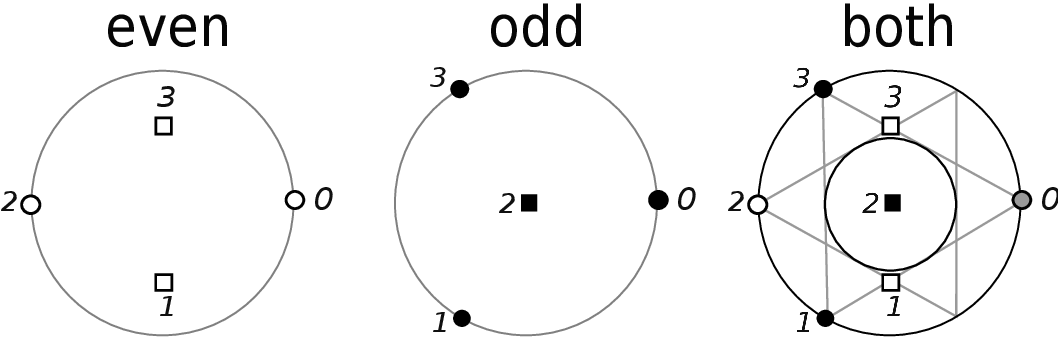}}
\newline
    {\bf Figure 2.1:\/} Even and odd avatars of the TBP.
  \end{center}

When we superimpose the two avatars we see some extra geometric
structure that is not relevant for our proof but worth mentioning.
The two circles
respectively have radii $1/2$ and $1$ and
the $6$ segments shown are tangent to the inner one.

\begin{lemma}
  \label{EVEN}
  Every avatar is isomorphic to an even avatar.
\end{lemma}

\startproof
We form a graph by joining
two points in a $5$-point configuration by an edge
if and only if their distance is at least $4/\sqrt 5$.
Every graph with an odd number of vertices
has a vertex of even degree. Hence,
our graph has a vertex of even degree.
We rotate
so that such a vertex is at $(0,0,1)$.  The
corresponding avatar is even.
 \endproof

 \subsection{The Hybrid Triples}
 \label{hybrid}

Now we introduce the potentials which we
use in order to understand the Riesz potentials.
Define
\begin{equation}
  \label{function0}
  G_k(r)=(4-r^2)^k.
\end{equation}
Also define
$$
G_5^{\flat} = G_5-25 G_1,$$
$$G_{10}^{\sharp\sharp} = G_{10}+28 G_5 + 102 G_2,$$
\begin{equation}
  \label{function1}
G_{10}^{\sharp} =  G_{10}+13 G_5 + 68 G_2 
\end{equation}

I found these hybrid triples experimentally.
They look rather arbitrary, but in fact they
are close to the unique choices which
function the right way in our proof.

 \subsection{The Big Domain}
 \label{biggg}

 We consider an avatar $\xi=(p_0,p_1,p_2,p_3)$.  We denote the
 $(x,y)$-coordinates as
 \begin{equation}
   p_k=(p_{k1},p_{k2}).
 \end{equation}
 We use this notation throughout the paper.
Suitably rotating about the origin,
we always normalize so that
$p_{02}=0$.
We define a domain $$\Omega \subset \R^8=(\R^2)^4$$ to be the set
of avatars $\xi$ satisfying the
following conditions.

\begin{enumerate}
\item $\xi$ is even.
  \item
    $\max(\|p_1\|,\|p_2\|,\|p_3\|) \leq \|p_0\|$ and
    $\max(\|p_1\|,\|p_2\|,\|p_3\|) < 3/2$.
\item $p_{12} \leq p_{22} \leq p_{32}$ and $p_{22} \geq 0$.
\item $p_{01} \in (0,2]$ and $p_{02}=0$.
\end{enumerate}
In Item 2 the second inequality is strict.
The Containment Theorem, in
\S \ref{containthm}, says that for
all the potentials we consider,
every minimizing avatar is isomorphic 
to one in $\Omega$.

\subsection{A Neighborhood of the TBP}
\label{definite}

Let $\xi_0$ denote the even avatar for the TBP.
When we string out the points of $\xi_0$, we get
\begin{equation}
  (1,\underline 0,0,-u,-1,0,0,u), \hskip 30 pt
  u=\frac{1}{\sqrt 3}.
\end{equation}
Let $\Omega_0$ denote the cube of side-length $2^{-16}$
centered at $\xi_0$.   When it comes time
to analyze the picture inside $\Omega_0$
we will drop the second (underlined) coordinate, because
it is always zero for our avatars. This will
let us consider $\Omega_0$ as a cube in
$\R^7$ rather than in $\R^8$.

\subsection{The Special Domain}
\label{special}

We let $\Upsilon \subset \R^8$
denote the set of avatars $(p_0,p_1,p_2,p_3)$ such that
\begin{enumerate}
\item $\|p_0\| \geq \|p_k\|$ for $k=1,2,3$. (This is consistent with
  the Big Domain.)
\item $512 p_0 \in [433,498] \times [0,0]$. (That is, $p_0 \in
  [433/512,498/512] \times \{0\}.$)
\item $512 p_1 \in [-16,16] \times [-464,-349]$.
\item $512 p_2 \in [-498,-400] \times [0,24]$.
\item $512 p_3 \in [-16,16] \times [349,464]$.
\end{enumerate}
It turns out that  $\Upsilon$ contains the avatars that compete with the TBP near
the phase transition exponent $\shin$.

\begin{center}
\resizebox{!}{3.5in}{\includegraphics{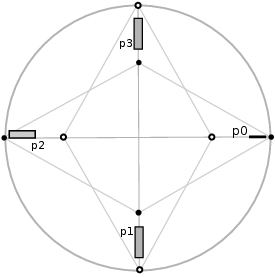}}
\newline
{\bf Figure 2.2:\/} The sets defining $\Upsilon$ compared with two TBP avatars.
\end{center}

For our computational algorithm in \S \ref{calcthm}, it
is useful to work with the symmetrized version $\Upsilon^*$ of
$\Upsilon$ at our disposal as well. This is the same
domain except that
$512 p_2 \in [-498,-400] \times [-24,24]$.
We have
\begin{equation}
  \label{symmU}
  \Upsilon^* \cap \Omega = \Upsilon \cap \Omega,
\end{equation}
so for our purposes the two domains are equivalent.

\subsection{Polynomials and Exponential Sums}

\label{oper}

\subsubsection{Positive Dominance}
\label{WPD}

See [{\bf S2\/}] and [{\bf S3\/}] for more details about the material
in this section.     Let $P \in \R[x_1,...,x_n]$ be a multivariable polynomial:
\begin{equation}
  P=\sum_{I} c_{I} X^{I}, \hskip 20 pt
  X^{I}=\prod_{i=1}^n x_i^{I_i}.
  \end{equation}
Given two multi-indices $I$ and $J$, we write
$I \preceq J$ if $I_i \leq J_i$ for all $i$.
Define
\begin{equation}
  \label{WPD0}
  P_{J}=\sum_{I \preceq J} c_{I}, \hskip 30 pt
  P_{\infty}=\sum_I c_{I}.
\end{equation}
We call $P$ {\it positive dominant\/} if $P_J>0$ for all $J$.

\begin{lemma}
  If $P$ is positive dominant, then $P>0$ on $[0,1]^n$.
\end{lemma}

\startproof
Suppose first that $n=1$. Write
$$
P(x)=a_0+a_1x+\cdots+a_mx^m,
\qquad
A_j=a_0+\cdots+a_j>0.
$$
The following identity makes it clear that $P(x)>0$ for $x \in [0,1]$.
$$P(x)=(1-x)(A_0+A_1 x +... + A_{m-1}x^{m-1}) + A_m x^m.$$
For the general case, let
$$
P=f_0+f_1x_n+\cdots+f_mx_n^m,
\qquad
f_j\in\R[x_1,\ldots,x_{n-1}].
$$
Each polynomial
$P_j=f_0+\cdots+f_j $
is positive dominant.
By induction on $n$, we have
$P_j>0$ on $[0,1]^{n-1}$.
Now fix $x_1,\ldots,x_{n-1}\in[0,1]$.
The resulting polynomial $Q(x_n)$ has successive partial sums
$A_j=P_j(x_1,...,x_{n-1})>0$ for $j=0,...,m$.
Hence $Q>0$ on $[0,1]$ by the $n=1$
case.
\endproof

We call $P$ {\it weak positive
  dominant\/} if $P_J \geq 0$ for all $J$ and $P_{\infty}>0$.

\begin{lemma}
  If $P$ is weak positive dominant then $P>0$ on $(0,1]^n$.
\end{lemma}

\startproof
The same argument as in the preceding lemma shows that
$P \geq 0$ on $[0,1]^n$ if $P$ is weak positive dominant.
(We don't even need $P_{\infty}>0$ for this.)
Choose $N$ larger than every exponent appearing in $P$ and
also so that $P_{\infty}>1/N$.   Then
$Q=P-(1/N) (x_1...x_n)^N$ is again weak
positive dominant.
Hence $Q \geq 0$ on $[0,1]^n$.  Hence $P>0$ on $(0,1]^n$.
\endproof

  \subsubsection{Polynomial Subdivision}
\label{subdivision}
  
      Let $P \in \R[x_1,...,x_n]$.
      For any $x_j$ and $k \in \{0,1\}$ we define
    \begin{equation}
      S_{{x_j},k}(P)(x_1,...,x_n)=P(x_1,...,x_{j-1},x_j^*,x_{j+1},...,x_n),
      \hskip 15 pt
      x_j^*=\frac{k}{2}+\frac{x_j}{2}.
     \end{equation}
      If $S_{{x_j},k}(P)>0$ on $(0,1]^n$ for $k=0,1$ then we also have
        $P>0$ on $(0,1]^n$.
    
        \subsubsection{Numerator Selection}
        \label{numplus}
        
        If $f=f_1/f_2$ is a bounded rational function on $[0,1]^n$,
              written so that $f_1,f_2$ have no common factors,
              we always choose $f_2$ so that $f_2(1,...,1)>0$.
              If we then show, one way or another,
              that $f_1>0$ on $(0,1]^n$ we can conclude that $f_2>0$ on
                $(0,1]^n$ as well.   The point is that $f_2$ cannot
                vanish, and hence cannot change
                  sign, because then $f$ blows up.
                  But then
                  $f>0$ on $(0,1]^n$.   We write ${\rm num\/}_+(f)=f_1$.

                  \subsubsection{Exponential Sums}
                 \label{ESUM}

    \begin{lemma}[Convexity]
      Suppose that $\alpha,\beta,\gamma \geq 0$ have the property that $\alpha+\beta \geq 2\gamma$.  Then
      $\alpha^s+\beta^s  \geq 2\gamma^s$ for all $s>1$, with equality iff $\alpha=\beta=\gamma$.
    \end{lemma}

    \startproof
By strict convexity of the function $f(x)=x^s$ we have
    $$\frac{\alpha^s+\beta^s}{2} \geq
    \bigg(\frac{\alpha+\beta}{2}\bigg)^s \geq \gamma^s.$$
    We get equality only when $\alpha=\beta=\gamma$.
    \endproof

    Let us recall Descartes' Rule of Signs.  Suppose
    $$P(t)=\sum_{i=1}^n c_i t^{e_i},$$
    where $e_1>...>e_n \in \N$.   The
    coefficients $c_1,...,c_n$ are all nonzero.
    Let $K$ be
    the number of sign changes in the sequence
    $c_1,...,c_n$ of coefficients.  Then
    $P$ has at most $K$ zeros in $(0,\infty)$.
    We use this result to prove a related
    fact about exponential sums.
        
    \begin{lemma}[Descartes]
    Let $0<r_1 \leq ... \leq r_n<1$ be a sequence of
    positive numbers.
    Let $c_1,...,c_n$ be a sequence of nonzero numbers.
    Define
    \begin{equation}
      E(s)= \sum_{i=1}^n c_i\ r_i^s.
    \end{equation}
    Let $K$ denote the number of sign changes in the
    sequence $c_1,...,c_n$ of coefficients.
    Then $E$ changes sign at most $K$ times on $\R$.
    \end{lemma}

    \startproof
    Suppose we have a counterexample. If we perturb the bases
    slightly, we still have a counterexample. So, we can perturb so
    that
    our counter-example has the form
       $r_i = r^{e_i}$ for some $r \in (0,1)$ and $e_1>...>e_n \in
    \N$.     Set $t=r^s$.    As $s$ ranges in $\R$, the variable $t$
    ranges in $(0,\infty)$.
    We are now considering the sign changes of the polynomial
    $$E(s)=P(t)= \sum_{i=1}^n c_i\ t^{e_i}.$$
    But $P(t)$ changes sign at most $K$ times on $(0,\infty)$ by Descartes' Rule of Signs.
    This gives us a contradiction.
    \endproof

\newpage

%% file: 3proof.tex
\section{Proof Assembly}

\noindent
{\bf Reading Guide:\/}
This chapter is for Reader 0.

\subsection{Interpolation}
\label{interp00}

We use the notation from \S \ref{prelim}.
Here is our main result about interpolation.

\begin{theorem}[Interpolation]
  Let $T_0$ be the TBP.
  Then
  \begin{enumerate}
     \item Suppose    $s \in (0,13]$ and $T$ is any
     $5$-point configuration. If we have
    $F(T_0)<F(T)$ for all 
    $F=G_4,G_5,G_6,G_{10}^{\sharp \sharp}$
    then
    ${R_s}(T_0)<{R_s}(T)$.
     \item Suppose    $s \in [13,15_+]$ and $T$ is any
     $5$-point configuration. If we have
    $F(T_0)<F(T)$ for all 
    $F=G_5^{\flat},G_{10}^{\sharp }$
    then
        ${R_s}(T_0)<{R_s}(T)$.
  \end{enumerate}
  \end{theorem}

  \subsection{The Containment Theorem}
  \label{containthm}

  Let $\xi_0$ be the even avatar for the TBP.   Let
  $\xi=(p_0,p_1,p_2,p_3)$ be another avatar.

\begin{theorem}[Containment]
  Let $F$ be any of
 $$G_4,\ \ G_5^{\flat},\ \ G_5,\ \ G_6,\ \ G_{10}^{\sharp},\ \ G_{10}^{\sharp \sharp}.$$
  If $F(\xi) \leq F(\xi_0)$ then $\xi$ is isomorphic to an
    avatar in $\Omega$.
  \end{theorem}

  \startproof
  Since the TBP has $6$
bonds of length $\sqrt 2$, and $3$ of length $\sqrt 3$, and
$1$ of length $\sqrt 4$, we have
$G_k(\xi_0)=6 \times 2^{k}+3.$
We use this formula to compute $F(\xi_0)$ for
the relevant choices of $F$.

Let $\xi=(p_0,p_1,p_2,p_3)$ be some other avatar.
For $j=0,1,2,3$, let
$\tau_j$ be the term in $F$ corresponding
to the pair $(p_j,p_4)$. That is
\begin{equation}
  \tau_j=F(\|\Sigma^{-1}(p_j)-(0,0,1)\|).
\end{equation}
For every potential of interest to us, the quantity
$\tau_j$ is monotone increasing in $\|p_j\|$ if
$\|p_j\| \geq 3/2$.   We use this fact to
compute $\tau_j$ for the relevant function
choices.  Here is a table of important
values for our proof.

\begin{center}
\begin{tabular}{c|c|c|c}
$F$ & $F(\xi_0)$ &
$\|p_j\|>3/2$ &
$\|p_j\|>2$ \\ \hline
$G_4$ & $99$ &
$\tau_j>58$ &
$\tau_j>104$ \\[3pt]
$G_5$ & $195$ &
$\tau_j>162$ &
$\tau_j>335$ \\[3pt]
$G_5^\flat+30$ & $120$ &
$\tau_j>123$ &
--- \\[3pt]
$G_6$ & $387$ &
$\tau_j>450$ &
--- \\[3pt]
$G_{10}^{\sharp}$ & $10518$ &
$\tau_j>26909$ &
--- \\[3pt]
$G_{10}^{\sharp\sharp}$ & $14361$ &
$\tau_j>31863$ &
---
\end{tabular}
\end{center}

From this chart we deduce
the following
{\it norm constraints\/}:  We have
$F(\xi_0)<F(\xi)$ if either 
$\|p_0\|>2$ or $\|p_0\|, \|p_j\| \geq 3/2$ for some $j=1,2,3$.

Again, we assume $F(\xi) \leq F(\xi_0)$.
By Lemma \ref{EVEN}, we can take $\xi$ to
be even.  
We re-order the points so that
$\|p_j\| \leq \|p_0\|$ for $j=1,2,3$.
We rotate about the origin, to arrange that $p_0$ lies in
the positive $x$-axis.   We reflect in the $X$-axis if
necessary to guarantee that at least two of the
points $p_1,p_2,p_3$ lie in the closed upper half-plane.
Finally, we re-order 
so that $p_{12} \leq p_{22} \leq p_{32}$.
Note that $p_{22} \geq 0$.
Given our norm constraints we must have
$p_{01} \in (0,2]$
and $\|p_j\|<3/2$ for $j=1,2,3$.
Hence $\xi \in \Omega$.
\endproof

  \begin{corollary}
    \label{exclude}
    Let $T_0$ be the TBP and
    let $T$ be a configuration that has no avatar isomorphic to one in
    $\Omega$.  Then
    $R_s(T_0)<{R_s}(T)$ for all $s \in (0,15_+]$.
  \end{corollary}

  \startproof
  This is an immediate corollary of the Interpolation Theorem and the
  Containment Theorem.
  \endproof
  
\subsection{Local Analysis}

Recall that $\Sigma^{-1}$ is inverse stereographic projection.
In this section we work in the small cube $\Omega_0$ defined
in \S \ref{definite}.
Here is the main local result we prove.

\begin{theorem}[Local Convexity]
  \label{hess}
    For $F=G_4,G_6,G_5^{\flat},G_{10}^{\sharp}$, the
  Hessian of $F \circ \Sigma^{-1}$ is positive
  definite at every point of $\Omega_0$.
\end{theorem}

\begin{corollary}
  \label{loc}
  Let $F$ be any of
 $$G_4,\ \ G_5^{\flat},\ \ G_5,\ \ G_6,\ \ G_{10}^{\sharp},\ \ G_{10}^{\sharp \sharp}.$$
  Then $\xi_0$, the TBP avatar, is the unique $F$-energy minimizer inside $\Omega_0$.
\end{corollary}

\startproof
Let $F$ be any of the functions from the Local Convexity Theorem. 
Let $\xi\in \Omega_0$ be other than $\xi_0$.
The Local Convexity Theorem combines with the vanishing gradient
to show
the restriction of $F \circ \Sigma^{-1}$ to the line segment $\gamma$
joining $\xi_0$ to $\xi$ is strictly convex and has $0$ derivative at $\xi_0$.
Hence we have  $F(\xi)>F(\xi_0)$.

It remains to deal with $G_5$ and $G_{10}^{\sharp\sharp}$.
We get the same result for $G_5$ because $G_5$ is a positive
combination of $G_5^{\flat}$ and $G_1$, and the
TBP, like any configuration whose center of mass is the
origin, is a global minimizer for $G_1$.
Likewise $G_{10}^{\sharp\sharp}$ is a positive
combination of $G_{10}^{\sharp}$ and $G_2$ and $G_5$,
and the TBP is a minimizer for $G_2$,
by Tumanov's result [{\bf T\/}].
\endproof

Combining Corollary \ref{loc} with
the Interpolation Theorem, we get:

\begin{corollary}
  \label{loc2}
  The
Phase Transition Theorem is true 
for $\Omega_0 \times (0,15_+]$.
In this region, there is no phase transition: The TBP is the
unique minimizer.
\end{corollary}

\subsection{Symmetrization}
\label{SYMM}

In this section we deal directly with the
Riesz potentials.  We deal with the configurations
in the special domain $\Upsilon$ from \S \ref{special}.
Let $\K_4$ denote the set of avatars
which are invariant under reflections in the
coordinate axes.
We describe a symmetrization
operation which maps $\Upsilon$ into $\K_4$.
Let $(p_0,p_1,p_2,p_3)$ be an
avatar with $p_0 \not = p_2$.
Define
{\small
\begin{equation}
  -p_2^*=p_0^*=(x,0), \hskip 8 pt
  - p_1^*=p_3^*=(0,y), \hskip 12 pt
  x=\frac{\|p_0-p_2\|}{2}, \hskip 8 pt
  y=\frac{\|\pi_{02}(p_1-p_3)\|}{2}.
\end{equation}
\/}
Here $\pi_{02}$ is the projection onto the
subspace perpendicular to
$p_0-p_2$. 
The avatar $(p_0^*,p_1^*,p_2^*,p_3^*)$ lies
in $\K_4$.
Note that our operation fixes avatars in $\K_4$.

\begin{theorem}[Symmetrization]
  \label{five}
  Let $(p_0,p_1,p_2,p_3) \in \Upsilon-\K_4$.
  Then $$R_s(p_0^*,p_1^*,p_2^*,p_3^*) <
  R_s(p_0,p_1,p_2,p_3)$$ when $s \geq 12$.
\end{theorem}

\subsection{Symmetric Configurations}
\label{CRIT}

Let $(x,y)$ denote the avatar of the form
\begin{equation}
  \label{PSI4}
      (x,0), \hskip 20 pt (0,-y), \hskip 20 pt (-x,0), \hskip 20 pt
      (0,y).
    \end{equation}
    We define the following subsets of such avatars.
$$
\Psi_4=\{(x,y):64(x,y)\in[43,64]^2\},
\hskip 30 pt
\Psi_4^\sharp=\{(x,y):64(x,y)\in[55,56]^2\}.
$$
We have $\Upsilon \cap\K_4\subset\Psi_4$.  Also
$\Psi_4^{\sharp}$ is a small subset of $\Psi_4$,
the spot where it turns out that all the action happens.

    Let $\Psi_8 \subset \Psi_4$ and $\Psi_8^{\sharp} \subset
    \Psi_4^{\sharp}$
    denote the diagonals, where $x=y$.
    
    Define the following map from $\Psi_4^{\sharp}$ into $\Psi_8$.
    \begin{equation}
  \label{ENDGAME}
  \sigma(x,y) = (z,z), \hskip 30 pt
  z=\frac{x+y+(x-y)^2}{2}.
\end{equation}

Note that $\sigma$ is the identity on
$\Psi_8$.  Here are the main results.

\begin{theorem}[Critical I]
 $R_s(\sigma(p)) < R_s(p)$ for $$(p,s) \in 
 (\Psi_4^{\sharp}-\Psi_8^{\sharp}) \times [14,16].$$
\end{theorem}

\begin{theorem}[Critical II]
  $R_s(\xi_0)<R_s(\xi)$ for 
    $$(\xi,s) \in (\Psi_4 \times [13,15]) \hskip 10 pt \cup \hskip 10 pt
    ((\Psi_4-\Psi_4^{\sharp}) \times [15,15_+]).$$
\end{theorem}

\begin{theorem}[Critical III]
  There exists $\shin \in (15,15_{+})$ such that
  \begin{enumerate}
  \item   $R_s(\xi_0)<R_s(\xi)$ for all $(\xi,s) \in
    \Psi_8^{\sharp} \times [15,\shin)$.
  \item    $R_s(\xi_0)>R_s(\xi)$ for some (fixed) $\xi \in
    \Psi_8^{\sharp}$ and all $s \in (\shin,15_{+}]$.
\item $R_s$ is uniquely minimized on $\Psi_8^{\sharp}$ for all $s \in (15,15_+]$.
  \end{enumerate}
\end{theorem}

\noindent
{\bf Remark:\/}
We call these theorems Critical, but in a certain sense, they are
just the endgame. Even without
these results we could deduce from everything else that a
minimizer for $R_s$ must be the TBP  for
each parameter $s \in (-2,0) \cup (0,13]$ and either the TBP or
an FP for each $s \in (13,15_+]$.
\newline

\begin{corollary}
  \label{symm2}
  The Phase Transition Theorem is true for $\Upsilon \times [13,15_+]$.
\end{corollary}

\startproof
Suppose that $s \in [13,15_+]$ and
$\xi \in \Upsilon$ is a minimizer.
We have $R_s(\xi) \leq R_s(\xi_0)$.
By the Symmetrization Theorem,
$\xi \in \Upsilon \cap \K_4 \subset \Psi_4$. 
The Critical Theorem II tells us that
$s \in (15,15_+]$ and $\xi \in \Psi_4^{\sharp}$.
The Critical Theorem I now says that
$\xi \in \Psi_8^{\sharp}$.
Statement 3 of the Critical Theorem III
says that $\xi$ is unique.

If $s \in [15,\shin)$ then
Statement 1 of the Critical Theorem III gives
a contradiction, and there is no
$\xi \in \Psi_8^{\sharp}$ which is a minimizer.
If $s \in (\shin,15_+]$ then Statements 2 and 3
of the Critical Theorem III say that
$\xi$ is the unique $R_s$-minimizer
in $\Psi_8^{\sharp}$.

Now consider the case that $s=\shin$.  It follows
from continuity that $\xi_0$ is a minimizer
for $R_s$ and also that the minimum value of
$R_s$ on $\Psi_8^{\sharp}$ equals $R_s(\xi_0)$.
Hence there is some unique $\xi \in \Psi_8^{\sharp}$
such that $R_s(\xi)=R_s(\xi_0)$.
\endproof

\subsection{Big Calculation}
\label{UPS}

\begin{theorem}[Calculation]
  The following is true.
  \begin{enumerate}
    \item
  The TBP is the unique minimizer for
  $G_4,G_5^{\flat},G_6$ amongst $5$-point configurations which have avatars
in $\Omega-\Omega_0$.
\item The TBP is the unique minimizer
  for $G_{10}^{\sharp}$ among $5$-point configurations which have
  avatars in
$\Omega-\Omega_0-\Upsilon$.
\item  The TBP has lower $G_{10}^{\sharp \sharp}$-energy
  than all of the configurations having avatars in $\Upsilon$.
      \end{enumerate}
\end{theorem}

\begin{corollary}
  \label{big}
  The TBP is the unique minimizer for
  $$G_4,\ \ G_5^{\flat},\ \ G_5,\ \ G_6, \ \ G_{10}^{\sharp \sharp}$$
  among configurations
  having
  avatars in $\Omega$.
  The TBP is the unique minimizer for
  $G_{10}^{\sharp}$ amongst configurations having avatars in
  $\Omega-\Upsilon$.
\end{corollary}

\startproof
Corollary \ref{loc} takes care of all the statements on $\Omega_0$.
The same argument as in
Corollary \ref{loc} deduces the result about $G_5$
from the result about $G_5^{\flat}$.
The Calculation Theorem and Corollary \ref{loc}
deal with
$G_{10}^{\sharp}$ on
$\Omega-\Upsilon$.
The same argument as in Corollary \ref{loc} now shows
that the TBP is the unique minimizer for $G_{10}^{\sharp \sharp}$
on $\Omega-\Upsilon$.  Statement 3 
finishes the job for $G_{10}^{\sharp\sharp}$.
\endproof

\noindent
{\bf Proof of the Phase Transition Theorem:\/}
Corollaries \ref{loc2} and \ref{symm2} reduce us to proving
the Phase Transition Theorem on
  \begin{equation}
    \label{BIGG}
    (\Omega-\Omega_0) \times (0,13]  \hskip 10 pt \cup \hskip 10 pt
    (\Omega-\Omega_0 - \Upsilon) \times [13,15_+].
  \end{equation}
Combining Corollary \ref{big} with the
Interpolation Theorem, we see that the
Phase Transition Theorem is true
on the domain in Equation \ref{BIGG}.
\endproof

\subsection{Proof of Theorem \ref{EXP}}
\label{POLYA}

Here we explain the changes in the argument above
needed to prove Theorem \ref{EXP}.
We replace the Interpolation Theorem
by the following result.

\begin{theorem}[Special Interpolation]
  Let $T_0$ be the TBP.
  Suppose that $s \in (-2,0)$ and $T$ is any
     $5$-point configuration. If we have
    $F(T_0)<F(T)$ for 
    $F=G_3,G_5$ then
        ${R_s}(T_0)<{R_s}(T)$.
      \end{theorem}
      This result was already known to
      Tumanov [{\bf T\/}].  The proof, which
      we omit, is essentially the same as our
      proof of Case 1 of Lemma A in
      \S \ref{A2proof}.  See \S \ref{discuss}.

      Theorem \ref{EXP} follows from Theorem
      \ref{GG} and the Special Interpolation Theorem.
      The only part of Theorem \ref{GG} we have not
      proved is the statement for $G_3$.
      As in the proof of the Containment Theorem,
      we have the following chart:
      \begin{center}
\begin{tabular}{c|c|c|c|c}
$F$ & $F(\xi_0)$ &
$\|p_j\|>3/2$ &
$\|p_j\|>4$ &
$\beta$ \\ \hline
$G_3$ & $51$ &
$\tau_j>21$ &
$\tau_j>53$ &
$ \geq 128/9$
\end{tabular}
\end{center}

Here $\beta$ denotes the $G_3$ energy for a $4$-point
configuration of $(\widehat p_0,\widehat p_1,\widehat p_2,\widehat
p_3)$.
This quantity is minimized for the regular tetrahedron.
Hence we have $G_3(\xi_0)<G_3(\xi)$ 
if either
$\|p_0\|>4$ or $\|p_0\|, \|p_j\|>3/2$ for some $j=1,2,3$.
In view of these norm constraints, the same proof as in the
Containment Theorem establishes the same result for
$G_3$, except for the larger domain $\Omega^*$ in which the
      condition $p_{01} \in (0,4]$ replaces the
      condition $p_{01} \in (0,2]$.  Running a version of
      the Calculation Theorem on $\Omega^*$ we get
      the following result.

      \begin{theorem}[Calculation II]
       The  TBP is the unique minimizer for $G_3$ amongst $5$-point
        configurations having avatars  in $\Omega^*-\Omega_0$.
      \end{theorem}

      Referring to \S \ref{local}, the conclusions of
      Lemmas L1 -- L3 also hold for the scaled
      potential $3G_3$. So, the analysis in \S \ref{local} finishes
      the job by establishing
      Corollary \ref{loc} for $G_3$.  Hence the TBP is the
unique global minimizer for $G_3$.
      
\newpage

%% file: 4hermite.tex
\section{The Interpolation Theorem}
\label{A2proof}

\noindent
{\bf Reading Guide:\/}
This chapter is for Reader 1.

\subsection{Main Result}
\label{HermiteMain}

Recall that $15_+=15+\frac{25}{512}$.
We let $R_s(T)$ be the Riesz $s$-potential of
a configuration $T$.
Referring to Equations \ref{function0} and \ref{function1}, we define
\begin{equation}
  P_1=(G_4,G_6), \hskip 30 pt
  P_2=(G_5,G_{10}^{\sharp \sharp}), \hskip 30 pt
  P_3=(G_5^{\flat},G_{10}^{\sharp}),
\end{equation}
\begin{equation}
  I_1=(0,6], \hskip 30 pt
  I_2=[6,13], \hskip 30 pt
  I_3=[13,15_+].
\end{equation}

We say that a pair $(\Gamma_3,\Gamma_4)$ of functions
{\it forces\/} the interval $I$ if the following is
true:  If $T$ is another $5$-point configuration such that
$\Gamma_k(T_0)<\Gamma_k(T)$ for $k=3,4$ then
$R_s(T_0)<R_s(T)$ for all
$s \in I$.

In this chapter we prove the following result, which immediately implies
the Interpolation Theorem from \S \ref{interp00} 
\begin{lemma}[A]
  \label{LEMMA}
  The following is true.
  \begin{enumerate}
  \item The pair $(G_4,G_6)$ forces $(0,6]$.
  \item The pair $(G_5,G_{10}^{\sharp\sharp})$ forces $[6,13]$.
  \item The pair $(G_5^{\flat},G_{10}^{\sharp})$ forces $[13,15_{+}]$.
  \end{enumerate}
\end{lemma}

\subsection{Reduction to Smaller Results}

It is convenient to define
\begin{equation}
  \Gamma_k=G_k, \hskip 30 pt k=0,1,2.
\end{equation}
Note that $\Gamma_0$ is just the constant function.
We say that a pair of functions $(\Gamma_3,\Gamma_4)$
{\it specially forces\/} $s>0$ if there is some
$a_0 \in \R$ and $a_1,a_2,a_3,a_4>0$ such that
\begin{equation}
  \label{COMBO0}
  \Lambda(x):=\sum_{i=0}^4 a_i\Gamma_i(x) \leq R_s(x), \hskip 30 pt
  \forall x \in [0,2],
\end{equation}
with equality at $x=\sqrt 2,\sqrt 3,\sqrt 4$.
What this means is that on $[0,2]$, the graph of $\Lambda$ lies under the
graph of $R_s$ but the two graphs touch at these three special
values.  We prefer writing $\sqrt 4$ instead of $2$ so that the
three values are more obviously related.

We say that $(\Gamma_3,\Gamma_4)$ {\it specially forces\/} the
interval $I$ if this pair specially forces all $s \in I$.
In this case, the constants will likely depend on $s$.

\begin{lemma}[A1]
  If $(\Gamma_3,\Gamma_4)$ specially
  forces $I$ then $(\Gamma_3,\Gamma_4)$ forces $I$.
  \end{lemma}

  \startproof
  Let $T_0$ be the TBP and let $T$ be some other $5$-point
  configuration.  We simplify the notation and write
  $F(T)={\cal E\/}_{F}(T)$.
  We assume
  \begin{equation}
    \label{force}
    \Gamma_j(T_0)<\Gamma_j(T), \hskip 30 pt j=3,4,
    \end{equation}
  and we want to show 
  that $R_s(T_0)<R_s(T)$ for all $s \in I$.

  Fix some $s \in I$.
  We claim that
  \begin{equation}
    a_j(s)\Gamma_j(T_0)  \leq a_j(s)\Gamma_j(T), \hskip 30 pt j=0,1,2,3,4.
  \end{equation}
  For $j=0$ this is obvious; we have equality.
   For $j=1$ this follows from the
  well-known fact that the minimizers of ${\cal E\/}_{\Gamma_1}$ are
  precisely those configurations whose center of mass is the origin.
  For $j=2$ this follows immediately from Tumanov's result [{\bf
    T\/}].
  For $j=3,4$ this follows from Equation \ref{force}.
  For the cases $j=1,2,3,4$ we are using the fact that $a_j(s)>0$.
  We also note that we have strict inequality when $j=3,4$.

  Let $\Lambda_s$ be the sum in
  Equation \ref{COMBO0} with respect to $a_j(s)$.
The quantities $\sqrt 2, \sqrt 3, \sqrt 4$ are the distances
which appear between pairs of points in $T_0$.  Therefore
$\Lambda_s(T_0)=R_s(T_0)$.  But then
$$R_s(T) \geq \Lambda_s(T)=
\sum_{j=0}^4 a_j(s)\Gamma_j(T) >
\sum_{j=0}^4 a_j(s) \Gamma_j(T_0)=
\Lambda_s(T_0)=R_s(T_0).$$
This completes the proof.
\endproof

\begin{lemma}[A2]
  For each $i=1,2,3$ the pair $P_i$ specially forces   $I_i$.
\end{lemma}

Lemma A is an immediate consequence of Lemma A1 and Lemma A2.
It remains to prove Lemma A2.

\subsection{Discussion}
\label{discuss}

Let me explain what made me search for these results and
how I found them.
The approach I take is related to linear-programming-type techniques
used in many works on energy minimization, such as [{\bf CK\/}].
Tumanov remarks in
[{\bf T\/}] (using somewhat different language)
that the pair $(G_3,G_5)$ forces
the parameter intervals $(-2,0)$ and $(0,2]$.  I found my
own proof of this. Subsequently, I
saw that $(G_4,G_6)$ specially forces $(0,6]$.
This is Case $i=1$ of Lemma A2.

The TBP is not the global
minimizer for $G_7,G_8,...$  If we want to
use expressions like $G_7$, etc, we need to
average it with $G_k$ for smaller values of
$k$ to give the TBP a chance of being the
minimizer.  I experimented with expressions
$a_1 G_{b_1}+a_2 G_{b_2}$ and wasn't having
much luck.  So, I then broadened the
search to the hybrid triples.  This worked
quite well.

My computer programs 
allow the reader to specify a hybrid triple, solve
for the coefficients, and then
check visually whether the properties for special forcing hold.
After fooling around for a while I hit on the specific
expressions that appear in Lemma A2.  Once I got
the expressions, it was a matter of computer algebra
to prove that the plots on my computer program
are indeed an accurate reflection of mathematical reality.

The proof of Lemma A2  relies on
interval arithmetic calculations in Java.
I give explicit expressions for everything
(with computer plots), so a really energetic
reader could find their own ways to verify
that the plots are accurate reflections of
mathematical reality.

\subsection{Techniques of Proof}
      \label{EO}

In our proofs below, we will need to deal with expressions of the form
\begin{equation}
  \label{summand}
  F(s)=\sum c_i s^{t_i} b_i^{s/2},
  \end{equation}
  where $b_i,c_i \in \Q$ and $t_i \in \N$ and
  $b_i>0$.
  \newline
  \newline
  \noindent
  {\bf Sandwiching Rationals:\/}
  For the term $b_i^{s/2}$ we use floating point point operations
  and then round outward to produce {\it sandwiching rational
  numbers\/} $x_{i0}, x_{i1}$ such that
  \begin{equation}
    \label{sand}
    x_{i0} < b_i^{s/2}< x_{i1}
  \end{equation}
We  verify Equation \ref{sand} in a way we explain next.
\newline
\newline
{\bf Expanding Out Method:\/}
Suppose we want to establish an inequality
like
$(\frac{a}{b})^{\frac{p}{q}}<\frac{c}{d}$,
where every number involved is a positive integer.
This inequality is true iff
$b^pc^q-a^pd^q>0.$
We check this using exact integer arithmetic.
The same idea works with $(>)$
in place of $(<)$.  We use this method
to verify Equation \ref{sand}.
\newline
\newline
{\bf Checking Positivity:\/}
To check the positivity of $F$ on
some interval $[s_0,s_1]$ we produce, for each term,
the $4$ sandwiching rationals $x_{i00}, x_{i10}, x_{i01},x_{i11}$.
Here $i$ is the index enumerating the terms of $F$, and
$x_{ijk}$ is the $j$th sandwiching rational computed with
respect to the $i$th term and the value $s_k$.  We then let
\begin{equation}
  y_i=\min_{j,k,\ell \in \{0,1\}} c_i\ s_{\ell}^{t_i}\ x_{ijk}.
\end{equation}
Then the value of the summand on $[s_0,s_1]$ is
greater than or equal to $y_i$.  This means that
$\sum y_i$ is
a lower bound for Equation \ref{summand} for
all $s \in [s_0,s_1]$.  On any  interval $I$  where
we want to show that Equation \ref{summand} is positive,
we pick the smallest dyadic interval $[0,2^m]$ that contains $I$
and then run the following subdivision algorithm.

\begin{enumerate}
\item Start with an ordered list $L$ of intervals. Initially $L=\{[0,2^m]\}$.
\item If $L$ is empty, then {\bf HALT\/}.  Otherwise 
  $L=L_0,...,L_n$.  We set $Q=L_n$.
  \item If either $Q \cap I = \emptyset$ or the method above
    shows that Equation \ref{summand} is positive on $Q$
  we delete $Q$ from $L$ and go to Step 2.
\item Otherwise we modify our list to
  be $L_0,...L_{n-1},Q_1,Q_2$, where $Q_1$ and $Q_2$ are the
  left and right halves of $Q$.
      Then we go  to Step 2.
\end{enumerate}
If this algorithm halts then
it constitutes a proof that $F(s)>0$ for all $s \in I$.
\newline

Here is another tool we will use below in the proof.
This next result is a repeat of a result
in \S \ref{WPD}.

\begin{lemma}[Positive Dominance]
  A real polynomial $a_0+a_1 t + ...  + a_n t^n$ is
  positive on $[0,1]$ provided that $A_i>0$ for
  $i=0,...,n$.  Here $A_i=a_0+...+a_i$.
    \end{lemma}

  \startproof
  Call the polynomial $P$.
  The following identity makes the result clear.
    $$P(x)=(1-x)(A_0+A_1 x +... + A_{n-1}x^{n-1}) + A_n x^n.$$
  This completes the proof.
  \endproof
  
  \subsection{Reduction of Lemma A2}

Referring to Equation \ref{COMBO0}
we solve the equations
\begin{equation}
  \Lambda_s(\sqrt m) = R_s(\sqrt m),  \hskip 8 pt m=2,3,4, \hskip 18 pt
  \Lambda'_s(\sqrt m) = R_s'(\sqrt m), \hskip 8 pt m=2,3.
\end{equation}
Here $f'$ denotes the derivative of $f$, a function defined on
$(0,2]$.  We don't need to constrain $f'(2)$.
For each $s$ this gives us a linear system with $5$ variables and
$5$ equations.  In all cases, our solutions have the following structure
\begin{equation}
  \label{MATX}
  (a_0,a_1,a_2,a_3,a_4)=M(
2^{-s/2}, 3^{-s/2}, 4^{-s/2} ,
s2^{-s/2}, s3^{-s/2})
\end{equation}
Here $M$ is a $5 \times 5$ matrix which we list out
in each of the $3$ cases.

\begin{lemma}[A21]
  For each $i=1,2,3$ the following is true.
  When $M$ is defined relative to the pair $P_i$ then
  the coefficients $a_1,a_2,a_3,a_4$ are positive functions on
  the interval $I_i$.
\end{lemma}

We want to see that the function
  \begin{equation}
    H_s=1-\frac{\Lambda_s}{R_s}.
  \end{equation}
  has interior local minima at $r=\sqrt 2, \sqrt 3$ on $(0,2]$.
  Letting $H_s'=dH_s/dr$, etc.,
\begin{equation}
\label{deriv}
H'_s(r)=  -r^{s-1} \psi_s, \hskip 30 pt \psi_s:=s\Lambda_s(r)+r\Lambda_s'(r).
\end{equation}
Using the general equation
$rG_k'(r)=2kG_k(r)-8kG_{k-1}(r)$,
we see that $\psi_s$
is a polynomial in $t=4-r^2$.

\begin{lemma}[A22]
  For each choice $P_j$ and each $s \in I_j$ the following is true.
  The function $\psi_s$ has $4$ simple roots in $[0,4]$.  Two of
  the roots are $1$ and $2$ and the other two respectively
  lie in $(0,1)$ and $(1,2)$.
\end{lemma}

Let us deduce Lemma A2.  Our construction and Lemma A21 immediately
take care
of everything in the special forcing definition except for the
inequality
that is equivalent to the statement that $H_s \geq 0$ on $(0,2]$.
We now justify this.  The roots of $\psi_s$ in $[0,4)$
  are in bijection
  with the roots of $H'_s$ in $(0,2]$ and simple roots are mapped
  to simple roots.
  In Cases 1,2,3 respectively we check that at the parameters
  $s=3,10,14$ the roots of $H_s'$ at $r=\sqrt 2$ and $r=\sqrt 3$
  correspond to local minima and the
  other two roots of $H_s'$ correspond to local maxima.  See the plots on the
  right side of Figures 4.1,4.2,4.3 below.
Since the roots
remain simple for all $s$ in the relevant interval, the nature of the roots cannot change
as $s$ varies. Hence $H_s$ has exactly $2$ interior local minima in $(0,2]$, at
$r=\sqrt 2, \sqrt 3$.
Furthermore $\lim_{r \to 0_+} H_s(r)=1$.
But then $H_s \geq 0$ on $(0,2]$. This completes the proof.

  \subsection{Data and Plots}

  Referring to Equation \ref{MATX}, we list out the matrix $M$ in
  each of the three cases and also show computer plots.
  The reader can interact with these plots and see others
  like (and unlike) them using our computer software.

  Here is Case 1.
  {\small
  $$
M=\frac{1}{792}\left[\matrix{
0 & 0 & 792 & 0 & 0  \cr
792 & 1152 & -1944 & -54 & -288  \cr 
-1254 & -96 & 1350 & 87 & 376 \cr
528 & -312 & -216 & -39 & -98  \cr
-66 & 48 & 18 & 6 & 10}\right]
=
\left[\matrix{
1 & 2 & 4 & 16 & 64 \cr
1 & 1 & 1 & 1 & 1 \cr
1 & 0 & 0 & 0 & 0 \cr
0 & 4 & 16 & 128 & 768 \cr
0 & 6 & 12 & 24 & 36
}\right]^{-1}.
$$
\/}
We list the inverse just because it is quite simple.
The left side of Figure 4.1 shows a graph of
$80a_1, 200 a_2, 2000 a_3, 10000 a_4$, 
considered as functions of the exponent $s$.
Here $a_1,a_2,a_3,a_4$ are colored darkest to lightest.
The completely unimportant positive multipliers
are present so that we get a nice picture.
On the left side of Figure 4.1, the thick
vertical segments are $s=0,1,2,3,4,5,6$.
 It turns out that $a_3$ goes
negative between $6$ and $6.1$, so the interval
$(0,6]$ is fairly nearly maximal for positivity.

The right side of Figure 4.1 plots the graph of
$H_3$ for $r \in [4/3,2]$.  Here
$\sqrt 2$ and $\sqrt 3$ are local minima rather
than local maxima or saddle points.

\begin{center}
  \resizebox{!}{2.3in}{\includegraphics{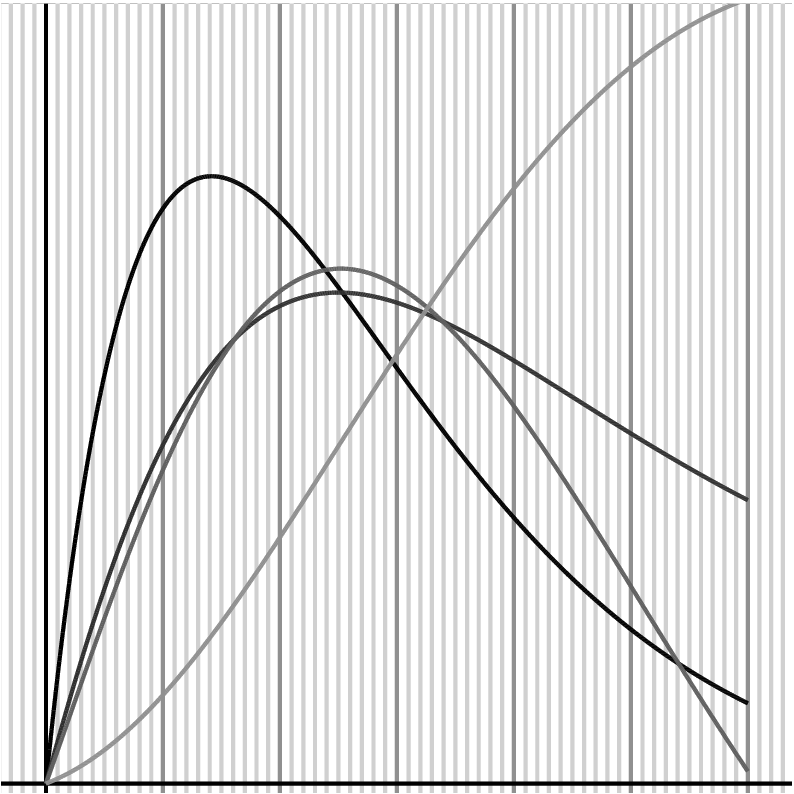}}
  \hspace{0.2in}
\resizebox{!}{2.3in}{\includegraphics{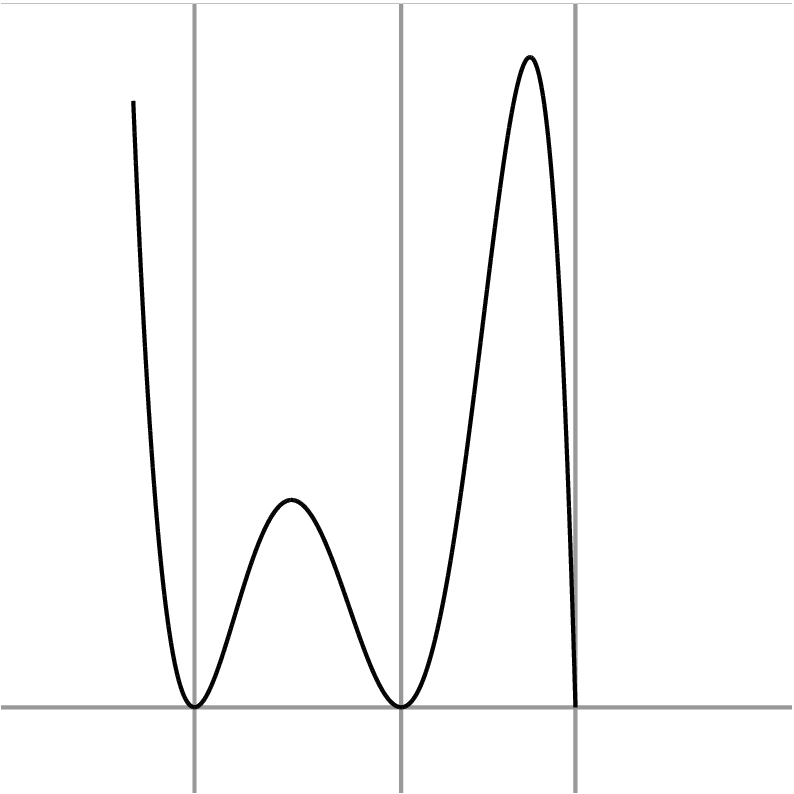}}
\newline
{\bf Figure 4.1:\/} Plots for Case 1.
\end{center}

We cannot directly apply our positivity algorithm
to Case 1 because this algorithm only works
for functions which have uniform positive
lower bounds. We will deal with Case 1 below.

Here is Case 2.

  \begin{equation}
M=\frac{1}{268536}\left[\matrix{
0& 0& 268536& 0& 0 \cr 
88440& 503040& -591480& -4254& -65728 \cr 
-77586& -249648& 327234& 2361& 65896\cr 
41808& -19440& -22368& -2430& -9076\cr 
-402& 264& 138& 33& 68}\right]
\end{equation}

Figure 4.2 does for Case 2 what Figure 4.1 does for Case 1.
This time the left hand side plots
$500 a_1,500a_2,5000a_3,500000 a_4$
for $s \in [6,13]$.
The thick vertical segments are
$s=6,7,8,9,10,11,12,13$.
The coefficients $a_1,a_2,a_3$ go negative for $s$ just 
a tiny bit larger than $13$.
The right hand side shows a plot of $H_{10}$
for $r \in [5/4,2]$.

\begin{center}
  \resizebox{!}{2in}{\includegraphics{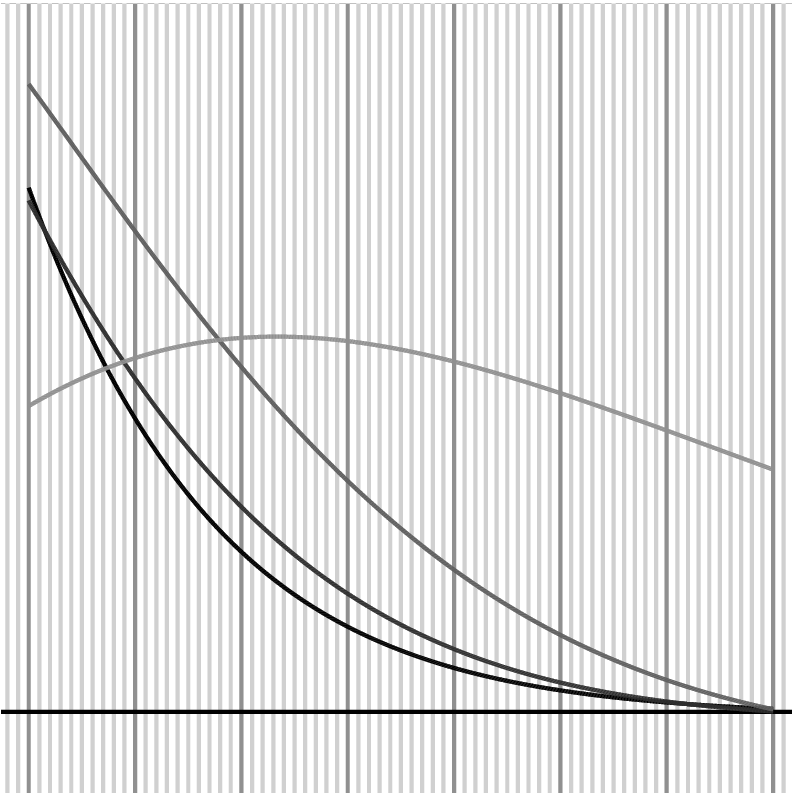}}
  \hspace{0.25in}
\resizebox{!}{2in}{\includegraphics{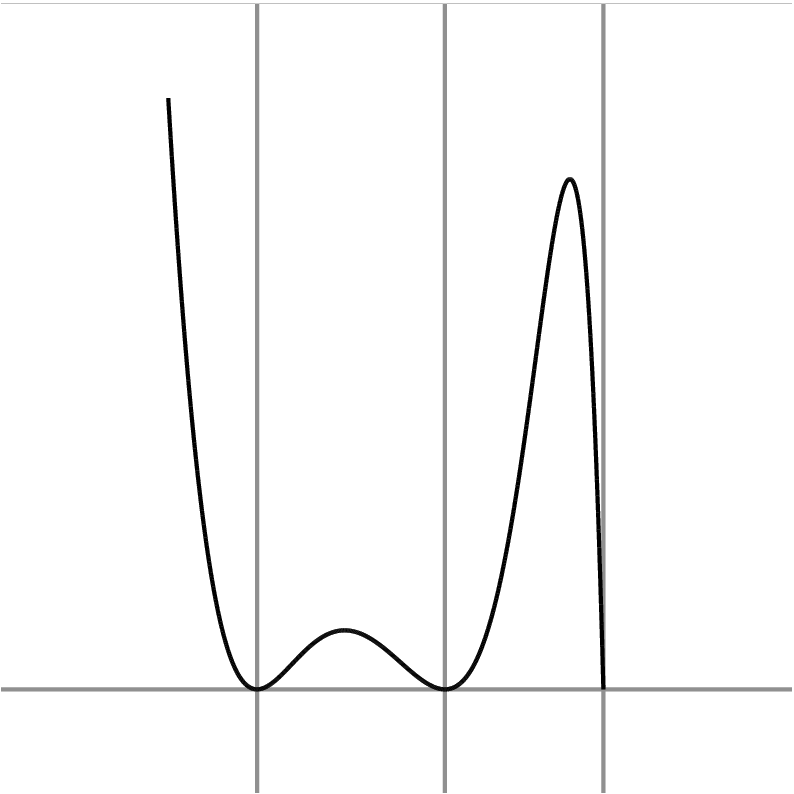}}
\newline
{\bf Figure 4.2:\/} Plots for Case 2.
\end{center}

For Case 2 we run the positivity algorithm
and show that for $k=1,2,3,4$ the
function $a_k(s)$ is positive on $[6,13]$,
as the plot indicates.

Here is Case 3.

 \begin{equation}
   \label{CASE4}
  M= \frac{1}{268536}\
\left[\matrix{
0&0&268536&0&0 \cr
982890&116040&-1098930&-52629&-267128  \cr
-91254&-240672&331926&3483&68208  \cr
35778&-15480&-20298&-1935&-8056  \cr
-402&264&138&33&68}\right]
\end{equation}

This matrix is quite similar to the one in the previous
case, because we are essentially still taking combinations
of $G_0,G_1,G_2,G_5,G_{10}$.  We are just grouping the
functions differently.
Figure 4.3 does for Case 3 what Figure 4.2 does for Case 2.
This time we plot $500a_1,15000a_2,20000a_3,500000a_4$
for $s \in [13,16]$.
The thick vertical segments are $s=13,14,15$.

The coefficients $a_1,a_2,a_3$ go negative for $s$ just 
a tiny bit larger than $15.05$.   In particular,
everything up to and including our cutoff of
$15+25/512$ is covered.
The right hand side shows a plot of $H_{14}$ for
$r \in [5/4,2]$.

\begin{center}
  \resizebox{!}{2in}{\includegraphics{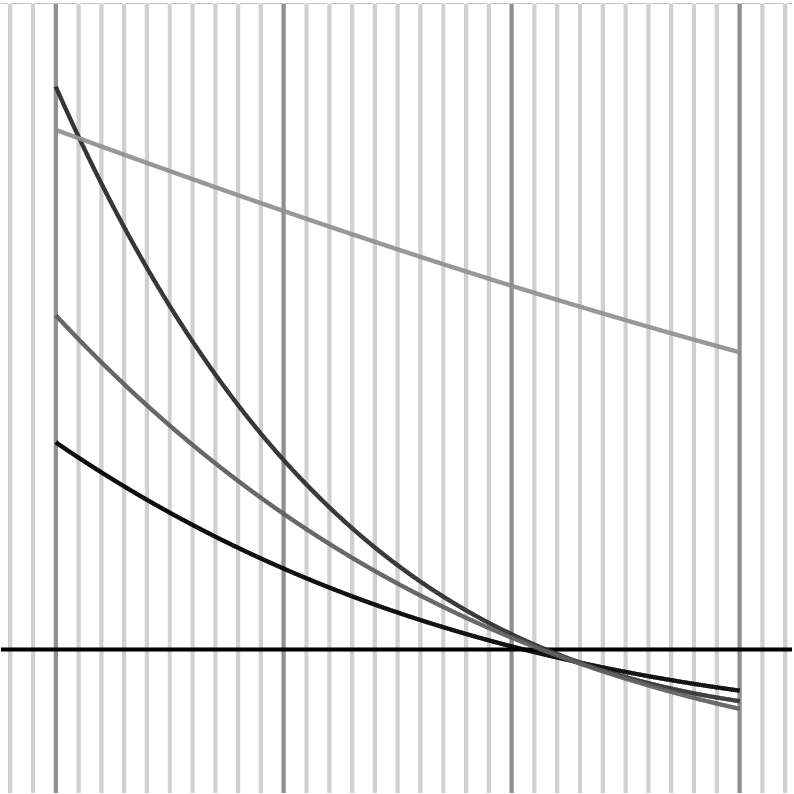}}
  \hspace{0.25in}
\resizebox{!}{2in}{\includegraphics{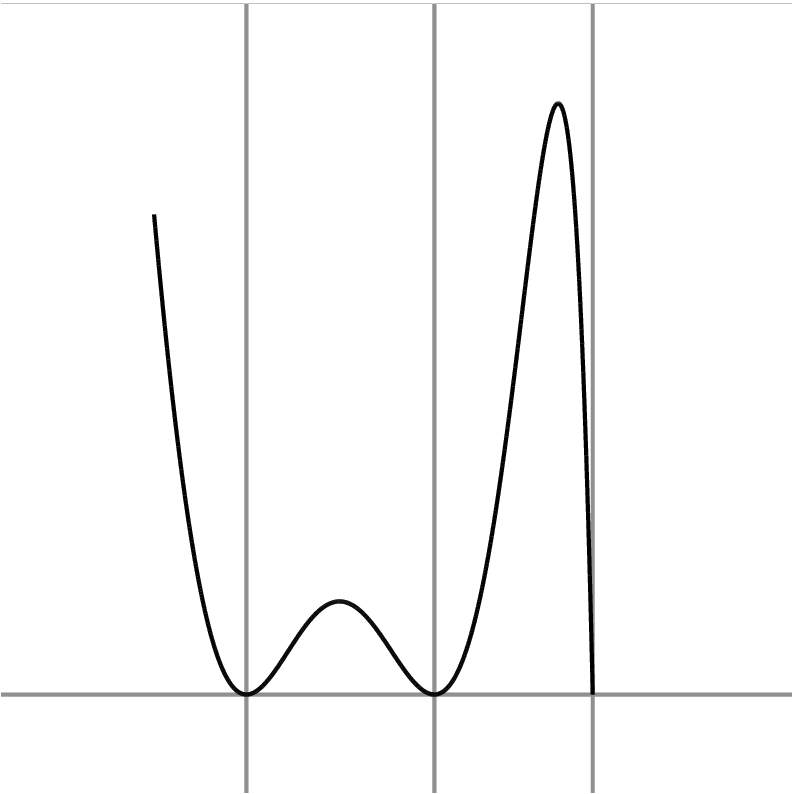}}
\newline
{\bf Figure 4.3:\/} Plots for Case 3.
\end{center}

For Case 3 we run the positivity algorithm
and show that for $k=1,2,3,4$ the
function $a_k(s)$ is positive on $[13,15_+]$,
as the plot indicates.

\subsection{Case 1 of Lemma A21}
\label{Polya2}

For Lemma A21, only Case 1 remains. We first add
two additional quantities we want to test, namely
   $\psi_s(0)$ and $\psi_s(4)$.   We have
  {\footnotesize
$$11\psi_s(0)=\left[\matrix{-88 \cr -128 \cr +216 \cr +6 \cr +32\cr +11}\right] \cdot
\left[\matrix{2^{-s/2} \cr 3^{-s/2} \cr 4^{-s/2} \cr s2^{-s/2} \cr s3^{-s/2} \cr s4^{-s/2}}\right],
\hskip 15 pt
\frac{11}{s} \psi_s(4)=\left[\matrix{-2112\cr +1664\cr +459\cr +219
    \cr 288 \cr 0}\right] \cdot
  \left[\matrix{2^{-s/2} \cr 3^{-s/2} \cr 4^{-s/2} \cr s2^{-s/2} \cr
      s3^{-s/2} \cr s 4^{-s/2}}\right]$$
  \/}
  These quantities have essentially the same form as
  the functions $a_j(s)$ for $j=1,2,3,4$.  The only
  difference is that they also involve the extra term $s 4^{-s/2}$.
  We run the positivity algorithm and show that all $6$ quantities
  are positive on $[1/4,6]$.
  
  Now we deal with $(0,1/4]$.
  All our (scaled) expressions have the form
  $$Y \cdot V(s), \hskip 30 pt V(s)=(2^{-s/2},3^{-s/2},4^{-s/2},s 2^{-s/2}, s 3^{-s/2},s 4^{-s/2}),$$
  for some integer vector $Y$.
Setting $\mu_s$ to be any of the $6$ component functions of $V(s)$, we
have
  \begin{equation}
  \sup_{s \in [0,1]} \bigg|\frac{\partial^6 \mu_s}{\partial s^6}
  \bigg|<1.
  \end{equation}
  Moreover the sum of the absolute values of the coefficients in each of the $Y$ vectors
  is at most $5000$.   This means that, when we take the $5$th order Taylor series
  expansion for $Y \cdot V(s)$, the error term is at most
  $$5000 \times 1 \times \frac{1}{6!} \times s^6<7s^6.$$
  We compute each Taylor series, set all non-leading positive terms to $0$, and
  crudely round down the other terms:
  $$792a_1(s): \hskip 30 pt 98 s - 69 s^2 +0 s^3 -6 s^4 + 0 s^5 - 7 s^6$$
  $$792a_2(s): \hskip 30 pt 14s - 3s^2 -2 s^3 + 0 s^4 -1 s^5 - 7 s^6.$$
  $$792a_3(s): \hskip 30 pt 1s + 0 s^2 -1 s^3 + 0 s^4 + 0 s^5 -7 s^6.$$
  $$792a_4(s): \hskip 30 pt.03 s + 0 s^2 + 0 s^3 - .01 s^4 + 0 s^5 -7 s^6.$$
  $$11\psi_s(0): \hskip 30 pt .08 s + 0 s^2 - .02 s^3 + 0 s^4 - .01 s^5 - 7 s^6.$$
  $$(11/s)\psi_s(4): \hskip 30 pt 11+ 0s + 0 s^2 -1 s^3 - 1 s^4 + 0 s^5 - 7 s^6.$$
  These under-approximations are all easily seen to be positive on $(0,1/4]$.

    \subsection{Proof of Lemma A22}
    \label{proofA22}

    \noindent
    {\bf Case 1:\/}
    We use the analysis surrounding Equation \ref{deriv}.
    The two expressions
    $\psi_s$ and $H_s'$ vanish at exactly the
    same values of $r \in (0,2)$.  We also noted
    that $\psi_s$ is a polynomial in $t=4-r^2$.
    We write $\psi_s(r)=f_s(t)$.
We compute that
  \begin{equation}
\label{monic1}
f_s(t)= t^6-\frac{48}{12+s} t^5 + ...
  \end{equation}
  We don't care about the other terms.
  Let $N$ be the number of roots of $H'_s$ for $r \in (0,2)$,
  counting multiplicity.  This is the same as 
  the number of roots of $f_s$ for $t = 4-r^2 \in (0,4)$,
  counting multiplicity.
  
  Since $f_s$ has degree $6$ we conclude that
  $N \leq 6$.
  By construction we have $H_s(\sqrt m)=H_s'(\sqrt m)=0$ for $m=2,3$ and
  $H_s(\sqrt 4)=0$.
  Since $H_s(\sqrt 2)=H_s(\sqrt 3)$ we must have
  $H_s'(r_{23})=0$ for some $r_{23} \in (\sqrt 2, \sqrt 3)$.
  Since $H_s(\sqrt 3)=H_s(\sqrt 4)$, we must have
  $H_s'(r_{34})=0$ for some
  $r_{34} \in (\sqrt 3,\sqrt 4)$.
  Now we know that $H_s'$ has roots at
  $r=\sqrt 2, r_{23}, \sqrt 3, r_{34}$.
  Correspondingly $f_s$ has roots $t_1=1$ and $t_2=2$
  and $t_{01} \in (0,1)$ and $t_{12} \in (1,2)$.
  The sum of all the roots
  of $f_s$ is $48/(12+s)<4$.
  Since $t_1+t_2+t_{01}+t_{12}>4$ we see that not all
  roots can be positive.  Hence $N<6$.
  Since $f_s$ is positive at $t=0,4$ we see that
  $N$ is even.  Hence $N=4$.  This means that
  the only roots of $f_s$ in $(0,4)$ are the $4$ roots
  we already know about.  Since these roots are distinct, they are
  simple roots.
  Hence the only roots of $H_s'$ in $(0,2)$ are the ones
  already mentioned, and they are all simple roots.
  \newline
  \newline
      {\bf Cases 2 and 3:\/}
      First of all, the functions $H_s$ are the same in Cases 2 and 3.
      This is not just a computational accident.  In both cases we
      are building $H_s$ from the functions $G_1,G_2,G_5,G_{10}$.
      So, we combine Cases 2 and 3 by proving that the common
      polynomial $\psi_s$ just has $4$ roots for each $s \in [6,16]$.
    I will describe a proof which took me quite a lot of
    experimentation to find. 
    
    The same analysis as in Case 1 shows that $\psi_s$ has
    roots at $1,2$, and in $(0,1)$ and in $(1,2)$.
    We just want to see that there are no other roots.

We can factor $\psi_s$ as $(t-1)(t-2)\beta_s$
where $\beta_s$ is a degree $8$ polynomial. Taking
derivatives with respect to $t$, we notice that
\begin{enumerate}
\item $\gamma_s=268536 \times 12^{s/2} \times (\beta_s''-\beta_s')$
  is positive for $(s,t) \in [6,16] \times [0,4]$.
\item $-\beta_s'(0)>0$ for all $s \in [6,16]$.
\item $\beta'_s(4)>0$ for all $s \in [6,16]$.
\end{enumerate}
In Statement 1, the $t$-dependence is only in the arguments
of $\beta_s'$ and $\beta_s''$.
Statement 1 shows in particular that $\beta_s'$ never has a double root.
This combines with Statements 2 and 3 to show that the number of roots
of $\beta_s'$ in $[0,4]$ is independent of
$s \in [6,16]$.  We check explicitly that $\beta_6'$ has only one
root in $[0,4]$. Hence $\beta'_s$ always has just one root.
But this means that $\beta_s$ has at
most $2$ roots in $[0,4]$.  This, in turn, means that
$\psi_s$ has at most $4$ roots in $[0,4]$.
Since we already know four distinct roots,
these are all the roots.

We claim that the two roots
$a,b$ of \(\beta_s\) are simple.
The reason is as follows.
\begin{itemize}
  \item By Rolle's Theorem,
    $\beta'_s$ has a root in $(a,b)$.
  \item We also know $\beta'_s$ has
    exactly one root in $[0,4]$.
  \end{itemize}
  Hence $\beta_s'$ does not vanish
  at $a$ or $b$.
  Since $\beta_s(1) \not =0$ and
  $\beta_s(2) \not =0$, 
the other two roots $t=1,2$ of
$\psi_s$ are simple.
Hence all four roots of \(\psi_s\) are simple.
This completes the proof modulo the $3$ statements.

Now we establish the $3$ statements.
We first give a formula for $\gamma_s$.  Define
matrices $M_3, M_4, M_6$ respectively as:

{\footnotesize
$$\left[\matrix{-546840& -1800480& 99720& -397440& -234600& -33120& 173880& -22080 \cr
    18366& 17112& 80766& 24288& 18630& 11592& 4830& -1104  \cr
    0& 0& 0& 0& 0& 0& 0& 0}\right]$$
  \/}
{\footnotesize
  $$
\left[
\begin{array}{cccccccc}
 -345600 & -1576320 & -509760 & -760320 & -448800 & -63360 & 332640 & -42240 \\
 -199296 & -698784 & 75216 & -149376 & -79960 & 5856 & 94920 & -12992 \\
 7104 & 8432 & 33960 & 11968 & 9180 & 5712 & 2380 & -544 \\
\end{array}
\right]
$$
\/}
{\footnotesize
  $$
\left[
\begin{array}{cccccccc}
 892440 & 3376800 & 410040 & 1157760 & 683400 & 96480 & -506520 & 64320 \\
 -73350 & -246888 & -228942 & -165792 & -110370 & -41688 & 27510 & -2064 \\
 1473 & 4092 & 10557 & 5808 & 4455 & 2772 & 1155 & -264 \\
\end{array}
\right]
$$
\/}

Define $3$ polynomials $P_3,P_4,P_6$ by the formula:
\begin{equation}
  P_k(s,t)=(1,s,s^2) \cdot M_k \cdot (1,...,t^7)=\sum_{i=0}^2 \sum_{j=0}^7 (M_k)_{ij} s^i t^j, \hskip 15 pt k=3,4,6.
\end{equation}
We have
\begin{equation}
 \gamma_s= P_3 3^{s/2} + P_4 4^{s/2} + P_6 6^{s/2}.
\end{equation}

To check the positivity of $\gamma_s$ we check that
each of the $16$ functions
\begin{equation}
  \gamma_s(v/4+u/4) = a_{v,0}(s)+a_{v,1}(s)u + ...  + a_{v,7}(s)u^7,
  \hskip 30 pt v=0,...,15
\end{equation}
satisfies the following condition:
$A_{v,k}(s)=a_{v,0}(s)+...+a_{v,k}(s)$ is positive for all $k=0,...,7$ and all $s \in [6,16]$.
The Positive Dominance Lemma now implies
that the corresponding polynomial is positive on $[0,1]$.

For each $v = 0,...,15$ and each $k=0,...,7$ we have a
$3 \times 3$ integer matrix $\mu_{v,k}$ such that
\begin{equation}
  \label{positive00}
  A_{v,k}(s)=(1,s,s^2) \cdot \mu_{v,k} \cdot (3^{s/2},4^{s/2},6^{s/2}).
\end{equation}
This gives $128$ matrices to check.  We get two more such matrices
from the conditions $-\beta_s'(0)>0$ and $\beta_s'(4)>0$.
All in all, we have to check that $130$ expressions of
the form in Equation \ref{positive00} are positive
for $s \in [6,16]$.  These expressions are all special
cases of Equation \ref{summand}, and we use the method
discussed above to show positivity in all $130$ cases.

\newpage

%% file: 5local.tex
\section{The Local Convexity Theorem}
\label{local}
\label{A111proof}

\noindent
{\bf Reading Guide:\/}
This chapter is for Reader 2. 
\newline

We use the notation from \S \ref{prelim}.
Recall that $\Sigma^{-1}$ is inverse stereographic projection.
The small domain $\Omega_0$ is defined in
\S \ref{definite}.  For each potential function $F$,
the symbol ${\cal E\/}_F$ denotes the $F$-energy
of configurations of $5$ points in $S^2$.
The map $\Sigma^{-1}$ carries avatars to
such configurations and so
${\cal E\/}_F \circ \Sigma^{-1}$ is a
function on our big space $\Omega$.
In this chapter we restrict these functions
to the small cube $\Omega_0$.
Here is the result we prove in this chapter.

\begin{theorem}[Local Convexity]
    For $F=G_4,G_6,G_5^{\flat},G_{10}^{\sharp}$, the
  Hessian of ${\cal E\/}_F \circ \Sigma^{-1}$ is positive
  definite at every point of $\Omega_0$.
\end{theorem}

\subsection{Discussion}

Before we launch into the proof, we discuss
why the proof has
the structure that it does.  We are interested
in showing that certain functions, essentially
the eigenvalues of the Hessian matrix of
various energy potentials, are positive
in a certain definite neighborhood.

Consider a toy version of this problem
where we want to show that a real
valued smooth function $f$ is positive
on an interval $[0,a]$.   We only want
to evaluate $f$ and its derivatives at $0$.
In general, this is a hopeless task, but
let us discuss it anyhow.

We evaluate $f(0)$ and we notice that it
is positive.  If we knew that $|f'|$ was small enough
on $[0,a]$ then we could use Taylor's Theorem with Remainder
to show that $f>0$ on $[0,a]$.
We could evaluate $|f'(0)|$ and observe that
it is quite small.  If we knew that
$|f''|$ was small on $[0,a]$ then we could again use
Taylor's Theorem to show that
$|f'|$ is small enough on $[0,a]$.

In general, we have a recursive problem
with no end in sight.   However, in our
specific situation we have some
good algebraic luck that saves us.
It turns out that certain 
{\it a priori\/} algebraic bounds on the
$n$th derivatives grow slowly in comparison
to the large constant $n!$ we divide by when
we apply Taylor's Theorem with Remainder.

We will use an algebraic trick to get
reasonable bounds on high derivatives
of the functions of interest to us, and then
we will use Taylor's Theorem with remainder
to promote these decent bounds on high
derivatives to excellent bounds on
the lower derivatives.

\subsection{Reduction to Simpler Statements}

Before we start we simplify our coordinates slightly.
All our avatars $$(x_0,x_1,x_2,x_3) \in \Omega_0 \subset \R^8$$ satisfy
$x_{02}=0$.  So, we drop the $x_{02}$ coordinate and instead think of
$\Omega_0$ as a subset of $\R^7$.

We consider $F$ to be any of the $4$ functions
$$
G_4, \hskip 15 pt G_6, \hskip 15 pt
G_5^{\flat} = G_5-25 G_1, \hskip 15 pt
2^{-5}G_{10}^{\sharp}=2^{-5}(G_{10}+13 G_5 + 68 G_2).$$
Scaling the last function by $2^{-5}$ makes our estimates more uniform.

Recall that $\Omega_0$ is the cube of side length $2^{-16}$
centered at the point
\begin{equation}
\xi_0=\bigg(1,0,\frac{-1}{\sqrt 3}, 
 -1,0, 0,\frac{1}{\sqrt 3}\bigg) \in \R^7
\end{equation}
In general, the point $(x_1,...,x_7)$ represents the  avatar
\begin{equation}
  p_0=(x_1,0), \hskip 3 pt
  p_1=(x_2,x_3), \hskip 3 pt
  p_2=(x_4,x_5), \hskip 3 pt
  p_3=(x_6,x_7).
\end{equation}
The quantity $F(x_1,...,x_7)$ is the $F$-potential of the
$5$-point
configuration $(\widehat p_0,...,\widehat p_4)$ associated to the  avatar under
inverse stereographic projection $\Sigma^{-1}$.
\begin{equation}
   F(x_1,...,x_7)=\sum_{i<j} F(\|\widehat p_i - \widehat
   p_j\|).
 \end{equation}
 Equation \ref{inversestereo} gives the formula
 for $\Sigma^{-1}$.

Let $HF$ be the Hessian of $F$.
The Local Convexity Theorem says $HF$ is positive definite in $\Omega_0$.
Let $\partial_J F$ be the (iterated)
    partial derivative of $F$ with respect to a
    multi-index $J=(j_1,...,j_7)$.
    Let $|J|=j_1+...+j_7$.
    Let
    \begin{equation}
      \label{thirdbound}
      M_N(F)=\sup_{|J|=N} M_J(F), \hskip 20 pt
      M_{J}(F)=\sup_{\xi \in \Omega_0} |\partial_J F(\xi)|.
\end{equation}
Let $\lambda(M)$ be the smallest eigenvalue of
a real symmetric matrix $M$.
The Local Convexity Theorem is an immediate consequence of the
following two lemmas.
    
    \begin{lemma}[L1]
      \label{L1}
  If $M_3(F)<2200 \lambda(HF(\xi_0))$ then
  $\lambda(HF(\xi))>0$ for all $\xi \in \Omega_0$.
\end{lemma}

\begin{lemma}[L2]
        \label{L2}
        $M_3(F)<2200 \lambda(HF(\xi_0))$.
  \end{lemma}

  \subsection{Proof of Lemma L1}

Let \begin{equation}
  H_0=HF(\xi_0), \hskip 30 pt
  H=HF(\xi), \hskip 30 pt
  \Delta=H-H_0.
\end{equation}

For any real symmetric matrix $X$ define the {\it Frobenius norm\/}:
\begin{equation}
  \|X\|_F=\sqrt{\sum_{ij}X_{ij}^2}.
\end{equation}
By the Cauchy-Schwarz inequality, we have
\begin{equation}
\sup_{\|v\|=1} \|X v\| \leq \|X\|_F.
\end{equation}

Given a unit vector $v \in \R^7$ we have
$H_0v \cdot v \geq \lambda(H_0).$ Hence
$$
Hv \cdot v = (H_0 v+\Delta v) \cdot v  \geq
H_0v \cdot v - |\Delta v \cdot v| \geq
\lambda(H_0) - \|\Delta v\| \geq
\lambda(H_0) - \|\Delta\|_F.
$$
So, to prove Lemma L1 we just need to establish the implication
$$
M_3<2200 \lambda(H_0) \hskip 15 pt \Longrightarrow \hskip 15 pt
\|\Delta\|_F<\lambda(H_0).
$$

Let $t \to \gamma(t)$ be the {\it unit speed parametrized\/}
line segment connecting $\xi_0$ to $\xi$ in $\Omega_0$.
Note that $\gamma$ has length 
$L \leq \sqrt 7 \times 2^{-17}$.
We write $\gamma=(\gamma_1,...,\gamma_7)$.
Let $H_t$ denote the Hessian of
$F$ evaluated at $\gamma(t)$.
Let $D_t$ denote the directional derivative along $\gamma$.

Now $\|D_t(H_t)\|_F$ is the speed of the path
$t \to H_t$ in $\R^{49}$, and
$\|\Delta\|_F$ is the Euclidean distance between the
endpoints of this path.  Therefore
\begin{equation}
  \label{speed}
\|\Delta\|_F \leq \int_0^{L} \|D_t(H_t)\|_F\ dt.
\end{equation}

Let $(H_t)_{ij}$ denote the $ij$th entry of $H_t$.
From the definition of directional derivatives, and
from the Cauchy-Schwarz inequality, we have

\begin{equation}
  (D_tH_t)_{ij}^2 =
  \bigg( \sum_{k=1}^7 \frac{d\gamma_k}{dt}
  \frac{\partial H_{ij}}{\partial \gamma_k} \bigg)^2
  \leq 7 M_3^2,
  \hskip 30 pt
  \|D_t(H_t)\|_F \leq 7^{3/2} M_3.
\end{equation}

The second inequality follows from summing the first one over all
$7^2$ pairs $(i,j)$ and taking the square root.  Equation
\ref{speed} now gives
\begin{equation}
\label{small2}
\|\Delta\|_F \leq
L \times 7^{3/2} M_3  \leq 49 \times 2^{-17} M_3<M_3/2200<\lambda(H_0).
\end{equation}
This completes the proof.

\subsection{Proof of Lemma L2}

Let $F$ be any of our functions.
Let $H_0=HF(\xi_0)$.

\begin{lemma}[L21]
  $\lambda(H_0)>39$.
      \end{lemma}

\startproof
Let $\chi$ be the characteristic polynomial of $H_0$.
This turns out to be a rational polynomial.  We check in Mathematica
    that the signs of the coefficients of
    $\chi(t+39)$ alternate, with the constant term nonzero.  Hence
    $\chi(t+39)$ has no negative roots.
    The file we use is {\tt L21.m\/}.
    \endproof

Recalling that $\xi_0 \in \R^7$ is the point representing the TBP, we define
 \begin{equation}
      \label{thirdbound2}
\mu_N(F)=\sup_{|I|=N}|\partial_I F(\xi_0)|.
\end{equation}
When $F$ is implicit in the discussion we write $\mu_N=\mu_N(F)$.
 
    \begin{lemma}[L22]
      For any of our functions we have the bound
\begin{equation}
  \label{eval}
  \mu_3<45893, \hskip 30 pt
  \frac{(7\times 2^{-17})^j}{j!}\mu_{j+3}<75, \hskip 15 pt
  j=1,2,3.
\end{equation}
\end{lemma}

\startproof
We compute this in Mathematica.
The file we use is {\tt  L22.m\/}.
\endproof

\begin{lemma}[L23]
  \label{L23}
  For any of our functions we have the bound
  $$  \frac{(7\times 2^{-17})^4}{4!}  M_7<37664.$$
\end{lemma}

\startproof
We give this proof in the next section.
\endproof

\begin{lemma}[L24]
We have
\begin{equation}
\label{BOUND}
M_3 \leq \mu_3+
\sum_{j=1}^3
\frac{(7 \times 2^{-17})^j}{j!}\mu_{j+3}
+\frac{(7 \times 2^{-17})^4}{4!}M_7
\end{equation}
\end{lemma}

\startproof
Choose any multi-index $J$ with $|J|=3$.
Let $\gamma$ be the line segment connecting
$\xi_0$ to any $\xi \in \Omega_0$.  We parametrize
$\gamma$ by unit speed and furthermore set $\gamma(0)=\xi_0$.  Let
$f(t)=\partial_J F \circ \gamma(t).$
The bound for $M_J$ follows from Taylor's Theorem with
remainder once we notice that
$$0 \leq t \leq \sqrt 7 \times 2^{-17}, \hskip 30 pt
\bigg|\frac{\partial^n f(0)}{\partial t^n}\bigg| \leq
(\sqrt 7)^n \mu _{n+3} \hskip 30 pt
\bigg|\frac{\partial^n f}{\partial t^n}\bigg| \leq
(\sqrt 7)^n M_{n+3}.$$
Since this works for all $J$ with $|J|=3$ we get the same bound for $M_3$.
\endproof

The lemmas above and Equation \ref{eval} imply
$$
M_3<45893+3\times 75+37664=83782
<85800=2200\times 39<2200\lambda(H_0).
$$
This completes the proof of Lemma L2.

\subsection{Proof of Lemma L23}
\label{combinatorics}

Now we come to the interesting part of the proof,
the one place where we need to go beyond
specific evaluations of our functions.
When $r,s \in \N \cup \{0\}$ and $r+s \leq 2d$ we have
\begin{equation}
  \label{basic00}
\sup_{(x,y) \in \R^2} \frac{x^r
  y^s}{(1+x^2+y^2)^d} \leq
(1/2)^{\min(r,s)}.
\end{equation}
One can prove Equation \ref{basic00} by factoring the expression
into pieces with quadratic denominators.
Here is a more general version.
Say that a function $\phi: \R^4 \to \R$ is {\it nice\/} if it has the
form
{\small
$$\sum_i \frac{C_i a^{\alpha_i} b^{\beta_i} c^{\gamma_i} d^{\delta_i}
}{(1+a^2+b^2)^{u_i}(1+c^2+d^2)^{v_i}},
\hskip 5 pt \alpha_i,\beta_i,\gamma_i,\delta_i \geq 0, \hskip 5 pt
\hskip 5 pt \alpha_i+\beta_i \leq 2u_i, \hskip 15 pt \gamma_i+\delta_i
\leq 2v_i.$$
\/}
It follows from Equation \ref{basic00} that
\begin{equation}
  \label{easy2}
  \sup_{\R^4} |\phi| \leq \langle \phi \rangle, \hskip 30 pt
  \langle \phi \rangle = \sum_i |C_i| (1/2)^{\min(\alpha_i,\beta_i)+\min(\gamma_i,\delta_i)}.
\end{equation}
Equation \ref{easy2} is useful to us because it allows us
to bound certain kinds of functions without having to evaluate
them anywhere.  We also note that
if $\phi$ is nice, then so is any iterated partial derivative of
$\phi$.   Indeed, the nice functions form a ring that is
invariant under partial differentiation.  This fact makes it
easy to construct nice functions.

For any $\phi: \R^n \to \R$
we define
    \begin{equation}
      \overline M_7(\phi)=\sup_{|J|=7} \overline M_J(\phi), \hskip 20 pt
      \overline M_{J}(\phi)=\sup_{\xi \in \R^n} |\partial_J (\phi)|.
    \end{equation}
    We obviously have
    \begin{equation}
M_7(F) \leq
\overline M_7(F).
\end{equation}

Recall that $\widehat p=\Sigma^{-1}(p)$, the inverse stereographic image of $p$.
Define 
\begin{equation}
f(a,b)=4-\|\widehat{(a,b)}-(0,0,1)\|^2=
\frac{4(a^2+b^2)}{1+a^2+b^2},
\end{equation}

  {\small 
\begin{equation}
g(a,b,c,d)=4-\|\widehat{(a,b)}-\widehat{(c,d)}\|^2=
\frac{4(1 + 2 a c + 2 b d + (a^2+b^2)(c^2+d^2))}
{(1 + a^2 + b^2) (1 + c^2 + d^2)}.
\end{equation}
\/}
Notice that $g$ is nice.  Hence $g^k$ is nice and
$\partial_I g^k$ is nice for any multi-index.
That means we can apply Equation \ref{easy2} to
$\partial_I g^k$.

${G_k}$ is a $10$-term expression involving
$4$ instances of $f^k$ and $6$ of $g^k$.  Each of these $10$ terms
occurs with coefficient $1$. However, each
variable appears in at most $4$ terms.  So, as soon as we
take a partial derivative, at least $6$ of the terms vanish.
Moreover, every partial derivative of $f$ is a limiting case of
the corresponding partial derivative of $g$.
From this, we see that
\begin{equation}
\label{bound1}
\overline M_7({G_k}) \leq 4 \times \overline M_7(g^k).
\end{equation}
Equation \ref{easy2} gives
\begin{equation}
  \label{bound2}
 4 \times  \overline M_7(g^k) \leq 4 \times \max_{|I|=7} \ \langle \partial_I g^k \rangle.
\end{equation}

Using this estimate, and
the Mathematica file {\tt L23.m\/}, we get
$$
  \max_{k \in \{1,2,3,4,5,6\}} \frac{(7 \times
    2^{-17})^4}{4!} \times  4 \times \overline M_7(g^k) \leq \frac{1}{2}.
  $$

  \begin{equation}
    \label{bound3}
2^{-5} \times \frac{(7 \times 2^{-17})^4}{4!} \times 4 \times
\overline M_7(g^{10}) \leq 37643.
\end{equation}

Together with the preceding estimate and the definition of
$G_{10}^{\sharp}$, this gives

$$
\frac{(7\times 2^{-17})^4}{4!}\,
M_7(2^{-5}G_{10}^{\sharp})
\leq
37643+2^{-5}\bigg(\frac{13}{2}+\frac{68}{2}\bigg)
<37645<37664.
$$

The remaining bounds in Lemma \ref{L23} follow directly from
Equations \ref{bound1}--\ref{bound3} and from the definitions
of our functions.

\newpage

%% file: 6symm.tex
\section{The Symmetrization Theorem}
\label{Bproof}

\noindent
{\bf Reading Guide:\/}
This is for Reader 3.    We prove
the Symmetrization Theorem from \S \ref{SYMM}.
We use two tricks. First, we write our
symmetrization operation as the composition
of three simpler things: rotation, then
horizontal symmetrization, then vertical symmetrization.
Second, and crucially, we partition
the $10$-term
sum into smaller sums, each involving just a
few terms, which are separately decreased
by these smaller operations.

\subsection{Rotation}

The domain $\Upsilon$ is defined in
\S \ref{special}. Let $X=(p_0,p_1,p_2,p_3)$ be an avatar in $\Upsilon$.
We let $X'$ be the planar
configuration which is obtained by
rotating $X$ about the origin so that
$p_0'$ and $p_2'$ lie on the same horizontal
line, with $p_0'$ lying on the right.
This operation does not change the $R_s$-energy.

\begin{lemma}
  \label{anglebound}
  The rotation angle $\theta$ needed to
  carry $X$ to $X'$ satisfies the inequality
  $0 \leq \theta <1/34$.
\end{lemma}

\startproof
Since $p_0$ lies on the $x$ axis and $p_{2}$ lies on or above it, we have to
rotate by a small amount counterclockwise to get $p'_0$ and $p'_2$ on the same
horizontal line.
Hence $\theta \geq 0$.
For the upper bound on $\theta$, write
{\small
$$512p_0=(a,0), \hskip 10 pt 512p_2=(b,c), \hskip 20 pt
a\in[433,498], \hskip 10 pt  b\in[-498,-400], \hskip 10 pt
c\in[0,24].$$
\/}
The condition that $\rho_\theta(p_0)$ and $\rho_\theta(p_2)$
have the same second coordinate gives
$$a\sin\theta=b\sin\theta+c\cos\theta,$$
and hence
$$\tan\theta=\frac{c}{a-b}.$$
This quantity is maximized when
$$a=433,\qquad b=-400,\qquad c=24.$$
Therefore
$$\tan\theta\leq \frac{24}{833}<\tan(1/34),$$
and hence $\theta<1/34$.
\endproof

For notational convenience, let $I_a=[-a,a]$ for any $a \geq 0$.
Let $\Upsilon'$ denote the domain of avatars $X'$ such that
(comparing with $\Upsilon$)
\begin{enumerate}
\item $\|p'_0\| \geq \|p'_k\|$ for $k=1,2,3$.
\item $512p'_0 \in [432,498] \times I_{16}$. (Compare $[433,498] \times I_0$.)
\item $512 p'_1 \in I_{32} \times [-465,-348]$. (Compare $I_{16} \times [-464,-349]$.)
\item $512 p'_2 \in [-498,-400] \times I_{16}$. (Compare $[-498,-400] \times [0,24]$.)
\item $512 p'_3 \in I_{32} \times [348,465]$. (Compare $I_{16} \times [349,464]$.)
\item $p_{02}'=p_{22}'$. (Compare $p_{02}=0$.)
\end{enumerate}

\begin{lemma}[B1]
  \label{symm0}
  If $X \in \Upsilon$ then $X' \in \Upsilon'$.
\end{lemma}

\startproof
Rotation about the origin does not change the norms, so $X'$
satisfies {\bf Condition 1}.   Next, we observe that
{\bf Condition 6} holds by construction.

Since $p'_{02}=p'_{22}$ and
$\|p'_0\|\geq \|p'_2\|$, we have
$|p'_{21}|\leq |p'_{01}|$.
Since
$$|p'_{01}|\leq \|p'_0\|=\|p_0\|\leq 498/512,$$
we have
$p'_{21}\geq -498/512$.
This gives us  {\bf Condition 4, Lower Bound}.

Now we check the other conditions.
Let $\theta$ be as in Lemma \ref{anglebound}.
Let $\rho_{\theta}$ denote the counterclockwise rotation through
the angle $\theta$.  By Lemma \ref{anglebound},
\begin{equation}
\label{trig}
512(1-\cos(\theta)) \in [0,1], \hskip 30 pt
512 \sin(\theta) \in [0,16].
\end{equation}

From Equation \ref{trig}, 
the map $512 p_0 \to 512 p_0'$ changes the first coordinate by
$512\delta_{01} \in [-1,0]$ and
$512\delta_{02} \in [0,16]$.  {\bf Condition 2} follows.
Next, we have $512 \delta_{21} \in [-1,0]$.   This gives
 {\bf Condition 4, Upper Bound}.

At this point we have checked {\bf Conditions 1, 2, 4, 6}.

{\bf Condition 3} follows from
$512 \delta_{11} \in [0,16]$ and $512 \delta_{12} \in [-1,1]$.
The first bound comes from $512 \sin(\theta)<16$.
For the second bound, let $\alpha$ be the signed angle
that $p_1$ makes with the negative $y$-axis.
This angle is largest in absolute value at a corner of the
constraint rectangle, namely when we have 
$512p_1=(\pm 16,-349)$.
Since
$$\tan(1/21)>16/349,$$
we have $|\alpha|<1/21$.
Also $\|p_1\|<1$.  Hence
$$
|512 \delta_{12}|
\leq 512 \max_{|x|\leq 1/21}
\bigg|\cos\bigg(x+\frac{1}{34}\bigg)-\cos(x)\bigg|<1.
$$
This gives {\bf Condition 3}.  The same argument
gives {\bf Condition 5}.
\endproof

\subsection{Horizontal and Vertical Symmetrization}

Now that we have rotated our configurations
slightly we will realize our symmetrization operation
as the composition of two subsidiary operations,
{\it horizontal symmetrization\/} and
{\it vertical symmetrization\/}.
\newline
\newline
\noindent
{\bf Horizontal Symmetrization:\/}
Given an avatar
$X' \in \Upsilon'$, there is
a unique configuration 
$X''$, invariant under
reflection in the $y$-axis, such that
$p_j'$ and $p_j''$ lie
on the same horizontal line for $j=0,1,2,3$, and
$p_0''$ lies on the right, and
$\|p_0''-p_2''\|=\|p_0'-p_2'\|$.
To get from $X'$ to $X''$ we are just changing
the first coordinates of the points -- i.e., we
are moving the points horizontally.
We call
this {\it horizontal symmetrization\/}.
In a straightforward way we see that
horizontal symmetrization maps $\Upsilon'$ into
$\Upsilon''$, the set of avatars
$p_0'',p_1'',p_2'',p_3''$ such that
\begin{enumerate}
\item $-512p_2'', 512 p''_0 \in [416,498] \times I_{16}$
\item $-512 p_1'', 512 p''_{3} \in I_0 \times  [348,465]$.
  \item $p''_{02}=p''_{22}$.
  \end{enumerate}
  Here is more justification for the first bound.
Horizontal symmetrization replaces the first coordinates
$p'_{01},p'_{21}$ by the symmetric pair
\[
\pm\frac{p'_{01}-p'_{21}}{2}.
\]
Since
$512p'_{01}\in[432,498]$ and
$512p'_{21}\in[-498,-400]$,
the new positive first coordinate lies in $[416,498]$.
\newline
\newline
{\bf Vertical Symmetrization:\/}
Recall from \S \ref{SYMM} that
${\bf K\/}_4$ is the set of avatars
which are symmetric with respect to
reflections in the coordinate axes.
Given a configuration 
$X''\in \Upsilon''$ there
is a unique configuration
$X'''\in {\bf K\/}_4$
such that $p_j''$ and $p_j'''$ lie on the same
vertical line for $j=0,1,2,3$. 
We call this operation {\it vertical symmetrization\/}.
All we are doing is sliding each pair $(p_0'',p_2'')$ and
$(p_1'',p_3'')$ rigidly up or down.
Here $X'''=X^*$ from \S \ref{SYMM}.

\subsection{Reduction to Three Lemmas}

Given an avatar $X=(p_0,p_1,p_2,p_3)$ define
\begin{equation}
  r_{ij}=\frac{1}{\|\Sigma^{-1}(p_i)-\Sigma^{-1}(p_j)\|}.
\end{equation}
Given a list
$L$ of pairs of points we define
$R_s(X,L)$ to be the sum of the $R_s$-potentials just over
the pairs in $L$.  E.g. $R_s(X,\{(0,2),(0,4)\})=r_{02}^s+r_{04}^s$.

We call $L$ {\it good\/} for $s$, and with
respect to one of the operations, if the operation does not increase
the value of $R_s(X,L)$.  We call $L$ {\it great\/} if
the operation strictly lowers $R_s(X,L)$ unless the operation
fixes $X$.
We mean to take the appropriate domains in all cases.
Here are three more lemmas we prove below.

\begin{lemma}[B2]
  \label{symm1}
  On $\Upsilon$, 
  $\{(0,2),(0,4),(2,4)\}$ and $\{(1,3),(1,4),(3,4)\}$ are good for all $s \geq 2$
  and with respect to symmetrization. Also, the
  union of these two lists is great.
\end{lemma}

\begin{lemma}[B3]
  \label{symm3}
  On $\Upsilon'$, the lists $\{(0,1),(1,2)\}$ and $\{(0,3),(3,2)\}$ are both good for all $s \geq 2$
  and with respect to horizontal symmetrization.
\end{lemma}

\begin{lemma}[B4]
  \label{symm4}
  on $\Upsilon''$, the lists
  $\{(0,1),(0,3)\}$ and $\{(2,1),(2,3)\}$ are both good for all $s \geq 12$
  and with respect to vertical symmetrization.
\end{lemma}

\noindent
{\bf Proof of the Symmetrization Theorem:\/}
We will deduce the Symmetrization Theorem from
Lemmas B2,B3,B4.   (Lemma B1 guarantees the
Lemmas B3 and B4 are applied to the correct domains.)
Note that all these lemmas are valid
for any power law exponent $s \geq 12$.

We have $10=6+4$ terms in the
energy formula we need to track under the symmetrization
operation.
Consider the following $6$ terms.
\[
(0,2),(0,4),(2,4),(1,3),(1,4),(3,4)
\]
Lemma B2 says that each of the two sums corresponding to
\[
\{(0,2),(0,4),(2,4)\},\qquad
\{(1,3),(1,4),(3,4)\}
\]
is non-increasing under symmetrization.
Hence the $6$-term sum is non-increasing
under symmetrization.

The remaining $4$ terms are
\[
(0,1),(1,2),(0,3),(2,3).
\]
Rotation from $X$ to $X'$ keeps all terms the same.
Lemma B3 says that the sum of these $4$ terms is non-increasing
under horizontal symmetrization $X'\to X''$.
Lemma B4 says  the same thing for
vertical symmetrization $X''\to X'''=X^*$.
Moreover, the union of the two lists in Lemma B2 is great, so if $X\notin {\bf K}_4$
the decrease is strict.
\endproof

\subsection{Proof of Lemma B2}
\label{B2proof}

Let $s_3=\sqrt 3/3$.
Inverse stereographic projection maps the
triangle with vertices $(\pm s_3,0)$ and $\infty$
to an equilateral triangle on $S^2$.
Avatars in $\Upsilon$ satisfy
$$\|p_0\|,  \|p_1\|, \|p_2\|, \|p_3\|, \frac{\|p_0-p_2\|}{2},  \frac{\|p_1-p_3\|}{2} \in
(s_3,1).$$

Let $(u,v)$ stand for either $(0,2)$ or $(1,3)$.
\begin{enumerate}
\item Let $a>0$ be such that
  $\|p_u-p_v\|/2 =s_3+a$.  Let
  $-q_u=q_v=(s_3+a,0)$.
  The points $q_u,q_v$ are symmetric w.r.t the $y$-axis.
  Also set $a_u=a_v=a$.
\item Choose $b_u,b_v$ with
  $0<b_u \leq a_u$ and $0<b_v \leq a_v$.
  Let $r_u=(-s_3-b_u,0)$ and $r_v=(s_3+b_v,0)$.
  Note that $\|r_u-r_v\| \leq \|q_u-q_v\|$.
\end{enumerate}
Up to rotation about the origin, our
symmetrization operation does the map
$(p_u,p_v) \to (r_u,r_v)$ for suitable $a_u,a_v,b_u,b_v$.
For our symmetrization operation we have the
additional properties $b_0=b_2=a$ and $b_1=b_3$,
but we want to consider the more general case as
part of our proof strategy.

Recall that $\widehat p$ is the image of $p$ under inverse
stereographic projection.
Lemma B2 is implied by:
$$
  \|\widehat r_u-\widehat r_v\|^{-s} +
  \|\widehat r_u-(0,0,1)\|^{-s} +
  \|\widehat r_v-(0,0,1)\|^{-s}  \leq$$
  \begin{equation}
    \label{maingoal}
  \|\widehat p_u-\widehat p_v\|^{-s} +
  \|\widehat p_u-(0,0,1)\|^{-s} +
  \|\widehat p_v-(0,0,1)\|^{-s}
\end{equation}
for all $s \geq 2$, 
  with equality iff $(r_u,r_v)=(p_u,p_v)$ up to rotation about the
  origin.    We will establish this in two steps.
  
  \begin{lemma}[B21]
  Let $s \geq 2$ and
  $$
  A_s=\|\widehat p_u-\widehat p_v\|^{-s} -
  \|\widehat q_u-\widehat q_v\|^{-s},$$
  {\small
    $$B_s= \|\widehat p_u-(0,0,1)\|^{-s} + \|\widehat p_v-(0,0,1)\|^{-s} -  \|\widehat q_u-(0,0,1)\|^{-s} - \|\widehat q_v-(0,0,1)\|^{-s}.$$
    \/}
  Then $A_s, B_s \geq 0$, with
  equality iff $p_u=q_u$ and $p_v=q_v$
  up to a rotation.
\end{lemma}

\startproof
Note that if $A_2>0$ then $A_s>0$ for all $s>0$.    If
$B_2>0$ then the Convexity Lemma implies that
$B_s>0$ for all $s>2$.  The
Convexity Lemma applies here because the last
two terms in the expression for $B_s$ are equal.

So, it suffices to prove that
$A_2,B_2>0$.
Rotating $p_u,p_v$ simultaneously
about the origin does not change any of the quantities
involved.  To simplify the calculation we make such a rotation, so
that
\begin{equation}
  \label{norma}
  p_u=(-x+h,y), \hskip 10 pt p_v=(x+h,y),
  \hskip 10 pt
  q_u=(-x,0), \hskip 10 pt q_v=(x,0).
\end{equation}
We compute
\begin{equation}
  A_2=\frac{h^4 + y^2 (2 + 2 x^2 + y^2) + 2h^2 (1 - x^2 + y^2)}{16 x^2},
  \hskip 15 pt B_2=\frac{y^2+h^2}{2}.
\end{equation}
Since $x \in (0,1)$ we have $A_2, B_2>0$ unless $h=y=0$.
\endproof

Define
\begin{equation}
  F_s(a_u,a_v) =
  \|\widehat q_u-\widehat q_v\|^{-s} +
  \|\widehat q_u-(0,0,1)\|^{-s} +
  \|\widehat q_v-(0,0,1)\|^{-s},
\end{equation}
Likewise define $F_s(b_u,b_v)$. This is the same expression
with respect to $\widehat r_u$ and $\widehat r_v$.
Finally, define
\begin{equation}
  E(s)=F_s(a_u,a_v)-F_s(b_u,b_v).
\end{equation}

\begin{lemma}[B22]
  For $s \geq 2$ we have
  $E(s) \geq 0$ with equality iff
     $b_u=a_u$ and $b_v=a_v$.
  \end{lemma}

  \startproof
  We first consider the case when $a_v=b_v$ and
  $b_u<a_u$.
With the file {\tt LemmaB22.m\/} we
compute that $\partial F_2/\partial a_u$ and
$-\partial F_{-2}/\partial a_u$ are both
rational functions of $a_u,a_v$ with all positive coefficients.
Therefore, we have $E(2)>0$ and $E(-2)<0$.

Referring to \S \ref{ESUM}, consider
the sign sequence for $E(s)$.
When $a_v=b_v$, the expression $E(s)$ is an exponential sum with $4$
terms.
When $a_u=a_v=0$ the points $\widehat q_u, \widehat q_v$ and
$(0,0,1)$ make an equilateral triangle on a great circle.  Hence,
when $a_u,a_v,b_u,b_v>0$ the point $\widehat q_u$ is closer to
$(0,0,1)$ than it is to $\widehat q_v$ both in its old location and in its new location.
The inward motion of the point $q_u$ increases the  shorter (corresponding spherical)
distance and
decreases the longer (corresponding spherical) distance. More to the point, our move decreases the
longer inverse-distance and increases the shorter inverse-distance.
Thus the sign sequence for $E(s)$
is $+,-,-,+$.

By Descartes' Lemma, $E(s)$ changes sign at most twice
and also $E(s)>0$ when $|s|$ is sufficiently large.
Since $E(-2)<0$ we see that
$E$ changes sign on $(-\infty,-2)$.  If
$E$ has a root in $(2,\infty)$ then in fact $E$ has at least $2$ roots
(counted with multiplicity) because it starts and ends positive on this interval.
But then $E$ has at least $3$ roots, counting multiplicity.  This is a
contradiction.
Hence $E(s)> 0$ for $s \geq 2$.

Interchanging the roles of $u$ and $v$, we see that
the same argument handles the case when $a_u=b_u$
and $b_v<a_v$.
Now consider the remaining case when
$b_u<a_u$ and $b_v<a_v$.  Using the
special cases we have already considered, we have
$$
F_s(a_u,a_v)>F_s(a_u,b_v)>F_s(b_u,b_v).
$$
This completes the proof.
\endproof

Lemmas B21 and B22 together prove that both lists are good.
Moreover, equality for both lists occurs only when
$X\in {\bf K}_4$.  Hence their union is great.

\subsection{Proof of Lemma B3}
\label{B3proof}

The domain $\Upsilon'$ is invariant under reflection in the
$X$-axis followed by interchanging the labels $1$ and $3$.
Thanks to this symmetry, it
suffices to prove Lemma B3 for the list $\{(0,3),(3,2)\}$.
To simplify some of our notation, we set
$q_0=p_0'$, $q_1=p_3'$, $q_2=p_2'$,
and
$q_0'=p_0''$, $q_1'=p_3''$, $q_2'=p_2''$.

We introduce the notation $q_1=(q_{11},q_{12})$, etc.
The horizontal symmetrization operation is
given by $$(q_0,q_1,q_2) \to (q_0',q_1',q_2'),$$ where
\begin{equation}
  q'_0=\bigg(\frac{q_{01}-q_{21}}{2},q_{02}\bigg), \hskip 20 pt
  q'_1=(0,q_{12}), \hskip 20 pt
  q'_2=\bigg(\frac{q_{21}-q_{01}}{2},q_{22}\bigg), \hskip 20 pt
\end{equation}
Note that $\|q_0'-q_1'\|=\|q_2'-q_1'\|$.
This means that the kind of inequality we are trying
to establish has the form $2A^s \leq B^s+C^s$ for
choices of $A,B,C$ which depend on the points involved.
Therefore, by the Convexity Lemma,
it suffices to prove that $\{(0,1),(1,2)\}$ is good for the parameter
$s=2$.

Let $D$ denote the set of triples of points $(q_0,q_1,q_2) \in (\R^2)^3$ such
that there is some $q_3$ such that $(q_0,q_1,q_2,q_3) \in \Upsilon'$.
 Most of our proof involves finding a concrete parametrization
of a subset of $\R^6$ that contains $D$.   Note that $D$ is
really a $5$ dimensional set, because $q_{22}=q_{02}$.
We will use parameters $a,b,c,d,e$ to parametrize a subset of
$\R^6$ that contains $D$.

We define
\begin{equation}
  \label{convention1}
  [a,b,t]=\frac{a(1-t)}{512}+\frac{bt}{512}.
\end{equation}
Hence $[a,b,\cdot]$ maps
the interval $[0,1]$ onto the interval with endpoints
$a/512$ and $b/512$ (in some order).
Given $(a,b,c,d,e) \in [0,1]^5$ and
$\sigma_1,\sigma_2 \in \{-,+\}$ we define
\begin{eqnarray}
  \begin{tabular}{lll}
$r_0=\big([+416,+498,a]+[0,49,e],[0,16 \sigma_1,b]\big);$\cr
$r_1=\big([0,32 \sigma_2,d],[348,465,c]\big);$\cr
$r_2=\big([-416,-498,a]+[0,49,e],[0,16 \sigma_1,b]\big);$
    \end{tabular}
  \end{eqnarray}
     We call this map $\phi_{\sigma_1,\sigma_2}$.
In these coordinates, horizontal symmetrization is the map
\begin{equation}
  (a,b,c,d,e) \to (a,b,c,0,0).
\end{equation}
We have two steps we need to take. First we really need to show that
we have parametrized a superset of $D$.  Second, we need to calculate
the energy change as a function of $a,b,c,d,e$ and check that it
decreases.

\begin{lemma}[B31]
  We have
  $$D \subset \phi_{+,+}([0,1]^5) \cup \phi_{+,-}([0,1]^5) \cup \phi_{-,+}([0,1]^5) \cup \phi_{-,-}([0,1]^5).$$
\end{lemma}

\startproof
Recall that $q_i=(q_{i1},q_{i2})$.
Let $D_{ij}$ denote the set of possible coordinates $q_{ij}$ that can arise
for points in $D$.  Thus, for instance $$D_{02} \subset [-16,16]/512.$$
Let $D_{ij}^*$ denote the set of possible coordinates $q_{ij}$ that can
arise from the union of our parametrizations.   By construction
$D_{i2} \subset D_{i2}^*$ for $i=0,1,2$ and
$D_{11} \subset D_{11}^*$.    The only coordinates not controlled by
$b,c,d$ and the sign choices $\sigma_1,\sigma_2$ are
the first coordinates $q_{01}$ and $q_{21}$ of $q_0$ and $q_2$
respectively.  We now deal with these.

Remembering that $q_{01} \geq |q_{21}|$, set
$x=512q_{01}$ and
$y=512q_{21}$.
Then
$$
432\leq x\leq 498, \qquad
-498\leq y\leq -400, \qquad
x\geq |y|.
$$
Define
$$
t=\frac{x+y}{2}, \qquad u=\frac{x-y}{2}.
$$
This gives $(x,y)=(t+u,t-u)$.
The inequalities above give
$ t \in [0,49]$ and $u \in [416,498]$.
Hence we can choose $a,e\in[0,1]$ so that
$$
[0,49,e]=\frac{t}{512}, \hskip 30 pt
[416,498,a]=\frac{u}{512},
$$
For these choices, the first coordinates of $r_0$ and $r_2$
in our parametrization are precisely $q_{01}$ and $q_{21}$.
\endproof

Using our coordinates above, we define
$$F_{\pm,\pm}(a,b,c,d,e) = \|\widehat q_0-\widehat q_1\|^{-2} + \|\widehat q_2-\widehat q_1\|^{-2},$$
\begin{equation}
  \Phi_{\pm,\pm}(a,b,c,d,e)={\rm num\/}_+(F_{\pm,\pm}(a,b,c,d,e)-F_{\pm,\pm}(a,b,c,0,0)).
\end{equation}
Here $q_0,q_1,q_2$ are the points which correspond to
$(a,b,c,d,e)$ under our map $\phi_{\pm,\pm}$ and
$\widehat q_0,\widehat q_1,\widehat q_2$ are their
images under inverse stereographic projection.
To finish our proof, we just have to
show that $\Phi_{\pm,\pm}(a,b,c,d,e) \geq 0$ on
$[0,1]^5$.   The following lemma, and continuity,
gives us this result.

\begin{lemma}[B32]
  For any sign choice, $\Phi_{\pm,\pm} \geq 0$ on  $[0,1]^5$.
\end{lemma}

\startproof
We let $\Phi_a=\partial \Phi/\partial a$, and likewise for
the other variables.  Iterating this notation, we
let $\Phi_{aa}$, etc., denote the second partials.

Let $\Phi$ be any of the $4$ polynomials.
 The file {\tt LemmaB32.m\/}  computes that
 \begin{enumerate}
  \item $\Phi$ and $\Phi_d$ and $\Phi_e$ are zero when $d=e=0$.
  \item $\Phi_{dd}$ and $\Phi_{ee}$ are weak positive dominant, hence nonnegative on $[0,1]^5$. 
  \item $\Phi_d+2\Phi_e$ is weak positive dominant, hence nonnegative on $[0,1]^5$.
  \end{enumerate}
  Let $Q_d \subset [0,1]^5$ be the sub-cube where $d=0$.
  Likewise define $Q_e$.
  We fix $(a,b,c)$ and consider the single variable
  function $\phi(d)=\Phi(a,b,c,d,0)$.
From Items 1 and 2 above,
$\phi(0)=\phi'(0)=0$ and $\phi''(d) \geq 0$.
Hence $\phi(d) \geq 0$ for $d \geq 0$.  Hence
$\Phi \geq 0$ on $Q_e$.
A similar argument shows that
$\Phi \geq 0$ on $Q_d$.

Any point in $(0,1)^5$ can be joined to a point in
$Q_{d} \cup Q_e$ by a line segment $L$ which is parallel
to the vector $(0,0,0,1,2)$.  From Item 3 above, $\Phi$ increases
along such a line segment as we move out of
$Q_d \cup Q_e$. Hence $\Phi  \geq  0$ on $[0,1]^5$.
\endproof

\subsection{Proof of Lemma B4}
\label{B4proof}

The set $\Upsilon''$ is invariant under reflection in the
$Y$-axis followed by interchanging the labels $0$ and $2$,
and under reflection in the $X$-axis followed by
interchanging the labels $1$ and $3$.
Thanks to these symmetries, it suffices to prove that
$\{(0,1),(0,3)\}$ is good for all $s\geq 12$,
and it suffices to consider the case when
$p_{02}''\geq0$.
That is, the point
$p''_0$ lies on or above the $X$-axis.
 For ease of notation set $q_k=p''_k$
and $q_k'=p_k'''$.   We are
considering the case when $q_{02} \geq 0$.

Let $D$ be the set of configurations
$(q_0,q_1,q_3)$ such that $q_{02} \geq 0$ and
$(q_0,q_1,q_2,q_3) \in \Upsilon''$ when
$q_2$ is the reflection of $q_0$ in the $Y$-axis.
Let $D_{\pm} \subset D$ denote
those configurations with
$\pm (q_{12}+q_{32}) \geq 0$. Obviously
$D=D_+ \cup D_-$.

The sets $D_{\pm}$ are $4$-dimensional subsets
of $(\R^2)^3$.  We parametrize a superset of
$D_{\pm}$ much as we
did in the proof of Lemma B3.
As in Equation \ref{convention1} we define
$$[a,b,t]=\frac{(1-t)a}{512} +
\frac{bt}{512}.$$
Given $(a,b,c,d) \in [0,1]^4$ and $\sigma \in \{+,-\}$ we define
\begin{eqnarray}
  \begin{tabular}{lll}
    $r_0=([416,498,b],[0,16,d])$; \cr
    $r_1=(0,-[348,465,a] + [0,59 \sigma,c])$; \cr
    $r_3=(0,+[348,465,a] + [0,59\sigma ,c])$;
  \end{tabular}
\end{eqnarray}
We call this map $\phi_{\sigma}$.
In these coordinates, the symmetrization operation is
$(a,b,c,d) \to (a,b,0,0)$.

\begin{lemma}[B41]
  $D_{\pm} \subset \phi_{\pm}([0,1]^4)$.
\end{lemma}

\startproof
This is just like the proof of Lemma B31.  The
only non-obvious point is why every pair $(q_{12},q_{32})$ is
reached by the map $\phi_{\pm}$. The essential point is
that for configurations in $D_{\pm}$ we
have $512 |q_{12}+q_{32}| \leq 2 \times 59$.
\endproof

Following the same idea as in the proof of Lemma B3, we define
\begin{equation}
  \label{phidef1}
  F_{s,\pm}(a,b,c,d)=\|\Sigma^{-1}(r_0) - \Sigma^{-1}(r_{1})\|^{-s}+ \|\Sigma^{-1}(r_0) - \Sigma^{-1}(r_{3})\|^{-s},
\end{equation}
\begin{equation}
  \label{phidef2}
  \Phi_{s,\pm}(a,b,c,d)={\rm num\/}_+(F_{s,\pm}(a,b,c,d)-F_{s,\pm}(a,b,0,0)).
\end{equation}
We can finish the proof by showing that
$\Phi_{2,+}\geq 0$ and $\Phi_{12,-}\geq 0$ on $[0,1]^4$. The Convexity Lemma then
takes care of all exponents greater than $2$ on $D_+$ and all exponents
greater than $12$ on $D_-$.  Notice the asymmetry in the
calculation. The
$(+)$ side is much less delicate.

\begin{lemma}[B42] $\Phi_{2,+} \geq 0$ on $[0,1]^4$.
\end{lemma}

\startproof
Let $\Phi=\Phi_{2,+}$.
Let $\Phi|_{c=0}$ denote the polynomial we get by
setting $c=0$. Etc.
Let $\Phi_c=\partial \Phi/\partial c$, etc.
The Mathematica file {\tt LemmaB42.m\/} computes that
$\Phi|_{c=0}$ and  $\Phi|_{d=0}$ and $\Phi_c + \Phi_d$ are weak
positive dominant.
Hence $\Phi \geq 0$ when $c=0$ or $d=0$, and the directional
derivative of $\Phi$ in the direction $(0,0,1,1)$ is
non-negative.
Hence $\Phi \geq 0$ on $[0,1]^4$.
\endproof

\begin{lemma}[B43] $\Phi_{12,-} \geq 0$ on $[0,1]^4$.
\end{lemma}

\startproof
The file {\tt LemmaB43.m\/} has the calculations.
    Let $\Phi=\Phi_{12,-}$.  This monster has $102218$ terms.
    \newline
    \newline
    {\bf Step 1:\/}
    Let $M$ denote the maximum coefficient of $\Phi$.  We let $\Phi^*$ be the polynomial we get by taking each
    coefficient $C$ of  $\Phi$ and replacing it with ${\rm floor\/}(10^{10} C/M)$.
          Note that if $\Phi^*$ is nonnegative on $[0,1]^4$
            then so is $\Phi$.
        \newline
        \newline
        {\bf Step 2:\/}
        Now $\Phi^*$ has $37760$ monomials in which the coefficient is $-1$.
        We check that each such monomial is divisible by one of $c^2$
        or $d^2$ or $cd$.  Let
               $\Psi=\Phi^{**} - 37760 (c^2+d^2+cd),$
    where $\Phi^{**}$ is obtained from $\Phi^*$ by setting all the $(-1)$ monomials to $0$.
    We have $\Psi \leq \Phi^*$ on $[0,1]^4$.  Hence, if $\Psi$ is non-negative on $[0,1]^4$ then so is $\Phi^*$.
    The polynomial $\Psi$ has $5743$ terms.
    \newline
    \newline
    {\bf Step 3:\/}
    We check that
            $\Psi_{aaa}$ is WPD and hence non-negative on $[0,1]^4$.  This  massive calculation
    reduces us to showing that
    the restrictions $\Psi|_{a=0}$ and $\Psi_a|_{a=0}$ and
    $\Psi_{aa}|_{a=0}$ are all non-negative on $[0,1]^3$.
    Consider
    \begin{equation}
      \label{hailmary}
       f|_{c=0}, \hskip 30 pt f|_{d=0} \hskip 30 pt  4f_c + f_d,
    \end{equation}
    We show that all three functions are WPD when either
    $f=\Psi_a|_{a=0}$ or $f=\Psi_{aa}|_{a=0}$.  This shows that $\Psi_a|_{a=0}$ and
    $\Psi_{aa}|_{a=0}$ are non-negative on $[0,1]^3$.   Also, we show
    that the first two functions are WPD when
    $f=\Psi|_{a=0}$.
    \newline
    \newline
    {\bf Step 4:\/}
   Let $f=\Psi|_{a=0}$ and
   $g=4f_c + f_d $.  We claim that $g\geq 0$ on $[0,1]^3$.
    We check that $g_d$ is WPD and hence non-negative on $[0,1]^3$.
    This reduces us to showing that $h=g|_{d=0}$ is non-negative on $[0,1]^2$.
    Here $h$ is a $2$-variable polynomial in $b,c$.  Referring to the
    operation in \S \ref{oper}, we check that the two
    subdivisions $S_{b,0}(h)$ and $S_{b,1}(h)$ are WPD.
    This proves that $h\geq 0$ on $[0,1]^2$.
    \endproof

\newpage

%% file: 7endgame.tex
\section{Symmetric Configurations}

\noindent
{\bf Reading Guide:\/}
This chapter is for Reader 4.
We prove the Critical Theorems from \S \ref{CRIT}. 
We use the notation from \S \ref{CRIT}.
This proof is rather technical and it appears to be
{\it ad hoc\/}.  However, I did extensive
numerical investigations of what was going on
with these functions before looking for a rigorous
proof.

\subsection{Critical Theorem I}

As in Equation \ref{ENDGAME}, we write
$(z,z)=\sigma(x,y)$.
Let $\phi: [0,1]^2 \to \Psi_4^{\sharp}$ be the affine map
which scales the coordinates by a factor of $1/64$.
We use coordinates $a,b$ on $[0,1]^2$ so that
$(x,y)=\phi(a,b)$.

For any rational function $F: \Psi_4^{\sharp} \to \R$ we define
\begin{equation}
  N_F(a,b)= \frac{{\rm num\/}_+((F - F \circ \sigma) \circ
    \phi)}{(a-b)^2}.
\end{equation}
See \S \ref{numplus}.
For all the choices of $F$ we make,
$N_F$ will be a polynomial.

Recall that $\Sigma^{-1}(p_4)=(0,0,1)$, and define
\begin{equation}
  \label{RIJ}
  r_{ij}=\frac{1}{\|\Sigma^{-1}(p_i)-\Sigma^{-1}(p_j)\|}.
\end{equation}
Note that $r_{ij}^s$ is a rational
function when $s$ is an even integer.

Let $R_s(x,y)$ be the $R_s$-energy of the avatar represented by $(x,y)$.
We write $R_s(x,y)=G_s(x,y)+H_s(x,y)$, where
\begin{equation}
  G_s= r_{02}^s + r_{13}^s, \hskip 20 pt
  H_s=2 r_{04}^s+2 r_{14}^s+ 4 r_{01}^s.
\end{equation}

\begin{lemma}[C1]
$G_s - G_s \circ \sigma \geq 0$ on
$\Psi_4^{\sharp} \times (2,\infty)$.
Equality happens only for configurations
which are FPs.
\end{lemma}

\startproof
The file {\tt LemmaC1.m\/} computes that
$N_{G_2}$ is a WPD polynomial.  Hence
$G_s - G_s \circ \sigma \geq 0$ for $s=2$. 
Let $\overline \zeta=\sigma(\zeta)$.
With the obvious notation, notice that
$\overline r_{02}=\overline r_{13}$.
Given this fact, and the form of
$G_s$, we now see from the Convexity
Lemma that our result extends from
$s=2$ to all $s>2$.  Finally, if
$\zeta$ is not an FP, then the
strict inequality part of the
Convexity Lemma gives
us $G_s(\zeta)>G_s(\overline \zeta)$
for $s>2$.
\endproof

\noindent
{\bf Remark:\/} We do not care about results for
the value $s=2$ in this chapter, but sometimes it
is useful to do our calculations at this value to
establish results we care about for larger $s$.
\newline

To finish the proof, we need to show
\begin{lemma}[C2] $H_s-H_s \circ \sigma \geq 0$ on
  $\Psi_4^{\sharp} \times [14,16]$.
\end{lemma}

We first prove two smaller lemmas and then deduce Lemma C2.

We will suppose, for the sake
of contradiction, that there exists some
$(x,y) \in \Psi_4^{\sharp}$ and some
$s_0 \in [14,16]$ such that 
\begin{equation}
  \label{assumebad}
  h(s_0)=H_{s_0}(x,y)-H_{s_0}(z,z)<0.
\end{equation}
Perturbing slightly, we can assume that
$(x,y)$ lies in the interior of $\Psi_4^{\sharp}$.
We study the single-variable function $h$.
The idea is to use Descartes' Lemma
from \S \ref{ESUM} to get a contradiction.
We first need some preliminary results.

\begin{lemma}[C21]
  \label{3root}
  $h$ has at least $3$ roots in $[2,16]$.
\end{lemma}

\startproof
The file {\tt LemmaC21.m\/}
computes that
$-N_{H_2}$ and $N_{H_{14}}$ and $N_{H_{16}}$ are all WPD polynomials.
Hence $h(2)<0$ and $h(14)>0$ and $h(16)>0$.
Given our assumption, Equation
\ref{assumebad}, we have $h(s_0)<0$.  So,
we see by the Intermediate Value Theorem that $h$ has roots
in each interval $(2,14)$ and $(14,s_0)$ and $(s_0,16)$.
Hence $h$ has at least $3$ roots in $[2,16]$.
\endproof
   
 Let $(p_0,p_1,p_2,p_3)$ and
    $(p'_0, p'_1, p'_2, p_3')$ respectively be the configurations
    corresponding to $(x,y)$ and $(z,z)=\sigma(x,y)$.  Without
    necessarily ordering the terms so that their bases occur in
    order of size, we have

    \begin{equation}
      \label{beta}
      h(s)=+ 2r_{04}^s - 4(r'_{04})^s + 2r_{14}^s + 4r_{01}^s - 4(r'_{01})^s.
    \end{equation}

    The next result gives us control on the ordering of these terms.
    For notational convenience, we set
    $r_0=r_{04}$ and  $r_1=r_{14}$ and
    $r_0'=r'_{04}$.
    
    \begin{lemma} [C22]
      $r_0,r_1,r'_0<1/\sqrt 2<r_{01},r_{01}'$ and $r_{01}<r_{01}'$.
    \end{lemma}

    \startproof
    We have $x,y,z \in (0,1)$.
    We compute
    {\small
    $$(1/2)-r_0^2=\frac{1-x^2}{4}>0, \hskip 15 pt
    (1/2)-r_1^2= \frac{1-y^2}{4}>0, \hskip 15 pt
      (1/2)-(r_0')^2=\frac{1-z^2}{4}>0,$$
    $$(r_{01})^2-(1/2)=\frac{(1-x^2)(1-y^2)}{4(x^2+y^2)}>0, \hskip 20 pt
      (r_{01}')^2-(1/2)=\frac{(1-z^2)^2}{8z^2}>0.$$
      \/}
    This proves the first statement.
    For the second statement, 
     {\tt LemmaC22.m\/} computes that $-N_{r_{01}^2}$ is a WPD polynomial.
     Hence $r_{01}'  \geq  r_{01}$.
     If we really had $r_{01}'=r_{01}$, then Equation \ref{beta} would
     only
     have $3$ terms.  There would then be at most $2$ sign changes
     and Lemma \ref{3root} would contradict Descartes' Lemma.
     We conclude that $r_{01}'>r_{01}$.
     \endproof

     Since the largest term in
     Equation \ref{beta} is  $- 4(r'_{01})^s$ we see that
     $h$ vanishes for some $s>16$.  Combining this
     with Lemma \ref{3root} and the fact that
     $h(16)>0$, we see that $h$ changes sign
     $4$ times on $[2,\infty)$.
          This is only possible if the sign sequence is
     $-+-+-$.  But this is impossible because there are
     three pluses and two minuses.  We have a contradiction.
     The only way out is that $h \geq 0$ on $[14,16]$.
     This proves Lemma C2 and thereby finishes the proof
     of the Critical Theorem I.

     \subsection{Critical Theorem II modulo Lemma C3}
     
     \noindent
     {\bf Derivative Bounds:\/}
As above, we identify $\Psi_4$ with a square in $\R^2$.
The point $(1,\sqrt 3/3)$, which is outside $\Psi_4$,
names the TBP.
We define
\begin{equation}
  \label{THETA}
  \Theta(s,x,y)=R_s(x,y)-R_s(1,\sqrt 3/3).
\end{equation}
Here we are comparing the $R_s$-energy of
an avatar in $\Psi_4$ to the $R_s$ energy of
the TBP.    Let $\Theta_x$ be the partial derivative
    of $\Theta$ with respect to $x$, etc.
    In \S \ref{BOUNDC} we establish the
    following bound.
\begin{lemma}[C3]
  $|\Theta_{xx}|, |\Theta_{yy}| \leq 4$
  and $|\Theta_{ss}| \leq 1/64$ 
    on $\Psi_4 \times [13,16]$.
\end{lemma}

\noindent
{\bf Blocks:\/}
We say that a {\it block\/} is a
rectangular solid of the form
\begin{equation}
  X = I \times Q \subset [0,16] \times  [0,1]^2.
\end{equation}
where $Q$ is a square and $I$ is an interval.
We define $|X|_1$ to be the length of $I$ and
$|X|_2$ to be the side length of $Q$.
Let $v(X)$ denote the  vertex set of $X$.

\begin{lemma}
  \label{C22}
  For any block
  $X \subset [13,16] \times \Psi_4$ we have
  $$\min_X \Theta \geq \min_{v(X)} \Theta -
  \bigg(\frac{|X|^2_1}{512} + |X|^2_2\bigg).$$
\end{lemma}

\startproof
Write $I=[s_0,s_1]$ and $Q=[x_0,x_1] \times [y_0,y_1]$.
Choose a point
$(s,x,y) \in X=I \times Q$.   Taylor's Theorem with remainder
(applied at the point of $[a,b]$ where $f$ is minimized)
implies that for any function $f: [a,b] \to \R$ and any $x \in [a,b]$
we have $$f(x) \geq \min(f(a),f(b))-\frac{1}{8} \max_{[a,b]} |f''|
\times |a-b|^2.$$
Applying this result $3$ times and using Lemma C3
we have
$$\Theta(s,x,y) \geq \min_{i} \Theta(s_i,x,y)- \frac{|I|^2}{512}
\geq 
\min_{i,j} \Theta(s_i,x_j,y) - \frac{|I|^2}{512}- \frac{|X|^2_2}{2}
\geq $$
$$
\min_{i,j,k} \Theta(s_i,x_j,y_k)
- \frac{|I|^2}{512}
- \frac{|X|_2^2}{2}
- \frac{|X|_2^2}{2}
=
\min_{v(X)} \Theta
- \frac{|X|_1^2}{512}
- |X|_2^2.$$
This completes the proof.
\endproof

\noindent
{\bf Verifying Inequalities:\/}
Suppose we want to establish an inequality
like
$$(\frac{a}{b})^{\frac{p}{q}}<\frac{c}{d},$$
where every number involved is a positive integer.
This inequality is true iff
$$b^pc^q-a^pd^q>0.$$
We check this using exact integer arithmetic.
The same idea works with $(>)$
in place of $(<)$.  We call this the
{\it expanding out method\/}.

More generally, we will want to verify inequalities like
\begin{equation}
  \label{verify0}
  \sum_{i=1}^{10} b_i^{-s} -\sum_{i=1}^{10} a_i^{-s/2}>C.
  \end{equation}
where all $a_i$ belong to the set $\{2,3,4\}$, and
$b_i,C,s$ are all rational.
More specifically $s \in [13,15_{+}]$ will be a
dyadic rational and $C$ will be positive.
The expression on the left will
be ${\cal E\/}_s(p)-{\cal E\/}_s(p_0)$ for
various choices of $p$, and the constant $C$ is related to the error term we define below.

Here is how we handle expressions like this.
For each index $i \in \{1,...,10\}$ we produce
rational numbers $A_i$ and $B_i$ such that

\begin{equation}
  \label{verify1}
  A_i>a_i^{-s/2},
  \hskip 30 pt
  B_i<b_i^{-s}.
\end{equation}

We use the expanding out method to check
these inequalities.  We then check that
\begin{equation}
  \label{verify2}
  \sum_{i=1}^{10} B_i- \sum_{i=1}^{10} A_i >C.
  \end{equation}
This last calculation is again done with
integer arithmetic.
Equations \ref{verify1} and \ref{verify2} together
imply Equation \ref{verify0}.
Logically speaking, the way that we produce the rational
$A_i$ and $B_i$ does not matter,
but let us explain how we find them in practice.
For $A_i$ we compute $2^{32} a_i^{-s/2}$ and round the
result up to the nearest integer $N_i$. We then
set $A_i=N_i/2^{32}$.  We produce $B_i$ in a similar way,
except that we round down rather than round up.
When we have verified Equation \ref{verify0}
in this manner we say that we have used the
{\it rational approximation method\/} to
verify Equation \ref{verify0}.
We will only need to make verifications like
this on the order of $20000$ times.
\newline
\newline
{\bf The Grading Step:\/}
We say that a rational number $p/q$ is
{\it dyadic\/} if $q$ is a power of $2$.
We say that a block, defined above,
is {\it dyadic\/} if all coordinates of all the block
vertices are dyadic rationals.

    We perform the following pass/fail evaluation of $X$.

\begin{enumerate}
\item If $I \subset [0,13]$ or $I \subset [15_+,16]$ or
  $Q \cap \Psi_4=\emptyset$,
  we pass $X$ because $X$ is irrelevant to
  the calculation.

\item If $Q \cap \Psi_4\neq\emptyset$ but
  $Q\not\subset\Psi_4$, we fail $X$.

\item If $s_0 \geq 15$ and $Q \subset \Psi_4^{\sharp}$ we pass $X$.

\item If $s_0<13$ and $s_1>13$ we fail $X$ because we don't want to
  make any computations which involve exponents less than $13$.

\item If $X$ has not been passed or failed, we
  try to use the rational approximation method to verify that
  $$\Theta(v)>|X|_1^2/512+|X|_2^2$$
  for each vertex $v$ of $X$.
  If we succeed at this, then we pass $X$.
  Otherwise we fail $X$.
\end{enumerate}

To prove the Critical Theorem II it suffices to
find a partition of $$[0,16] \times [0,1]^2$$ into
blocks which all pass the evaluation.
\newline
\newline
    {\bf Subdivision:\/}
    Let $X=I \times Q$.   Here is the rule we use
    to subdivide $X$:
    If $16|X|_2>|X|_1$ we subdivide $X$ along $Q$ dyadically,
    into $4$ pieces. Otherwise we subdivide $X$ along $I$,
    into two pieces.  This method takes advantage of the
    lopsided form of Lemma C22 and produces a small partition.
    \newline
    \newline
    {\bf Running the Algorithm:\/}
    We perform the following algorithm.
    \begin{enumerate}
\item
  We start with a list $L$ of blocks.  Initially $L$ has the single member
  $[0,16] \times [0,1]^2$.
\item We let $B$ be the last block on $L$.  We grade $B$. If $B$ passes, we
  delete $B$ from $L$. If $L=\emptyset$ then {\bf HALT\/}.
  If $B$ fails, we delete $B$ from $L$ and append to $L$ the subdivision of $B$.
  Then we go back to Step 2.
\end{enumerate}

For the calculation, I used the Java program
in {\bf 4SymmetricConfigs/Main\/}.
      When I run the algorithm, it halts with success
    after $21655$ steps and in about $1$ minute.
    The partition it produces has
    $14502$ blocks.
    
    This establishes the Critical Theorem II
    modulo the proof of Lemma C3.
    In the next section we prove
    Lemma C3 and also some derivative
    bounds needed for the Critical Theorem III.

    \subsection{Critical Theorem III}
    \label{BOUNDC}

    Let $\Theta$ be the function from the previous section.
    
\begin{lemma}[C31]
        On $[13,16] \times \Psi_4$ we have
      $\Theta_{xx}, \Theta_{yy}, \Theta_{xy}>0$.
    \end{lemma}

\startproof
We prove this for $\Theta_{xx}$ and $\Theta_{xy}$.
The case of $\Theta_{yy}$ follows from this and symmetry.
Setting $u=s/2$ we compute
\begin{equation}
\label{maineq}
      {\cal E\/}_s(x,y)=A(x,s)+A(y,s)+2B(x,s)+2B(y,s)+4C(s,x,y),
    \end{equation}
    $$A(x)=a(x)^u, \hskip 30 pt
    B(x)=b(x)^u, \hskip 30 pt C(x)=c(x)^u,$$
    $$
a(x) =
\frac{(1+x^2)^2}{16x^2}  \hskip 15 pt
b(x) =
\frac{1+x^2}{4}  \hskip 15 pt 
c(x,y) =
\frac{(1+x^2)(1+y^2)}{4(x^2+y^2)}
$$
Hence
\begin{equation}
  \Theta_{xx}=A_{xx}+2B_{xx}+4C_{xx}, \hskip 30 pt
  \Theta_{xy}=4C_{xy}.
  \end{equation}
  For each choice of $F=A,B,C$ we have
  {\small
\begin{equation}
\label{monoX}
F_{xx}=u(u-1)f^{u-2}f_x^2 + uf^{u-1}f_{xx}, \hskip 14 pt
C_{xy}=u(u-1) c^{u-2} c_x c_y + u c^{u-1} c_{xy}.
\end{equation}
\/}
Our notation is such that $f=a$ when $F=A$, etc.

We compute
$$a_{xx}=\frac{3+x^4}{8x^4}>0, \hskip 20 pt
b_{xx}=\frac{1}{2}, \hskip 20 pt
c_{xx}=\frac{(1-y^4)(3x^2-y^2)}{2(x^2+y^2)^3} \geq 0.$$
$$c_x=\frac{x(y^4-1)}{2(x^2+y^2)^2} \leq 0, \hskip 10 pt
  c_y=\frac{y(x^4-1)}{2(x^2+y^2)^2} \leq 0, \hskip 10 pt
c_{xy}=\frac{2xy(1+x^2y^2)}{(x^2+y^2)^3}>0.$$
Equation \ref{monoX} combines with all this to prove that
$\Theta_{xx}>0$ and $\Theta_{xy}>0$ on $[13,16] \times \Psi_4$.
\endproof

\noindent
{\bf Proof of Statement 3 of the Critical Theorem III:\/}
We actually prove a broader result. 
Let $s \in [13,16]$.  We know by Lemma C31 that
  all the second partials of $\Theta$ are positive
  at each point of $\Psi_8^{\sharp}$ for this value of $s$.
  But then the restriction of $\Theta$ to
  $\Psi_8^{\sharp}$, at this parameter $s$, is a strictly convex function.
  This shows that for each $s \in [13,16]$ the restriction
  of $\Theta$ to $\Psi_8^{\sharp}$ has a unique minimizer.
  In particular, this is true on the smaller
  interval $(\shin,15_+]$.
  \endproof

\noindent
{\bf Proof of Lemma C3:\/}
We keep the notation from the proof of
Lemma C31.
We first consider $\Theta_{xx}$.
We already know $\Theta_{xx}>0$ on our
domain.
An easy exercise in
calculus shows that $f \in (0,3/5)$ on
$\Psi_4$ for each $f=a,b,c$.
From this bound, we see that
the expression in Equation \ref{monoX} is
decreasing as a function of $u$ for $u \geq 6$.
(Recall that $u=s/2$.)
Hence it suffices to prove that
$4-\Theta_{xx} \geq 0$ on
$\{12\} \times [43/64,1]^2$.

We define $\phi(t)=(43/64)(1-t)+t.$  The file
{\tt LemmaC3.m\/} computes that for $s=12$
the polynomial
  $\Phi={\rm num\/}_+(4-\Theta_{xx} \circ \phi)$
  is WPD and hence non-negative on $[0,1]^2$. Hence
  $4-\Theta_{xx} \geq 0$ when $s=12$ and $(x,y) \in \Psi_4$.
  The same bound for $\Theta_{yy}$ follows from symmetry.

  Now we consider $\Theta_{ss}$.
Let $\psi(s,b)=b^{-s}$.  Let $\beta=(1.3,\sqrt 2,\sqrt 3)$ and also let
  $\gamma=(440,753,4184)$. We first establish the following bound:
\begin{equation}
    \label{boundz}
    0<\max_{b \geq \beta_j} \psi_{ss}(s,b) \leq 1/\gamma_j, \hskip 30 pt j=1,2,3, \hskip 30 pt
    \forall s \geq 13.
  \end{equation}
As a function of $s$, and for $b>1$ fixed,
$\psi_{ss}(s,b)=b^{-s}\log(b)^2$
is decreasing. Hence,
it suffices to prove Equation \ref{boundz} when $s=13$.
Choose $b \geq 1.3$.
The equation $\psi_{ssb}(13,b)=0$
has its unique solution in $[1,\infty)$ at
the value $b=\exp(2/13)<1.3$.  Moreover, the function
$\psi_{ss}(13,b)$ tends to $0$ as $b \to \infty$.
  Hence the restriction of the function
  $b \to \psi_{ss}(13,b)$ to $[b,\infty)$ takes its maximum value at $b$.
  Evaluating at $b=1.3,\sqrt 2, \sqrt 3$ we get Equation \ref{boundz}.

For $x,y \in [43/64,1]$ we easily check the inequalities
$$
a(x)^{-1} \geq 3, \hskip 15 pt
b(x)^{-1} \geq 2, \hskip 15 pt
c(x,y)^{-1} \geq (1.3)^2.
$$
The quantities on the left are the reciprocals of
the functions $a,b,c$ defined in
the proof of Lemma C31; they are the square
distances of the various pairs of points in the corresponding
configuration on $S^2$. From this analysis we conclude that
  the $10$ distances associated to a $5$-point
  configuration parametrized by a point in
  $\Psi_4$ are at least $1.3$, and at least $6$
  of them are at least $\sqrt 2$, and at least
  $2$ of them are at least  $\sqrt 3$.  The same
  obviously holds for the TBP.

Now, $10$ of the $20$ terms comprising $\Theta_{ss}(s,x,y)$ are
positive and $10$ are negative.  Also, for the terms of the
same sign, all $10$ of them are less than $1/440$, and at
least $6$ of them are less than $1/753$, and at least
$2$ of them are less than $1/4184$.  Hence, by Equation \ref{boundz}, we have
the final bound $|\Theta_{ss}| \leq (4/440)+(4/753)+(2/4184)<1/64.$
\endproof

With the proof of Lemma C3, we have finished the proof of
the Critical Theorem II.

All that remains is to prove
Statements 1 and 2 of the Critical Theorem III.  We
first prove the derivative bounds we need for this and
then we give the final argument.  Let

\begin{equation}
  I=\bigg[\frac{55}{64},\frac{56}{64}\bigg].
  \end{equation}

\begin{lemma}[C4]
  $  \Theta_{tts}(15,t,t)<0, \hskip 30 pt \forall t \in I.$
\end{lemma}

The file {\tt LemmaC4.m\/} does the calculations for this proof.
Because the $s$-energy of the TBP does not depend on the
$t$-variable, we have
\begin{equation}
  \Theta_{stt}(15,t,t)=2 A_{stt}|_{s=15}+4 B_{stt}|_{s=15} + 4
  C_{stt}|_{s=15}:=
  \alpha(t)+\beta(t)+\gamma(t).
\end{equation}

We write $f \sim f^*$ if
$$\frac{f}{f^*}= 2^u t^v (1+t^2)^w (2+t^2+t^{-2})^x$$
for exponents $u,v,w,x\in \R$.  In this case, $f$ and $f^*$ have the same sign.
\newline
\newline
{\bf Step 1:\/}
Taking $(u,v,w,x)=(-14,0,11/2,0)$ we have $\beta \sim -\beta^*$,
$$\beta^*(t)= (-2 + 30 \log(2)) + t^2 (-58 + 420 \log(2)) - 
  15 (1+14 t^2) \log(1+t^2).
  $$
  
Noting that
$$
0.69<\log(2)<0.7,
\qquad
\log(1+t^2)<\log(2)<0.7,
$$
we have
$$
\beta^*(t)
>
18.7+231.8t^2-10.5(1+14t^2).
$$
Hence
$$
\beta^*(t)>8.2+84.8t^2>0.
$$
Therefore $\beta<0$ on $I$.
\newline
\newline
{\bf Step 2:\/}
Taking $(u,v,w,x)=(-41/2,-16,12,1/2)$ we have
$\gamma \sim -\gamma^*$. Here

$$
\gamma^*(t)=
(-31 + 360 \log(2))
+
\underline{t^2 (56 - 585 \log(2))}
+
t^4 (-29 + 315 \log(2))+
$$

$$
15 (-8 + 13 t^2 - 7t^4)
\log(2 + t^2 + t^{-2}).
$$
We will explain the meaning of the underlined
term below.

We have $\gamma^*(55/64)>2^4$.
Next we show that
$\gamma^*_t>-2^{10}$ on $I$.
We have
$$
P(t)=-8+13t^2-7t^4, \qquad
q(t)=2+t^2+t^{-2}.
$$

On $I\subset(0,1)$ we have
$P<0$ and $P'>0$ and $q>1$ and $q'<0$.
Hence

$$
\frac{d}{dt}\big(P\log q\big)
=
P'\log q+P\frac{q'}{q}>0.
$$

Also, the $t^4$ term has positive derivative on $I$.
Thus the underlined $t^2$ term is the only term in
$\gamma^*$ having negative derivative on $I$.
This gives us our estimate
$\gamma^*_t>-2^{10}$  on $I$.
Since $I$ has length $2^{-6}$, the function
$\gamma^*$ cannot decrease by more than $2^4$ as we move
from $55/64$ to any other point of $I$.  Hence
$\gamma^*>0$ on $I$, and therefore $\gamma<0$ on $I$.
  \newline
  \newline
{\bf Step 3:\/}
Taking $(u,v,w,x)=(-29,-14,10,3/2)$ we have
$\alpha \sim -\alpha^*$, where
$$
\alpha^*(t)=\gamma^*(t)+\delta^*(t),
\hskip 20 pt
\delta^*(t)=15\log(2)\,(8-13t^2+7t^4).
$$

Setting $u=t^2$, the quadratic
$7u^2-13u+8$
has discriminant
$$
13^2-4\cdot 7\cdot 8=-55<0.
$$
Since its leading coefficient is positive, it is positive
for all $u$. Hence $\delta^*(t)>0$ on $I$.

From Step 2 we already know that $\gamma^*(t)>0$ on $I$.
Therefore $\alpha^*(t)>0$ on $I$.
Hence $\alpha<0$ on $I$.
  \endproof

\begin{lemma}
  \label{C5}
  For any $\xi \in \Psi_8^{\sharp}$ let
  $\Theta(s,\xi)={\cal E\/}_s(\xi)-{\cal E\/}_s(\xi_0)$.
  Then $\Theta_s<0$
    for $s \in [15,15_+]$.
\end{lemma}

\startproof
Let $t_0=55/64$ be the left endpoint of the interval $I$.
We compute that
\begin{equation}
  \Theta_{st}(15,t_0,t_0)<0, \hskip 30 pt
  \Theta_s(15,t_0,t_0)<-2^{-7}.
\end{equation}
The previous lemma now tells us that
\begin{equation}
  \frac{d}{dt} \Theta_{st}(15,t,t)= \Theta_{tts}<0, \hskip 30 pt
  \forall t \in I.
\end{equation}
The last two equations therefore combine to show that
\begin{equation}
  \Theta_{s}(15,t,t)<-2^{-7}.\hskip 30 pt \forall t \in I.
\end{equation}

We also have the bound
$|\Theta_{ss}| \leq 2^{-6}$ on $[13,16] \times \Psi_4$.
Hence
\begin{equation}
  |\Theta_{ss}| \times |15_+-15| \leq 2^{-6} \times \frac{25}{512} < 2^{-7}.
  \end{equation}
Hence $\Theta_s(s,t,t)$ varies by less than
$2^{-7}$ as $s$ ranges in $[15,15_+]$.
Hence $\Theta_s(s,t,t)<0$ for all $s \in [15,15_+]$ and all $t \in I$.
\endproof

\noindent
{\bf Proof of Statements 1 and 2 of the Critical Theorem III:\/}
  By the Critical Theorem II, we have
$\Theta(15,*)>0$ on $\Psi^{\sharp}_8$.
We compute $\Theta(15_{+},x,x)<0$ for $x=445/512 \in [55,56]/64$.
Combining this with Lemma \ref{C5}, we see (by continuity and compactness)
that there exists
a smallest parameter $\shin \in (15,15_+)$ such that
$\Theta(\shin,p^*)=0$ for some $p^* \in \Psi^{\sharp}_8$.
For $s>\shin$, Lemma \ref{C5} now says that
$\Theta(s,p^*)<0$.
This establishes Statements 1 and 2 of the Critical Theorem III.
\endproof

  \newpage

%% file: 8energy.tex
\section{The Energy Theorem}
\label{ENG}

\noindent
{\bf Reading Guide:\/}  This chapter is for
Readers 5 and 6.   For Reader 5, we prove the
Energy Theorem in \S \ref{LemmaEproof}.  For Reader 6,
we use the Energy Theorem in our big computation
in \S \ref{calcthm}.

\subsection{Background Definitions}
\label{BLOCKBACK}

We first give some background definitions and then we
give our main result.
\newline
\newline
\noindent
{\bf Energy Hybrids:\/}
Let $\Gamma_1=-G_1$.
We say that an {\it  energy hybrid\/} is a potential of the form
      \begin{equation}
        \label{functionlist}
      F=c_1 \Gamma_1 + \sum_{k=2}^m c_k G_k, \hskip 20 pt
      G_k(r)=(4-r^2)^k, \hskip 20 pt
      c_1,...,c_m  \in \N \cup \{0\}.
    \end{equation}

    \noindent
    {\bf Remark:\/}
    For applications, we only have one instance where $c_1>0$.
    This is for the potential $F=G_5-25G_1=25\Gamma_1+G_5$.
    In this case, $c_1=25$ and $c_5=1$ and
    all other terms are $0$.  Also, everything works for real
    coefficients but we only care about the integer case.
    \newline
    
    We normalize our avatars so that
    $p_0$ lies on the positive $X$-axis.  In this
    way, $p_{02}=0$ and we just drop this coordinate entirely.
    We then string out the remaining coordinates
    and thereby
   identify an avatar with a point in
    $\R^7=\R \times (\R^2)^3$.
    Thus we think of the potential
    ${\cal E\/}_F$ as a function on $\R^7$.
    \newline
    \newline
      \noindent
      {\bf Blocks:\/}
We define a {\it block\/} to be a
product of the form
\begin{equation}
  \label{BLOCK}
  B=Q_0 \times Q_1 \times Q_2 \times Q_3,
\end{equation}
where $Q_0$ is a segment in $[0,\infty)$  and $Q_1,Q_2,Q_3$ are
squares with sides parallel to the coordinate axes.
We call $Q_0,Q_1,Q_2,Q_3$ the {\it components\/} of the block.

We say that $B$ is $r$-{\it acceptable\/} if, for all $i=1,2,3$, the following is true.
\begin{itemize}
\item $Q_i \subset [-r,r]^2$.
\item $Q_i$ is contained in a closed quadrant of $\R^2$.
\end{itemize}
We say that $B$ is {\it acceptable\/} if it is $(3/2)$-acceptable.

\subsection{The Main Result}

 We let $\cal Q$ denote the set of components
 of acceptable blocks.
  We also let
  $\{\infty\}$ be a member of $\cal Q$.
  We let $v(\{\infty\})=\{\infty\}$. Then
  $Q^{\bullet}$ in this case is $\{(0,0,1)\}$.

 We first define $4$ basic
 measurements we take of members in $\cal Q$.
 \newline
 \newline
 {\bf 0. The Flat Approximation:\/}
 Let $v(Q)$ denote the set of vertices of a square
 or segment $Q$.
    Given $Q \in \cal Q$ we define
\begin{equation}
  Q^{\bullet} = {\rm Convex\ Hull\/}(\Sigma^{-1}(v(Q)))
\end{equation}
$Q^{\bullet}$ is either the point $(0,0,1)$, a chord of $S^2$
or else a convex planar quadrilateral with
vertices in $S^2$ that is inscribed in a circle.  We let
$d_{\bullet}$ be the diameter of $Q^{\bullet}$.  The quantity
$d_{\bullet}^2$ is a rational function of the vertices of $Q$.
      \newline
      \newline
    {\bf 1. The Hull Approximation Constant:\/}
    We think of $Q^{\bullet}$ as the linear
    approximation to
\begin{equation}
  \widehat Q=\Sigma^{-1}(Q).
\end{equation}
The constant we define here turns out to measure the
distance between $\widehat Q$ and $Q^{\bullet}$.
When $Q=\{\infty\}$ we define $\delta(Q)=0$.
Otherwise, let
\begin{equation}
  \chi(D,d)=\frac{d^2}{4D} + \frac{(d^2)^2}{4D^3}.
\end{equation}
This weird function turns out to be an upper bound
to a more geometrically meaningful non-rational function that computes
the distance between a chord of length $d$ of a circle
of diameter $D$ and the arc of the circle it subtends.

When $Q$ is a  segment we define
\begin{equation}
  \label{hull1}
  \delta(Q)=\chi(2,\|\widehat q_1-\widehat q_2\|).
\end{equation}
Here $q_1,q_2$ are the endpoints of $Q$.
When $Q$ is a square we define
\begin{equation}
  \label{hull2}
  \delta(Q)=\max(s_1,s_3)+\max(s_2,s_4), \hskip 30 pt
  s_j=\chi(1,\|\widehat q_j - \widehat q_{j+1}\|).
\end{equation}
Here $q_1,q_2,q_3,q_4$ are the vertices of $Q$ and the
indices are taken cyclically.
These are rational computations because $\chi(1,d)$ and $\chi(2,d)$
are polynomial in $d^2$.
\newline
\newline
    {\bf 2. The Dot Product Estimator:\/}
    By way of motivation, we
        point out that if $V_1,V_2 \in S^2$ then
    $$G_k(\|V_1-V_2\|) = (2+2 V_1 \cdot V_2)^k.$$

    Now suppose that $Q_1,Q_2 \in \cal Q$.
We set $\delta_j=\delta(Q_j)$.
Given any $p \in \R^2 \cup \infty$ let
$\widehat p=\Sigma^{-1}(p)$.   Define
\begin{equation}
Q_1\cdot Q_2=\max_{i,j}(\widehat q_{1i} \cdot \widehat q_{2j})
+(\tau) \times (\delta_1+\delta_2+\delta_1\delta_2).
\end{equation}
Here $\{q_{1i}\}$ and $\{q_{2j}\}$ respectively are the vertices of
$Q_1$ and $Q_2$.  The constant $\tau$ is $0$ if one of
$Q_1$ or $Q_2$ is $\{\infty\}$ and otherwise $\tau=1$.
Finally, we define
\begin{equation}
  T(Q_1,Q_2)=2+2(Q_1 \cdot Q_2).
\end{equation}
\noindent
{\bf 3. The Local Error Term:\/}
For $Q_1,Q_2 \in \cal Q$ and $k \geq 2$ we define
\begin{equation}
\label{EPSILON}
\epsilon_k(Q_1,Q_2)=
\frac{1}{2} k(k-1) T^{k-2}d_1^2 +
2kT^{k-1} \delta_1,\end{equation}
$$d_1=d_{\bullet}(Q_1), \hskip 15 pt
\delta_1=\delta(Q_1), \hskip 15 pt
T=T(Q_1,Q_2).
$$
When $k=1$ we define
\begin{equation}
  \label{EPSILON2}
  \epsilon_1(Q_1,Q_2)=2\delta_1.
\end{equation}

We mention two features of this quantity:
\begin{itemize}
\item   In general, $\epsilon_k(Q_1,Q_2) \not = \epsilon_k(Q_2,Q_1)$.
\item When $Q_1=\{\infty\}$ we have $d_1=\delta_1=0$.  So,
  $\epsilon_k(\{\infty\},Q_2)=0$.
\end{itemize}

We will explore the geometry of the expression
$\epsilon_k(Q_1,Q_2)$ in the next chapter.  
Here is a brief geometric explanation of this quantity.
The first term on the right in Equation \ref{EPSILON} comes
from the analysis of the flat approximation and
the second term comes from the analysis
of the difference between the flat approximation
and the actual subset of the sphere.

\newpage

\noindent
{\bf 4. The Global Error Estimate:\/}
We always set $Q_4=\{\infty\}$.
Given a block $B=Q_0 \times Q_1 \times Q_2 \times Q_3$ let
\begin{equation}
\label{errorsum}
      {\bf ERR\/}_k(B)=\sum_{i=0}^3 {\bf ERR\/}_k(B,i), \hskip 30 pt
{\bf ERR\/}_k(B,i)=\sum_{j \not =i}  \epsilon_k(Q_i,Q_j).
\end{equation}
In the last sum, $j \in \{0,1,2,3,4\}-\{i\}$.
More generally, when $$F=c_1\Gamma_1+\sum_{k=2}^m c_k G_k$$ is as in Equation
\ref{functionlist}, we define
  \begin{equation}
  {\bf ERR\/}_F(B)=\sum_{i=0}^3 {\bf ERR\/}_F(B,i), \hskip 20 pt
  {\bf ERR\/}_F(B,i)=\sum_{k=1}^m c_k\ {\bf ERR\/}_k(B,i)
\end{equation}

We prove the following result in the next chapter.

\begin{theorem} [Energy]
Let $B$ be an acceptable block.  Let $F$ be an energy
hybrid.  Then
\label{ENERGYX}

$$
\min_{p\in B} {\cal E\/}_F(p) \geq
\min_{p\in v(B)} {\cal E\/}_F(p)-{\bf ERR\/}_F(B).
$$
\end{theorem}

\noindent
{\bf Remark:\/}
Recall that acceptable blocks are
$3/2$-acceptable.   Referring to Equation
\ref{hull2}, our proof of the Energy Theorem
for $r$-acceptable blocks works with
$\chi(f(r),d)$ in place of $\chi(1,d)$
where $f(r)=4/(1+r^2)$.    For $r=3/2$
we could have used $16/13$ in place
of $1$ but we wanted to keep it simple.

\newpage

%% file: 9energy_proof.tex
\section{Proof of the Energy Theorem}
\label{LemmaEproof}

\noindent
{\bf Reading Guide:\/}
This chapter is for Reader 5.

\subsection{Guide to the Proof}

Our proof of the Energy Theorem splits into two halves,
an algebraic part and a geometric part.
The geometric part establishes
a local result in terms of the Local Error Term
from \S 8.

Now we discuss the algebraic part of the proof.
We say that an {\it averaging system\/}
for a member of $\cal Q$ is a collection
of maps
$\lambda_1,\lambda_2,\lambda_3,\lambda_4: Q \to [0,1]$
such that
$$\sum_{i=1}^4 \lambda_i(z)=1, \hskip 30 pt \forall\ z \in Q.$$
This is analogous to a partition of unity, except that
the functions need not vary continuously.
In case $Q$ is a segment, we list the vertices
as $q_1,q_1,q_2,q_2$ and we have
$\lambda_1=\lambda_2$ and $\lambda_3=\lambda_4$.
In case $Q=\{\infty\}$ we list the vertices
as $\infty,\infty,\infty,\infty$ and we have
$\lambda_j=1/4$ for $j=1,2,3,4$.

We say that an {\it averaging system\/}
for $\cal Q$ is a choice of averaging system
for each member $Q$ of $\cal Q$.  The
averaging systems for different members
need not have anything
to do with each other.
We will construct a particular averaging system for $\cal Q$
which compares the energy of an arbitrary
configuration with the energies of the vertex configurations.
The algebraic part of the proof promotes a local estimate
about the averaging system to the
Energy Theorem.  Our averaging system is a
geometric construction.  It works the same for all
energy potentials.

There is an important algebraic point in our proof.
The Energy Theorem is stated in terms of minima,
but we get a stronger result in terms of averaging.
Because the averaging works the same way for
all potentials, it will suffice to
prove the basic estimates just
for the basic energy potentials
$\Gamma_1,G_2,G_3,...$ and add  the resulting
estimates.
We minimize {\it after\/} averaging.

\subsection{Reduction to a Local Result}

We fix the function $F=G_k$ for some $k \geq 2$ or else $F=\Gamma_1=-G_1$.
We write ${\cal E\/}={\cal E\/}_F$.
When $F=\Gamma_1$ we set $k=1$.
We let $\epsilon=\epsilon_k$, as
in Equation \ref{EPSILON}.
Our algebraic argument would work for
any choice of $F$, but we need to use the choices above
to actually get the averaging system
we need.  Let
$q_{1,1},q_{1,2},q_{1,3},q_{1,4}$ be the vertices of $Q_1$.
Again, when $Q_1$ is a square, these points go in cyclic order around the square.

\begin{lemma}[E1]
  There exists an averaging system on
  $\cal Q$ with the following property:
Let $Q_1,Q_2$ be members of $\cal Q$.
Given any $z_1 \in Q_1$ and $z_2 \in Q_2$ we have
\begin{equation}
  \sum_{i=1}^4 \lambda_i(z_1)F(\|\widehat q_{1,i}-\widehat z_2\|) -
  F(\|\widehat z_1-\widehat z_2\|) \leq \epsilon(Q_1,Q_2).
\end{equation}
\end{lemma}
What we are doing is estimating the cost of replacing
a calculation involving $z_1$ with a calculation involving
the vertices of $Q_1$.  We prove this result below.

We are interested in $5$-point configurations but we will
work more generally so as to elucidate the general structure of the argument.
We suppose that we have the acceptable block
$B=Q_0 \times ... \times Q_N$.   For our application,
we take $N=4$ and $Q_4=\{\infty\}$.
The vertices of $B$ are indexed by
a multi-index $$I=(i_0,...,i_N) \in \{1,2,3,4\}^{N+1}.$$
Given such a multi-index, which amounts to a choice
of vertex in each component member of the block,
we define (as always, {\it via\/} inverse stereographic projection)
the energy of the corresponding vertex configuration:

\begin{equation}
{\cal E\/}(I)=
{\cal E\/}(q_{0,i_0},...,q_{N,i_N})
\end{equation}

Here is one more piece of notation. Given
$z=(z_0,...,z_N) \in B$ and a multi-index $I$
we define
\begin{equation}
  \lambda_I(z)=\prod_{j=0}^N \lambda_{i_j}(z_j).
\end{equation}
Here $\lambda_{i_j}$ is defined relative to the
averaging system on $Q_j$.

Now we state our main global result.
The global result uses the existence of an efficient
averaging system.  That is, it relies on the Lemma E1.

\begin{lemma}[E2]
\label{ENERGY2}
Let $z=(z_0,...,z_N) \in B$.
Then
\begin{equation}
  \label{AVE}
  \sum_{I} \lambda_I(z) {\cal E\/}(I)-
      {\cal E\/}(z)
\leq \sum_{i\not = j} \epsilon(Q_i,Q_j).
\end{equation}
The lefthand sum is taken over all multi-indices.
\end{lemma}

We have proved Lemma E2 relative to one of
the potentials $$F=\Gamma_1,G_2,G_3,...$$
However, we use one and the same averaging
system for all of these potentials.
The inequality in Lemma E2 is additive:
Multiplying the individual inequalities by the
positive constants and then adding them up,
we see that Lemma E2 is true for any energy hybrid.

Now let us deduce the Energy Theorem from Lemma E2.
Our quantities in this deduction are defined relative
to some arbitrary energy hybrid.
Notice that
\begin{equation}
  \sum_{I} \lambda_I(z)=
\prod_{j=0}^N \bigg(\sum_{a=1}^4 \lambda_a(z_j)\bigg)=1.
\end{equation}

Choose some $z=(z_0,...,z_N) \in B$
which minimizes ${\cal E\/}$.  We have
$$
0 \leq \min_{p \in v(B)} {\cal E\/}(p)-\min_{p \in B}
{\cal E\/}(p)=\min_{p \in v(B)} {\cal E\/}(p)-{\cal E\/}(z) \leq^* $$
\begin{equation}
\sum_{I} \lambda_I(z){\cal E\/}(I)-
{\cal E\/}(z) \leq  \sum_{i\not = j}\epsilon(Q_i,Q_j).
\end{equation}
Here, writing $F=c_1\Gamma_1+\sum_{k =2}^m c_k G_k$ as
in Equation \ref{functionlist}, we have set
$$\epsilon(Q_i,Q_j)=\sum_{k=1}^m c_k \epsilon_k(Q_i,Q_j).$$
The starred inequality comes from the fact that a minimum is less than or
equal to a convex average.

Setting $N=4$ and noting that
$\epsilon(Q_4,Q_i)=0$ for all $i$, we
see that this last expression is
${\bf ERR\/}_F(B)$.    Thus, the Energy Theorem
stated at the end of the last chapter follows
from Lemma E2.

\subsection{Proof of Lemma E2}

Now we deduce the global Lemma E2 from the local Lemma E1.
First suppose $N=1$.
In this case,
we have a block $B=Q_0 \times Q_1$.
Setting $\epsilon_{ij}=\epsilon(Q_i,Q_j)$, Lemma
E1 gives us
\begin{equation}
  F(\|\widehat z_0-\widehat z_1\|) \geq
\sum_{\alpha=1}^4 \lambda_{\alpha}(z_0)
F(\|\widehat q_{0\alpha}-\widehat z_1\|)-\epsilon_{01}.
\end{equation}
Applying Lemma E1 to the pair of points
$(z_1,q_{0\alpha}) \in Q_1 \times Q_0$ we have
\begin{equation}
F(\|\widehat z_1-\widehat q_{0\alpha}\|) \geq
\sum_{\beta=1}^4 \lambda_{\beta}(z_1)
F(\|\widehat q_{1\beta}-\widehat q_{0\alpha}\|)
-\epsilon_{10}.
\end{equation}
Plugging the second equation into the first and using
$\sum \lambda_{\alpha}(z_0)=1$, we have
\[
F(\|\widehat z_0-\widehat z_1\|) \geq -\epsilon_{01}
+\sum_{\alpha,\beta}
\lambda_{\alpha}(z_0)\lambda_{\beta}(z_1)
\left[
F(\|\widehat q_{1\beta}-\widehat q_{0\alpha}\|)
-\epsilon_{10}
\right]=
\]
\begin{equation}
  \label{AVE1}
 \sum_{\alpha,\beta}
 \lambda_{\alpha}(z_0)\lambda_{\beta}(z_1)
 F(\|\widehat q_{1\beta}-\widehat q_{0\alpha}\|)
 -(\epsilon_{10}+\epsilon_{01}).
\end{equation}
Equation \ref{AVE1} is equivalent
to Equation \ref{AVE} when $N=1$.
\newline

Now we do the general case and suppose that $N \geq 2$.
We rewrite Equation \ref{AVE1} as follows:
\begin{equation}
\label{AVE2}
F(\|\widehat z_0-\widehat z_1\|) \geq
\sum_{A}
\lambda_{A_0}(z_0)\lambda_{A_1}(z_1)\
F(\|\widehat q_{0A_0}-\widehat q_{1A_1}\|)
-(\epsilon_{01}+\epsilon_{10}).
\end{equation}
The sum is taken over multi-indices $A$ of length $2$.

We also observe that
\begin{equation}
\label{AVE3}
\sum_{I'} \lambda_{I'}(z')=1, \hskip 30 pt
z'=(z_2,...,z_N).
\end{equation}
The sum is taken over all multi-indices
$I'=(i_2,...,i_N)$.
Therefore, if we hold $A=(A_0,A_1)$ fixed, we have
\begin{equation}
\lambda_{A_0}(z_0)\lambda_{A_1}(z_1) =
\sum_{I''} \lambda_{I''}(z).
\end{equation}
The sum is taken over all multi-indices of length $N+1$ which
have $I_0=A_0$ and $I_1=A_1$.
Combining these equations, we have
\begin{equation}
F(\|\widehat z_0-\widehat z_1\|) \geq
\sum_I \lambda_{I}(z)
F(\|\widehat q_{0I_0}-\widehat q_{1I_1}\|)
-(\epsilon_{01}+\epsilon_{10}).
\end{equation}
The same argument works for other pairs of indices, giving
\begin{equation}
\label{AVE4}
F(\|\widehat z_i-\widehat z_j\|) \geq
\sum_I \lambda_I(z)
F(\|\widehat q_{iI_i}-\widehat q_{jI_j}\|)
-(\epsilon_{ij}+\epsilon_{ji}).
\end{equation}

Let us restate this as
$X_{ij}-Y_{ij}\leq Z_{ij},$
where
\[
X_{ij}=
\sum_I \lambda_I(z)
F(\|\widehat q_{iI_i}-\widehat q_{jI_j}\|),
\hskip 10 pt
Y_{ij}=F(\|\widehat z_i-\widehat z_j\|),
\hskip 10 pt
Z_{ij}=\epsilon_{ij}+\epsilon_{ji}.
\]
When we sum $Y_{ij}$ over all $i<j$ we get the second term in Equation
\ref{AVE}.
When we sum $Z_{ij}$ over all $i<j$ we get the third term in Equation
\ref{AVE}.
When we sum $X_{ij}$ over all $i<j$ we get
\[
\sum_{i<j}
\bigg(
\sum_I \lambda_I(z)
F(\|\widehat q_{iI_i}-\widehat q_{jI_j}\|)
\bigg)
=
\sum_I \sum_{i<j}
\lambda_I(z)
F(\|\widehat q_{iI_i}-\widehat q_{jI_j}\|)
=
\]
\[
\sum_I \lambda_I(z)
\Bigg(
\sum_{i<j}
F(\|\widehat q_{iI_i}-\widehat q_{jI_j}\|)
\Bigg)
=
\sum_I \lambda_I(z){\cal E\/}(I).
\]
This is the first term in Equation \ref{AVE}.
This proves Lemma E2 modulo the proof of
Lemma E1.

\subsection{The Efficient Averaging System}

The rest of the chapter is devoted to proving Lemma E1.
Lemma E1 posits the existence of what
we call an efficient averaging system.
Here we define it.  Recall that $Q^{\bullet}$ is the
convex hull of the vertices
$\widehat q_1,\widehat q_2,\widehat q_3,\widehat q_4$ of
$\widehat Q=\Sigma^{-1}(Q)$.  In the square case,
these vertices go in cyclic order around $Q^{\bullet}$.

What we want from the
system is that for any $z^{\bullet} \in Q^{\bullet}$
\begin{equation}
  z^{\bullet}=\sum_{i=1}^4 \lambda_i(z^{\bullet}) \widehat q_i.
\end{equation}
Here is how we arrange this.
\begin{itemize}
\item If $z^{\bullet}$ lies in the
convex hull of
$\widehat q_1$, $\widehat q_2$, $\widehat q_3$,
then we let $\lambda_1(z^{\bullet})$, $\lambda_2(z^{\bullet})$, $\lambda_3(z^{\bullet})$
be barycentric coordinates on this triangle and we set
$\lambda_4(z^{\bullet})=0$.
\item  If $z^{\bullet}$ lies in the
convex hull of
$\widehat q_1$, $\widehat q_3$, $\widehat q_4$,
then we let $\lambda_1(z^{\bullet})$, $\lambda_3(z^{\bullet})$, $\lambda_4(z^{\bullet})$
be barycentric coordinates on this triangle and we set
$\lambda_2(z^{\bullet})=0$.
\end{itemize}
This definition agrees on the overlap, which is the line segment
joining $\widehat q_1$ to $\widehat q_3$.
Our definition favors the diagonal $q_1q_3$ over
the diagonal $q_2q_4$ but this does not matter.

To get our averaging system on $Q \in \cal Q$ we do
the following.  When $Q$ is a square, we define
\begin{equation}
  \label{CONVEX}
  \lambda_j(z)=\lambda_j(z^{\bullet}),
\end{equation}
where $z^{\bullet}$ is some choice of
point in $Q^{\bullet}$ which is \underline{closest}
to $\widehat z$.  If there are several closest points
we pick the one (say) which has the smallest first coordinate.

When $Q$ is a segment, listed as
$q_1,q_1,q_2,q_2$ we let
$z^{\bullet}$ be some point in $Q^{\bullet}$ closest to $\widehat z$.
There is some $t \in [0,1]$ such that
$z^\bullet=t\widehat q_1+(1-t)\widehat q_2$.
We set
$$
\lambda_1(z)=\lambda_2(z)=\frac{t}{2},
\hskip 30 pt
\lambda_3(z)=\lambda_4(z)=\frac{1-t}{2}.
$$
If $Q=\{\infty\}$, we set $\lambda_i=1/4$ for $i=1,2,3,4$.

We prove Lemma E1 with respect to the averaging system above.
Let $F$ be either $G_k$ for some $k \geq 2$ or else $F=\Gamma_1=-G_1$.
In this case we set $k=1$.
For convenience we expand out the statement of
Lemma E1.

\begin{lemma}[E1]
  The efficient averaging system on
  $\cal Q$ has the following property.
Let $Q_1,Q_2$ be members of $\cal Q$.
Given any $z_1 \in Q_1$ and $z_2 \in Q_2$ we have
\begin{equation}
  \sum_{i=1}^4 \lambda_i(z_1)F(\|\widehat q_{1,i}-\widehat z_2\|) -
  F(\|\widehat z_1-\widehat z_2\|) \leq
\frac{1}{2} k(k-1) T^{k-2}d_1^2 +
2kT^{k-1} \delta_1.
\end{equation}
\end{lemma}
Here
$\delta_1$ and $d_1$ respectively are the Hull Approximation constant and diameter of $Q_1$, and
\begin{equation}
  \label{dotdef}
  T=2+2(Q_1 \cdot Q_2), \hskip 20 pt
Q_1\cdot Q_2=\max_{i,j}(\widehat q_{1,i} \cdot \widehat q_{2,j})
+(\tau) \times (\delta_1+\delta_2+\delta_1\delta_2).
\end{equation}
$\tau=0$ or $\tau=1$ depending on whether one of $Q_1,Q_2$ is $\{\infty\}$.
We are maximizing over the dot product of the vertices and then either
adding an error term or not.
When $k=1$ the expression in Lemma E1 massively simplifies.
We address this in the remark below.
Define

\begin{equation}
  X_{\bullet}=
  (2+2 z_1^{\bullet} \cdot \widehat z_2)^k \hskip 10 pt {\rm or\/} \hskip 10 pt
  -2-2 z_1^{\bullet} \cdot \widehat z_2.
\end{equation}

Lemma E1 is an immediate consequence of the following two results.

\begin{lemma}[E11]
$\sum_{i=1}^4 \lambda_i(z_1)F(\|\widehat q_{1,i}-\widehat z_2\|) -
X_{\bullet} \leq \frac{1}{2}k(k-1)T^{k-2}d_1^2.$
\end{lemma}

\begin{lemma}[E12]
$X_{\bullet}-F(\|\widehat z_1- \widehat z_2\|) \leq 2k T^{k-1}\delta_1.$
\end{lemma}

\noindent
{\bf Remark:\/}
When $k=1$ the expression on the right in Lemma E11 and the corresponding
expression in Lemma E1 are meant to be $0$.  In this case
we will justify this by showing that the expression on the
left side in Lemma E11 is $0$.
  \newline

  The rest of the chapter is devoted to proving
  Lemmas E11 and E12.  We first prove two helpful
  analytic estimates A1 and A2.  Then we prove
  two helpful geometric estimates G1 and G2.  Then we
  prove Lemmas E11 and E12.

  \subsection{Analytic Estimates}

 We first prove two purely
  analytic lemmas. The first one explains the first term
  in our local error expression. The second one deals with the
  second term.
  
\begin{lemma}[A1]
Let $M \geq 2$ and $k=2,3,...$.
Suppose \begin{itemize}
\item $x_1,...,x_M \geq 0$.
\item $\sum_{i=1}^M \lambda_i=1$ and $\lambda_i \geq 0$ for all $i$.
\end{itemize}
Then
\begin{equation}
\label{INEQ}
0 \leq \sum_{i=1}^M \lambda_i x_i^k-
\bigg( \sum_{i=1}^M \lambda_i x_i\bigg)^k \leq 
\frac{1}{8} k(k-1) x_+^{k-2}\ (x_+-x_-)^2.
\end{equation}
Here $x_+=\max(x_i)$ and $x_-=\min(x_i)$.
\end{lemma}

\startproof
Without loss of generality, we
may relabel the pairs $(x_i,\lambda_i)$ so that
$x_1 \leq \cdots \leq x_M$.
This gives $x_-=x_1$ and $x_+=x_M$.

The lower bound in Lemma A1 is a trivial consequence of
convexity, and both bounds are trivial when
$k=1$.  So, we take $k=2,3,4,...$ and prove
the upper bound.
Suppose first that $M \geq 3$.
We have at least one degree of freedom when we keep 
$\sum \lambda_i x_i$ constant, fix the $x$-values,
and try to vary
$\{\lambda_j\}$ so as to maximize the left
hand side of Equation \ref{INEQ}.  The right hand
side does not change when we do this, and 
the left hand side varies linearly. Hence,
the left hand side is maximized when
$\lambda_i=0$ for some $i$. But then
any counterexample
to the lemma for $M \geq 3$ gives rise to
a counterexample for $M-1$.
(The point here is that the deletion of the
corresponding $x_i$ does not change the
left hand side and cannot increase the right hand side.)
Hence, it suffices to prove the inequality
when $M=2$.  

In the case $M=2$, we set
$a=\lambda_1$.  If $x_2=0$ then
$x_1=x_2=0$ and the result is trivial.
So, assume $x_2>0$.
Both sides of the inequality in
Lemma A1 are homogeneous of
degree $k$, so it suffices to consider the
case when $x_2=1$.  We set $x=x_1 \in [0,1)$.
Equation \ref{INEQ} becomes the statement
that $f(x) \leq g(x)$ where
\begin{equation}
\label{McM}
f(x)=
(a x^k + 1-a)-(a x + 1-a)^k;
\hskip 20 pt
g(x)=\frac{1}{8}k(k-1)(1-x)^2.
\end{equation}

Since $f(1)=g(1)=f'(1)=g'(1)=0$ the
Cauchy Mean Value Theorem (applied twice)
tells us that for any $x \in [0,1)$ there are values
$y<z \in [x,1]$ such that
  \begin{equation}
    \label{ref00}
    \frac{f(x)}{g(x)}=\frac{f'(y)}{g'(y)}=
\frac{f''(z)}{g''(z)}=
4az^{k-2}\bigg[1-a\bigg(a+\frac{1-a}{z}\bigg)^{k-2}\bigg] \leq 4a(1-a) \leq 1.
\end{equation}
This completes the proof.
\endproof

\noindent
{\bf Remark:\/} An anonymous referee suggested
the proof that $f(x) \leq g(x)$. This proof is better
than my original proof.
\newline

  We remind the reader of the weird function $\chi(D,d)$ and we
  introduce
  a more geometrically meaningful function
  \begin{equation}
    \chi(D,d)=\frac{d^2}{4D} + \frac{d^4}{4D^3}, \hskip 30 pt
        \chi^*(D,d)=\frac{1}{2}(D-\sqrt{D^2-d^2}).
      \end{equation}

      Let us compare these functions.

      \begin{lemma}[A2]
        $\chi^*(D,d) \leq \chi(D,d)$ for all $d \in [0,D]$.
\end{lemma}

\startproof
  By homogeneity, it suffices
to prove the result when $D=1$.  To simplify the algebra we define
$A=2 \chi(1,d)-1$ and $A^*=2\chi^*(1,d)-1$.
We compute
$$4A^2-4(A^*)^2=d^4(d-1)(d+1)(d^2+3).$$
Hence, the sign of $A-A^*$ does not change on
$(0,1)$.  We check that $A>A^*$ when $d=1/2$.
Hence $A>A^*$ on $(0,1)$.  This implies
the inequality.
\endproof

\subsection{Geometric Estimates}

Now we prove two geometric estimates.
Let $\delta(Q)$ be the hull approximation constant
for $Q \in \cal Q$, as defined
in Equation \ref{hull1} or Equation \ref{hull2}.

\begin{lemma}[G1]
  Let $Q \in \cal Q$ be any component of an acceptable block.
  Then every point of $\widehat Q$ is within
  $\delta(Q)$ of  $Q^{\bullet}$.
\end{lemma}

\startproof
We first consider the case when $Q \subset [0,\infty)$ is a segment.
Here $\widehat Q$ is the arc of a great circle
  and $Q^{\bullet}$ is the chord of the arc
  joining the endpoints of this arc.
  Let $d$ be the length of $Q^{\bullet}$.
  The point of $\widehat Q$ farthest from
  $Q^{\bullet}$ is the midpoint of this $\widehat Q$.
  Let $x$ be the distance between the midpoint
  of $\widehat Q$ and the midpoint of
  $Q^{\bullet}$.  From elementary geometry,
  $x(2-x)=(d/2)^2$.
    Solving for $x$ we find that $x=\chi^*(2,d)$.
    Lemma A2 finishes the proof.

    Now we consider the case where $Q$ is a square.
    This is the only place in the proof of the Energy Theorem where we
    use  the condition that $Q \in [-3/2,3/2]^2$.   Recall also that
    $Q$ is contained
  in a closed quadrant.  Let $z \in Q$ be an arbitrary point.
  We label the vertices of $Q$
    such that $q_1q_2$ and $q_3q_4$ are the
    horizontal segments.
Let $L$ be the vertical line through $z$ and let
$z_{12},z_{34}$  respectively be points where
$L$ intersects $q_1q_2$ and $q_3q_4$.

  If $M$ is a horizontal line intersecting $Q$ then
  the circle $\Sigma^{-1}(M \cup \infty)$ has diameter at least $1$.
  The point is that this circle contains $(0,0,1)$ and also
  $\Sigma^{-1}(0,y)$ for some $|y| \leq 3/2$.  In fact the diameter is
  at least $4/\sqrt{13}$.  The same goes for vertical lines
  intersecting $Q$.

Define
$d_j=\|\widehat q_j-\widehat q_{j+1}\|$ with the indices taken
cyclically.  Let $L$ be the vertical line whose
$X$-coordinate equals $c$.
Then $L\cap Q$ has endpoints $(c,a)$ and $(c,b)$.
Let $\sigma$ be the segment joining
$\Sigma^{-1}(c,a)$ to $\Sigma^{-1}(c,b)$.
Then
\[
 \|\Sigma^{-1}(c,a)-\Sigma^{-1}(c,b)\|
 =
 \frac{2|a-b|}
 {\sqrt{(1+c^2+a^2)(1+c^2+b^2)}}.
\]
 Since $Q$ lies in a quadrant,
the quantity above is maximized when $L$ is one of
the two vertical sides of $Q$.
Hence $\sigma$ has length at most $\max(d_2,d_4)$.
At the same time, $\Sigma^{-1}(L \cap Q)$ is
contained in a circle of diameter at least $1$.
The same argument as in the segment case now shows that
there is a point $z^* \in \sigma$ which is within
$t_{24}=\max(\chi(1,d_2),\chi(1,d_4))$
of $\widehat z$.  Here we have used the monotonicity
of $\chi(D,d)$ as a function of $D$, and the fact
that $D \geq 1$.

The endpoints of $\sigma$ respectively are on the
spherical arcs obtained by mapping the top and
bottom edge of $Q$ onto $S^2$ via $\Sigma^{-1}$.
By the same reasoning as in the Segment case,
one endpoint of $\sigma$ is within
$\chi(1,d_1)$ of a point on the corresponding
edge of $\partial Q^{\bullet}$
and the other endpoint of $\sigma$ is within
$\chi(1,d_3)$ of a point on the opposite
edge of $\partial Q^{\bullet}$.
But that means that each endpoint of $\sigma$ is
within
$t_{13}=\max(\chi(1,d_1),\chi(1,d_3))$
of a point in $Q^{\bullet}$.
But then every point of the segment $\sigma$ is within $t_{13}$ of
some point of the line segment joining these
two points of $Q^{\bullet}$.
In particular, there is a point $z^{\bullet} \in Q^{\bullet}$
which is within $t_{13}$ of $z^*$.  The triangle
inequality completes the proof.
\endproof

\begin{lemma}[G2]
  For $j=1,2$ let $\gamma_j$ be a point on a line segment connecting a point of
  $\widehat Q_j$ to a closest point on $Q_j^{\bullet}$.
  Then
  $\gamma_1 \cdot \gamma_2  \leq Q_1 \cdot Q_2$.
\end{lemma}

\startproof
See Equation \ref{dotdef}  for the definition of $Q_1 \cdot Q_2$.
We first treat the case $\tau=1$, meaning that neither $Q_1$ nor $Q_2$ is $\{\infty\}$.
Since the dot product is bilinear,
for any
$q_1^\bullet\in Q_1^\bullet$ and
$q_2^\bullet\in Q_2^\bullet$ we have
\begin{equation}
  \label{max1}
  q_1^{\bullet} \cdot q_2^{\bullet} \leq
  \max_{i,j}(\widehat q_{1,i} \cdot \widehat q_{2,j}).
\end{equation}
By Lemma G1, and by hypothesis, we can find points
$z_1^{\bullet}$ and 
$z_2^{\bullet}$ such that
$$\gamma_1 = z_1^{\bullet} + h_1, \hskip 30 pt
\gamma_2 = z_2^{\bullet} + h_2, \hskip 30 pt
\|h_j\| \leq \delta_j.$$

Since $z_j^\bullet$ lies in the convex hull of points of $S^2$,
we have $\|z_j^\bullet\|\leq 1$.
But then by the triangle inequality and the Cauchy-Schwarz inequality
$$
|(\gamma_1 \cdot \gamma_2)-(z_1^{\bullet} \cdot z_2^{\bullet})| \leq
|z_1^{\bullet} \cdot h_2| +
|z_2^{\bullet} \cdot h_1| + |h_1 \cdot h_2| \leq \delta_1+\delta_2 + \delta_1 \delta_2.$$
This combines with Equation \ref{max1} to complete the proof when $\tau=1$.

Now we consider the case $\tau=0$.
Without loss of generality assume that $Q_2=\{\infty\}$.
If $q_1=(x,y)$ then, by Equation \ref{inversestereo},
\[
\widehat q_1\cdot (0,0,1)
=
1-\frac{2}{1+x^2+y^2}.
\]
This quantity increases with $x^2+y^2$.
Since $Q_1$ lies in a quadrant, $x^2+y^2$ is maximized
at a vertex of $Q_1$.
At the same time, the maximum of $q_1^{\bullet} \cdot (0,0,1)$,
for $q_1^{\bullet} \in Q_1^{\bullet}$ is achieved when $q_1^{\bullet}$ is a
vertex of $Q_1^{\bullet}$.  But then our lemma is true for the endpoints
of the segment containing $\gamma_1$.  Since the dot product with $(0,0,1)$
varies linearly along this line segment, the same result is true for
all points on the line segment.
\endproof

\subsection{Lemmas E11 and E12}

\noindent
{\bf Proof of Lemma E11:\/}
Suppose first $F=-G_1$.  We hold $\widehat z_2$ fixed.
For $v \in \R^3$ define
\begin{equation}
L(v)=-2-2v\cdot\widehat z_2.
\end{equation}
For $\widehat q\in S^2$ we have
$L(\widehat q)=F(\|\widehat q-\widehat z_2\|)$.
Lemma E11, in this special case, says that
\[
\sum_{i=1}^4 \lambda_i(z_1)L(\widehat q_{1,i})
-L(z_1^\bullet)=0.
\]
But this follows from Equation \ref{CONVEX} and the (bi) linearity of the
dot product.

Now we deal with the case where $F=G_k$ for $k \geq 2$.
Recall that $q_{1,1},q_{1,2},q_{1,3},q_{1,4}$ are the vertices of
$Q_1$, listed in cyclic order.
Recall also that $z_j$ is an arbitrary point in $Q_j$
for $j=1,2$.
Let $\lambda_i=\lambda_i(z_1)$.  
We set
\begin{equation}
x_i=4-\|\widehat q_{1,i}-\widehat z_2\|^2=2+2\widehat q_{1,i} \cdot \widehat z_2, \hskip 30 pt i=1,2,3,4.
\end{equation}
Note that $x_i \geq 0$ for all $i$.

We have
\begin{equation}
\label{term1}
\sum_{i=1}^4 \lambda_i(z_1)F(\|\widehat q_{1,i}-\widehat z_2\|)=\sum_{i=1}^4 \lambda_i x_i^k,
\end{equation}
\begin{equation}
\label{term2}
X_{\bullet}=(2+2 z_1^{\bullet} \cdot \widehat z_2)^k=
\bigg(\sum_{i=1}^4 \lambda_i \times (2+2 \widehat q_{1,i} \cdot \widehat z_2)\bigg)^k=
\bigg(\sum_{i=1}^4 \lambda_i x_i\bigg)^k.
\end{equation}

By Equation \ref{term1}, Equation \ref{term2}, and
the case $M=4$ of Lemma A1, we have
{\footnotesize
\begin{equation}
  \label{boundX}
  \sum_{i=1}^4 \lambda_i(z_1)F(\|\widehat q_{1,i}-\widehat z_2\|) - X_{\bullet}=
  \sum_{i=1}^4 \lambda_i x_i^k - \bigg(\sum_{i=1}^4 \lambda_i x_i\bigg)^k
\leq
\frac{1}{8} k(k-1)x_+^{k-2}(x_+-x_-)^2.
\end{equation}
\/}
Let $\alpha$ and $\beta$ be such that
$x_{\alpha}=x_+$ and $x_{\beta}=x_-$.
By Lemma G2
\begin{equation}
\label{subs1}
x_{\alpha}=2+2(\widehat q_{1,\alpha} \cdot \widehat z_2) \leq T.
\end{equation}
Since $d_1$ is the diameter of  $Q_1^{\bullet}$, and
$\widehat z_2$ is a unit vector,

\begin{equation}
\label{subs2}
x_\alpha-x_\beta
=
2 \widehat z_2 \cdot
(\widehat q_{1,\alpha}-\widehat q_{1,\beta})
\leq
2\|\widehat q_{1,\alpha}-\widehat q_{1,\beta}\|
\leq 2d_1.
\end{equation}

Plugging Equations \ref{subs1} and
\ref{subs2} into Equation
\ref{boundX}, we get Lemma E11.
\endproof

\noindent
{\bf Proof of Lemma E12:\/}
Lemma G1 implies that $\|\widehat z_1 - z_1^{\bullet}\| \leq \delta(Q_1)=\delta_1$.
Let $\gamma_1$ denote the unit speed line segment connecting
$z_1^{\bullet}$ to $\widehat z_1$.  The length $L$ of
$\gamma_1$ is at most $\delta_1$, by Lemma G1.
So,
$\gamma_1(0)=z_1^{\bullet}$ and $\gamma_1(L)=\widehat z_1.$
Define
\begin{equation}
f(t)=
\bigg(2+2 \widehat z_2 \cdot \gamma_1(t)\bigg)^k \hskip 5 pt
{\rm or\/} \hskip 5 pt - 2 - 2 \widehat z_2 \cdot \gamma_1(t),
\end{equation}
depending on the case.  The argument we give works equally
well more generally when we use $F=\pm G_k$.

We have
$f(0)=X_{\bullet}$ and
$f(L)=F(\|\widehat z_1-\widehat z_2\|)$.
Hence 
\begin{equation}
\label{integ}
X_{\bullet}-F(\|\widehat z_1-\widehat z_2\|)=f(0)-f(L), \hskip 30 pt L \leq \delta_1.
\end{equation}

Since $z_1^\bullet$ lies in the convex hull of points of $S^2$,
we have $\|z_1^\bullet\|\leq 1$.  Hence every point
$\gamma_1(t)$ on the segment joining $z_1^\bullet$ to
$\widehat z_1$ also satisfies $\|\gamma_1(t)\|\leq 1$.
Therefore
$\widehat z_2\cdot \gamma_1(t)\geq -1$, which implies that
$2+2\widehat z_2\cdot\gamma_1(t)\geq 0$.

Combining the Chain Rule, the Cauchy-Schwarz inequality, and Lemma G2, we have
$$
|f'(t)|=
\bigg|\big(2 \widehat z_2 \cdot \gamma_1'(t)\big) \times
k\bigg(2+2\widehat z_2 \cdot \gamma_1(t)\bigg)^{k-1}\bigg| \leq $$
$$2k\bigg|\big(2+2\widehat z_2 \cdot \gamma_1(t)\big)\bigg|^{k-1} \leq
2k \big(2+2(Q_1 \cdot Q_2)\big)^{k-1}=2kT^{k-1}.
$$
In short
\begin{equation}
  \label{diff}
  |f'(t)| \leq 2kT^{k-1}.
  \end{equation}
Lemma E12 follows from Equation \ref{diff}, Equation \ref{integ}, and
integration.
\endproof

\newpage

%% file: 10calc.tex
\section{The Calculation Theorem}
\label{calcthm}

\noindent
{\bf Reading Guide:\/}
This chapter is for Reader 6.
We prove the Calculation Theorem
from \S \ref{UPS}

\subsection{More about Blocks}

Here we discuss blocks in more detail.
Blocks are defined in Equation \ref{BLOCK}, and
$r$-acceptable blocks are defined immediately thereafter.
\newline
\newline
{\bf Almost Acceptable Blocks:\/}
Let $B$ be a block.  Recall that $B$
is {\it acceptable\/} if $B$ is $3/2$-acceptable.
We say that
$B$ is {\it almost acceptable\/} if
$B$ is $2$-acceptable. 
Our algorithm distinguishes between acceptable
and almost acceptable blocks.
This is a carefully considered
design choice.
\newline
\newline
{\bf Dyadic Subdivision:\/}
We let $S_0(B)$ denote the two blocks obtained by bisecting
$Q_0$.  For indices $k=1,2,3$, we let $S_k(B)$ denote the four
blocks obtained by cutting $Q_k$ in half in both coordinate
directions.   We call a block {\it rational\/} if all the coordinates
of all the vertices are dyadic rational numbers.
We will start with the block $[0,2] \times ([-2,2]^2)^3$ and iteratively
subdivide it. This will produce only rational blocks.  In one
exceptional case we will start with $[2,4] \times ([-2,2]^2)^3$.
Iterated subdivision likewise produces rational blocks.
\newline
\newline
{\bf Vertex Norm Property:\/}
      Let $B=Q_0 \times Q_1 \times Q_2 \times Q_3$ be an almost acceptable block.
      Recall that $Q_i$ is axis-aligned and contained in a closed
      quadrant for $i=1,2,3$.
      For such squares, the minimum
      and maximum norm of a point in the square are realized at a vertex.
      Thus, we check that a square lies inside (respectively outside)
      a disk of radius $r$ centered at the origin by checking
      that the squared-norm of each vertex is at most (respectively
      at least) $r^2$.  This is really the reason why we have included
      the ``contained in a closed quadrant'' condition.

\subsection{The Four Calculation Ingredients}

      A {\it rational block computation\/} is
      a finite calculation, only involving the arithmetic operations
      and min and max. The output is either {\bf yes\/} or {\bf no\/}.
      Let $\Omega$ and $\Omega_0$ and $\Upsilon$ be as in the Calculation Theorem.
      \newline
      \newline
      {\bf Ingredient 1:\/}
        We describe a rational block computation $C_1$ such that an output of
        {\bf yes\/} for a block $B$ implies that $B \subset \Omega_0$.
        Define intervals

{\footnotesize
\begin{equation}
  I_0=[-2^{-17},2^{-17}], \hskip 10 pt
  I_1=[1-2^{-17},1+2^{-17}] \hskip 10 pt
  2^{30} I_{\sqrt 3/3}=[619916940,619933323]
\end{equation}
\/}
$I_{\sqrt 3/3}$ is a rational interval that is just barely contained
inside the interval of length $2^{-16}$ centered at $\sqrt 3/3$.
Define
\begin{equation}
  \Omega_{00} = (I_1 \times \{0\}) \times (I_0 \times -I_{\sqrt 3/3})  \times (-I_1 \times I_0)
  \times (I_0 \times I_{\sqrt 3/3}).
\end{equation}
We have $\Omega_{00} \subset \Omega_0$, though just barely.
      There are $128$ vertices of $B$.  We simply check whether each of these
      vertices is contained in $\Omega_{00}$.
      If so then we return {\bf yes\/}.  Otherwise we return {\bf
        no\/}.  We scale
      up all the coordinates by $2^{30}$ so that this test just involves
      integer comparisons.
      \newline
      \newline
      {\bf Ingredient 2:\/}
     We describe a rational block computation $C_2$ such that
    an output of {\bf yes\/} for an almost acceptable block $B$
    means that $B \cap \Omega=\emptyset$.
    We return {\bf no\/} if we don't return {\bf yes\/} in
    any
    of the calculations.
        \newline
        \newline
        {\bf A (avatar parity):\/}
        Let $a_{<}$ and $a_{>}$ respectively be the number of
        indices $j$ such that all vertices of $Q_j$ have norm
        less than or greater than $1/2$.
      We return {\bf yes\/} if $a_{<}$ is odd and $a_{>}+a_{<}=4$.  By the Vertex
      Norm Property
      every avatar in $B$ is odd and hence outside
      $\Omega$.
      \newline
      \newline
      {\bf B (norm ordering):\/}
      We return {\bf yes\/} if there is some $j \in \{1,2,3\}$ such that
      all vertices of $Q_j$ have norm greater than $\max Q_0$.  By the
      V.N.P.,
      we have $\|p_j\|>\|p_0\|$ for all avatars in $B$.
      Such avatars are not in $\Omega$.
      \newline
      \newline
     {\bf C (max norm bound):\/}
      We return {\bf yes\/} if there is some $j \in \{1,2,3\}$ such
      that all vertices of $Q_j$ have norm at least $3/2$.
      By the V.N.P.,
      every avatar $(p_0,p_1,p_2,p_3)$ in
      $B$  has $\|p_j\| \geq 3/2$ and hence lies outside
      $\Omega$. 
      \newline
      \newline
{\bf D (height ordering):\/}
      We write $I \leq J$ to indicate that all values in an interval
      or singleton $I$ are less than or equal to all values in an
      interval
      or singleton $J$.
      For each $j=0,1,2,3$ we let
      $YQ_j$ be the projection of $Q_j$ onto the  $Y$-axis.
      We return {\bf yes\/} if
      $YQ_1> YQ_2$ or $YQ_1> YQ_3$ or $YQ_2
        > YQ_3$ or $YQ_{2}< 0$.
      This condition guarantees that all avatars in 
      $B$ lie outside $\Omega$.
        \newline
      \newline
      {\bf Ingredient 3:\/}
        We describe a rational block computation $C_3^{\sharp}$ such that an output of
        {\bf yes\/} for a block $B$ implies that $B \cap \Omega \subset \Upsilon$.
        Likewise, there exists a rational block computation $C_3^{\sharp \sharp}$ such that an output of
        {\bf yes\/} for a block $B$ implies that $B$ is disjoint from
        $\Upsilon$.
        \newline
        Referring to \S \ref{special} we use the symmetrized domain
        $\Upsilon^*$ in place of $\Upsilon$ for $C_3^{\sharp}$.
      For $C_3^\sharp$ we return {\bf yes\/} if all the vertices of $B$
     lie in the rectangular product specified in the definition of
     $\Upsilon^*$.   We are ruling out blocks $B$ that lie in
     $\Upsilon^*$.
      Since we only consider configurations in $\Omega$, this implies that
      $B\cap\Omega\subset\Upsilon \cap \Omega$.  See
      Equation \ref{symmU}.   We discuss why we use $\Upsilon^*$
      instead of $\Upsilon$ at the end of \S \ref{discussXX}.

For $C_3^{\sharp \sharp}$ we return {\bf yes\/} if one of the factors
of $B$ is disjoint from the corresponding factor of this rectangular
product.
These are exact integer tests involving scaled-up versions of the
vertex coordinates.
See the discussion of implementation below.
            \newline
      \newline
      {\bf Ingredient 4:\/}
       For any function $F$ given by Equation \ref{functionlist}, 
        we describe a rational block computation $C_{4,F}$ such that an output of
        {\bf yes\/} for an almost acceptable block $B$ implies that the minimum of ${\cal E\/}_F$ on
        $B$ is at least ${\cal E\/}_F(\xi_0)+2^{-50}$.
                \newline

                First, we return {\bf no\/} for any almost acceptable block that is not
                acceptable.   We also return the smallest $j$ such
                that $Q_j \not \subset [-3/2,3/2]^2$.
                This is a short calculation done with
                exact integers -- scaled up versions of
                the coordinates.  Thus we only
                apply the Energy Theorem (\S \ref{ENG})  to acceptable blocks.
                \newline

                At this point we have an acceptable block $B$.
Let $F$ be an
energy hybrid.
Let $[F]$ denote the $F$-potential of the TBP.
If
\begin{equation}
  \min_{p \in v(B)} {\cal E\/}_F(p)-{\bf ERR\/}_F(B) \geq [F]+2^{-50}
\end{equation}
we return {\bf yes\/}. Otherwise we return {\bf no\/}. We also
return the index
$i$ such that ${\bf ERR\/}_F(B,i)$ is the largest.
In case of a tie, which probably never happens,
we pick the highest such index.

\subsection{The Algorithm}

      Here is the main calculation.
      \begin{enumerate}
        
      \item We start with the list $L=\{[0,2] \times ([-2,2]^2)^3\}$.
        This is a single block.
        
      \item If $L=\emptyset$ then {\bf HALT\/}.  Otherwise let $B=Q_0 \times Q_1 \times Q_2 \times Q_3$
        be the last block of $L$.
        
      \item If $Q_0$ has length greater than $1$, or
        any $Q_i$, for $i=1,2,3$, has side length greater than
           $1$ we delete $B$ from $L$ and
        append to $L$ the blocks of the subdivision $S_k(B)$ where $k$
        is the lowest index of a component having the offending
        property.  Then we go to Step 2.

      \item At this step, and for later steps,
        $B$ is almost acceptable because its
      square components lie in $[-2,2]^2$, and the subdivision
    process has forced them to lie in closed quadrants.
    If $C_1(B)={\bf yes\/}$ we remove $B$ from $L$ and go to Step 2.
    We are eliminating blocks in $\Omega_{00}$.

    \item If  $C_2(B)={\bf yes\/}$ we remove $B$ from $L$ and go to
      Step 2.
      Here we are eliminating $B$ when
      $B \cap \Omega=\emptyset$.
    
  \item If $F=G_{10}^{\sharp}$ and $C_3^{\sharp}(B)={\bf yes\/}$ we remove $B$ from $L$ and
        go to Step 2.
        If $F=G_{10}^{\sharp \sharp}$ and $C_3^{\sharp \sharp}(B)={\bf yes\/}$ we remove $B$ from $L$ and
         go to Step 2.

      \item  If $C_{4,F}(B)={\bf yes\/}$ then we remove $B$ from $L$ and go
        to Step 2.  Here we have verified that the $F$-energy of any avatar in $B$
        is at least $[F]+2^{-50}$.
        
      \item If $C_{4,F}(B)={\bf no\/}$ and the integer is
        $k$, then we delete $B$ from $L$ and
        append to $L$ the blocks of the subdivision $S_k(B)$ and
        return to step 2.
        
      \end{enumerate}

      \noindent
      {\bf Exceptional Run:\/}
      Referring to \S \ref{POLYA}, we make two runs for $G_3$, one
      starting with
      $[0,2] \times ([-2,2]^2)^3$ and one starting with
      $[2,4] \times
      ([-2,2]^2)^3$.

       \subsection{Discussion of the Implementation}
      \label{discussXX}

      The directory {\bf 6BigCalculation/Main\/} contains the Java program used
      in the big calculation.
      Now we give more
      information about the calculation.
      \newline
      \newline
\noindent
{\bf Representing Blocks:\/}
We represent the coordinates of blocks by the {\it long\/} data type
scaled up by $2^{30}$.  This makes all (scaled) coordinates integers.
      Our algorithm never does so many
      subdivisions that it defeats this method of
      representation.
      In all but the main step in the algorithm above we compute with
      exact integers.  When the calculation (such as squaring a
      {\tt long\/}) could cause an overflow error, we first recast
      the {\tt longs\/} as a {\tt BigIntegers\/} in Java and then do
      the calculations.
      \newline
      \newline
      {\bf Interval Arithmetic:\/}
      For the main step of the algorithm we
      use interval arithmetic.  We use the same implementation
      as we did in [{\bf S1\/}], where we explain it in detail.
      Here is how it works in brief. If we have a
      calculation involving numbers $r_1,...,r_n$, then we produce
      intervals $I_1,...,I_n$ with dyadic rational numbers represented
      exactly by the computer such that $r_i \in I_i$ for $i=1,...,n$.
      We then perform the usual arithmetic operations on the intervals,
      rounding outward at each step.  The final output of the calculation,
      an interval, contains the result of the actual calculation.
      \newline
      \newline
          {\bf Guess and Check:\/}
          Here is how we speed up the calculation.
      When we do Steps 7 -- 8, we
      first do the calculation $C_{4,F}$ using floating point operations.
      If the floating version returns an integer, we use this integer to subdivide
      the box and return to step 2.  If $C_{4,F}$ says {\bf yes\/} then
      we retest the box using the interval arithmetic.  In this way, we
      only pass a box for which the interval version says {\bf yes\/}.
      This way of doing things keeps the calculation rigorous but speeds
      it up by using the interval arithmetic as sparingly as
      possible.
      \newline
      \newline
      {\bf Weak versus Strict:\/}
      In Ingredient 2 we mostly use strict inequalities.  For
      instance, Test {\bf D\/}  says {\bf yes\/} to a block $B$ if all
      avatars $\xi =(p_0,p_1,p_2,p_3)$ satisfy
      $p_{22}<0$.   We could implement a weaker test where
      we say {\bf yes\/}  if $p_{22} \leq 0$ for all $\xi \in B$.
      Using weak inequalities like this speeds up the calculations
      by a factor of about $3/2$, and for a huge calculation like
      the $G_{10}^{\sharp}$ calculation, this is very convenient.
      However, the weaker tests create edge cases, avatars
      lying in the boundaries of blocks which pass Step 5 of
      the algorithm
        yet which  belong to (the boundary of) $\Omega$.  It turns out
        that these edge cases are included in other blocks that
        are passed at Step 7 and thus are not actually
      minimizers.  However, the argument which proves this
      is rather tedious.  I opted for a longer calculation
      so that I could omit this tedious argument.
      \newline
      \newline
      {\bf Symmetrized Domain:\/}
      In Ingredient 3,
      we use $\Upsilon^*$ in place of
      $\Upsilon$ because our strict Ingredient 2 tests
      do not eliminate blocks having a
      square component lying below the $X$-axis and having
      its top edge on the $X$-axis on the boundary of the
      $p_2$-component of $\Upsilon$, the leftmost rectangle in
      Figure 2.2.      The corresponding avatars are potentially isomorphic to
      ones
      in $\Upsilon$.  They are irrelevant to the calculation
      and need to be eliminated somehow.  This is where the choice
      of $\Upsilon^*$ does its job.  Were we to use the weak
      versions of the Ingredient Tests, our algorithm would run
      with $\Upsilon$ in place of $\Upsilon^*$.

      \subsection{Record of the Calculation}

      If the algorithm reaches the {\bf HALT\/} state for a given choice of $F$, this
      constitutes a proof that the corresponding statement of the
      Calculation Theorem is true.  In fact this happens in all cases.
      Here I give an account of one time I ran the
      computations to completion in Sept 2026 using the
            computer
      discussed at the end of the introduction.
      In listing the calculations I will give some of the
      calculation statistics:
     M.L. stands for the maximum list length,
     M.D. stands for the maximum subdivision depth, and
P.S.  stands for partition size.  The first 2 lines are relevant for
the Theorem \ref{EXP} and the remaining lines are
relevant for the Phase Transition Theorem.

\begin{center}
\begin{tabular}{c|c|r|r|r|r}
Potential & Domain & Running Time & P.S. & M.L. & M.D.\\
\hline
$G_3$ & $\Omega-\Omega_0$ & 0 hr 35 min 07 sec & 16646301 & 154 & 22\\
$G_3$ & $\Omega^*-\Omega$ & 0 hr 02 min 05 sec & 1132513 & 36 & 6\\
$G_4$ & $\Omega-\Omega_0$ & 0 hr 41 min 22 sec & 19430121 & 155 & 22\\
$G_5^{\flat}$ & $\Omega-\Omega_0$ & 3 hr 04 min 56 sec & 56673218 & 164 & 23\\
$G_6$ & $\Omega-\Omega_0$ & 1 hr 50 min 47 sec & 47269653 & 165 & 23\\
$G_{10}^{\sharp}$ & $\Omega-\Omega_0-\Upsilon$ & 19 hr 45 min 10 sec & 193528137 & 169 & 24\\
$G_{10}^{\sharp\sharp}$ & $\Upsilon-\Omega_0$ & 0 hr 08 min 53 sec & 1619236 & 66 & 9
\end{tabular}
\end{center}

\newpage

%% file: refs.tex
\section{References}

\noindent
[{\bf A\/}] N. N. Andreev,
{\it An Extremal Property of the Icosahedron\/},
East J. Approx. {\bf 2\/} (1996), no. 4, pp. 459--462.
\vskip 8 pt

\noindent
[{\bf BBCGKS\/}] Brandon Ballinger, Grigoriy Blekherman, Henry Cohn,
Noah Giansiracusa, Elizabeth Kelly, Achill Sch\"urmann,
{\it Experimental Study of Energy-Minimizing Point Configurations on Spheres\/},
Experimental Mathematics {\bf 18\/} (2009), no. 3, pp. 257--283.
\vskip 8 pt

\noindent
[{\bf BDHSS\/}] P. G. Boyvalenkov, P. D. Dragnev, D. P. Hardin,
E. B. Saff, M. M. Stoyanova,
{\it Universal Lower Bounds for Potential Energy of Spherical Codes\/},
Constructive Approximation {\bf 44\/} (2016), pp. 385--415.
\vskip 8 pt

\noindent
[{\bf BHS\/}] A. V. Bondarenko, D. P. Hardin, E. B. Saff,
{\it Mesh Ratios for Best Packings and Limits of Minimal Energy Configurations\/},
Acta Math. Hungar. {\bf 142\/} (2014), no. 1, pp. 118--131.
\vskip 8 pt

\noindent
[{\bf C\/}] Harvey Cohn,
{\it Stability Configurations of Electrons on a Sphere\/},
Mathematical Tables and Other Aids to Computation
{\bf 10\/} (1956), no. 55, pp. 117--120.
\vskip 8 pt

\noindent
[{\bf CK\/}] Henry Cohn and Abhinav Kumar,
{\it Universally Optimal Distributions of Points on Spheres\/},
J.A.M.S. {\bf 20\/} (2007), pp. 99--148.
\vskip 8 pt

\noindent
[{\bf CCD\/}] The Riesz Energy Database, online website:
\newline
https://www-wales.ch.cam.ac.uk/$\sim$wales/CCD/Thomson/site/index.html
\vskip 8 pt

\noindent
[{\bf DLT\/}] P. D. Dragnev, D. A. Legg, and D. W. Townsend,
{\it Discrete Logarithmic Energy on the Sphere\/},
Pacific Journal of Mathematics {\bf 207\/} (2002), no. 2,
pp. 345--357.
\vskip 8 pt

\noindent
[{\bf F\"o\/}] L. F\"oppl,
{\it Stabile Anordnungen von Elektronen im Atom\/},
J. f\"ur die Reine Angew. Math. {\bf 141\/} (1912),
pp. 251--302.
\vskip 8 pt

\noindent
[{\bf HZ\/}] Xiaorong Hou and Junwei Shao,
{\it Spherical Distribution of 5 Points with Maximal Distance Sum\/},
Discrete Comput. Geom. {\bf 46\/} (2011), no. 1,
pp. 156--174.
\vskip 8 pt

\noindent
[{\bf I\/}] IEEE Standard for Binary Floating-Point Arithmetic
(IEEE Std 754-1985),
Institute of Electrical and Electronics Engineers, July 26, 1985.
\vskip 8 pt

\noindent
[{\bf KY\/}] A. V. Kolushov and V. A. Yudin,
{\it Extremal Dispositions of Points on the Sphere\/},
Anal. Math. {\bf 23\/} (1997), pp. 143--146.
\vskip 8 pt

\noindent
[{\bf MKS\/}] T. W. Melnyk, O. Knop, W. R. Smith,
{\it Extremal Arrangements of Points and Unit Charges on the Sphere:
Equilibrium Configurations Revisited\/},
Canadian Journal of Chemistry {\bf 55\/} (1977),
pp. 1745--1761.
\vskip 8 pt

\noindent
[{\bf RSZ\/}] E. A. Rakhmanov, E. B. Saff, and Y. M. Zhou,
{\it Electrons on the Sphere\/},
in {\it Computational Methods and Function Theory\/},
R. M. Ali, St. Ruscheweyh, and E. B. Saff, Eds. (1995),
pp. 111--127.
\vskip 8 pt

\noindent
[{\bf S0\/}] R. E. Schwartz,
{\it Divide and Conquer: A Distributed Approach to $5$-Point Energy Minimization\/},
Research Monograph (preprint, 2023).
\vskip 8 pt

\noindent
[{\bf S1\/}] R. E. Schwartz,
{\it The Five-Electron Case of Thomson's Problem\/},
Experimental Mathematics {\bf 22\/} (2013), no. 2,
pp. 157--186.
\vskip 8 pt

\noindent
[{\bf S2\/}] R. E. Schwartz,
{\it The Projective Heat Map\/},
Mathematical Surveys and Monographs {\bf 219\/},
American Mathematical Society, 2017.
\vskip 8 pt

\noindent
[{\bf S3\/}] R. E. Schwartz,
{\it Lengthening a Tetrahedron\/},
Geometriae Dedicata {\bf 174\/} (2014),
pp. 121--144.
\vskip 8 pt

\noindent
[{\bf S4\/}] R. E. Schwartz,
{\it Five Point Energy Minimization: A Synopsis\/},
Constructive Approximation {\bf 51\/} (2020), no. 3,
pp. 537--564.
\vskip 8 pt

\noindent
[{\bf SK\/}] E. B. Saff and A. B. J. Kuijlaars,
{\it Distributing Many Points on a Sphere\/},
Math. Intelligencer {\bf 19\/} (1997), no. 1,
pp. 5--11.
\vskip 8 pt

\noindent
[{\bf Th\/}] J. J. Thomson,
{\it On the Structure of the Atom: an Investigation of the
Stability and Periods of Oscillation of a Number of Corpuscles Arranged
at Equal Intervals around the Circumference of a Circle with Application
of the Results to the Theory of Atomic Structure\/},
Philosophical Magazine, Series 6, {\bf 7\/} (1904), no. 39,
pp. 237--265.
\vskip 8 pt

\noindent
[{\bf T\/}] A. Tumanov,
{\it Minimal Biquadratic Energy of $5$ Particles on $2$-Sphere\/},
Indiana Univ. Math. Journal {\bf 62\/} (2013),
pp. 1717--1731.
\vskip 8 pt

\noindent
[{\bf W\/}] S. Wolfram,
{\it The Mathematica Book\/}, 4th ed.,
Wolfram Media/Cambridge University Press,
Champaign/Cambridge (1999).
\vskip 8 pt

\noindent
[{\bf Y\/}] V. A. Yudin,
{\it Minimum Potential Energy of a Point System of Charges\/}
(Russian), Diskret. Mat. {\bf